\documentclass{article}
\usepackage[utf8]{inputenc}

\newcommand{\change}[1]{{#1}} 
\newcommand{\newchange}[1]{{#1}} 
\newcommand \conv {\mathrm{conv}}

\title{When Deep Learning Meets Polyhedral Theory: A Survey}
\author{
Joey Huchette\\{\footnotesize Google Research, USA} \and
Gonzalo Mu\~{n}oz\\{\footnotesize ~~Universidad de Chile, Chile~~} \and 
Thiago Serra\\{\footnotesize University of Iowa, USA} \and
Calvin Tsay\\{\footnotesize Imperial College London, UK}
}
\date{September 2025}

\usepackage[comma]{natbib}
\newcommand{\samplesize}{D}
\newcommand{\sampleinput}{\tilde{\vx}}
\newcommand{\sampleoutput}{\tilde{\vy}}
\newcommand{\functolearn}{\hat{f}}
\newcommand{\weights}{\mW}
\newcommand{\biases}{\vb}
\newcommand{\hiddenlayers}{L}
\newcommand{\inputdimension}{{n_0}}
\newcommand{\outputdimension}{{n_{\hiddenlayers+1}}}
\newcommand{\layerwidth}{n}
\newcommand{\maxlayerwidth}{{n_{\max}}}
\newcommand{\parameterset}{\Theta}
\newcommand{\parameternumber}{N}

\newcommand{\NNfunction}[3]{f(#1,#2,#3)}
\newcommand{\loss}{\ell}
\newcommand{\ERMfunction}{\mathcal{L}}
\newcommand{\setofones}{Y}

\usepackage{amsmath,amsfonts,bm}

\def\1{\bm{1}}

\def\vb{{\bm{b}}}

\def\vd{{\bm{d}}}

\def\vh{{\bm{h}}}

\def\vt{{\bm{t}}}
\def\vu{{\bm{u}}}
\def\vv{{\bm{v}}}
\def\vw{{\bm{w}}}
\def\vx{{\bm{x}}}
\def\vy{{\bm{y}}}
\def\vz{{\bm{z}}}

\def\evh{{h}}

\def\mA{{\bm{A}}}

\def\mI{{\bm{I}}}

\def\mL{{\bm{L}}}

\def\mT{{\bm{T}}}

\def\mV{{\bm{V}}}
\def\mW{{\bm{W}}}
\def\mX{{\bm{X}}}

\DeclareMathAlphabet{\mathsfit}{\encodingdefault}{\sfdefault}{m}{sl}
\SetMathAlphabet{\mathsfit}{bold}{\encodingdefault}{\sfdefault}{bx}{n}

\def\sI{{\mathbb{I}}}

\def\sL{{\mathbb{L}}}

\def\sS{{\mathbb{S}}}

\newcommand{\R}{\mathbb{R}}

\usepackage{braket}
\newcommand{\gr}{\texttt{gr}}

\usepackage{stmaryrd}
\usepackage{mathtools}

\usepackage{array}
\usepackage{multirow}

\usepackage{amssymb}

\usepackage{todonotes}
\newtheorem{theorem}{Theorem}
\newtheorem{definition}{Definition}
\newtheorem{example}{Example}
\newtheorem{conjecture}{Conjecture}

\usepackage{lscape}

\usepackage{adjustbox}
\usepackage{url}

\usepackage{tikz-network}
\begin{document}

\maketitle

\begin{abstract}
\noindent 
In the past decade, deep learning became the prevalent methodology for predictive modeling thanks to the remarkable accuracy of deep neural networks in tasks such as computer vision and natural language processing. 
Meanwhile, the structure of neural networks converged back to simpler representations based on piecewise constant and piecewise linear functions such as the Rectified Linear Unit~(ReLU), 
which became the most commonly used type of activation function in neural networks. 
That made certain types of network structure ---such as the typical fully-connected feedforward neural network--- amenable to analysis 
through polyhedral theory and to the application of methodologies such as Linear Programming~(LP) and Mixed-Integer Linear Programming~(MILP) for a variety of purposes. 
In this paper, 
we survey the main topics emerging from this fast-paced area of work, 
which brings a fresh perspective to understanding neural networks in more detail as well as to applying linear optimization techniques to train, verify, and reduce the size of such networks.  
\end{abstract}

\section{Introduction}
\label{sec:intro}
\begin{figure}
\centering
\includegraphics[width=0.5\textwidth]{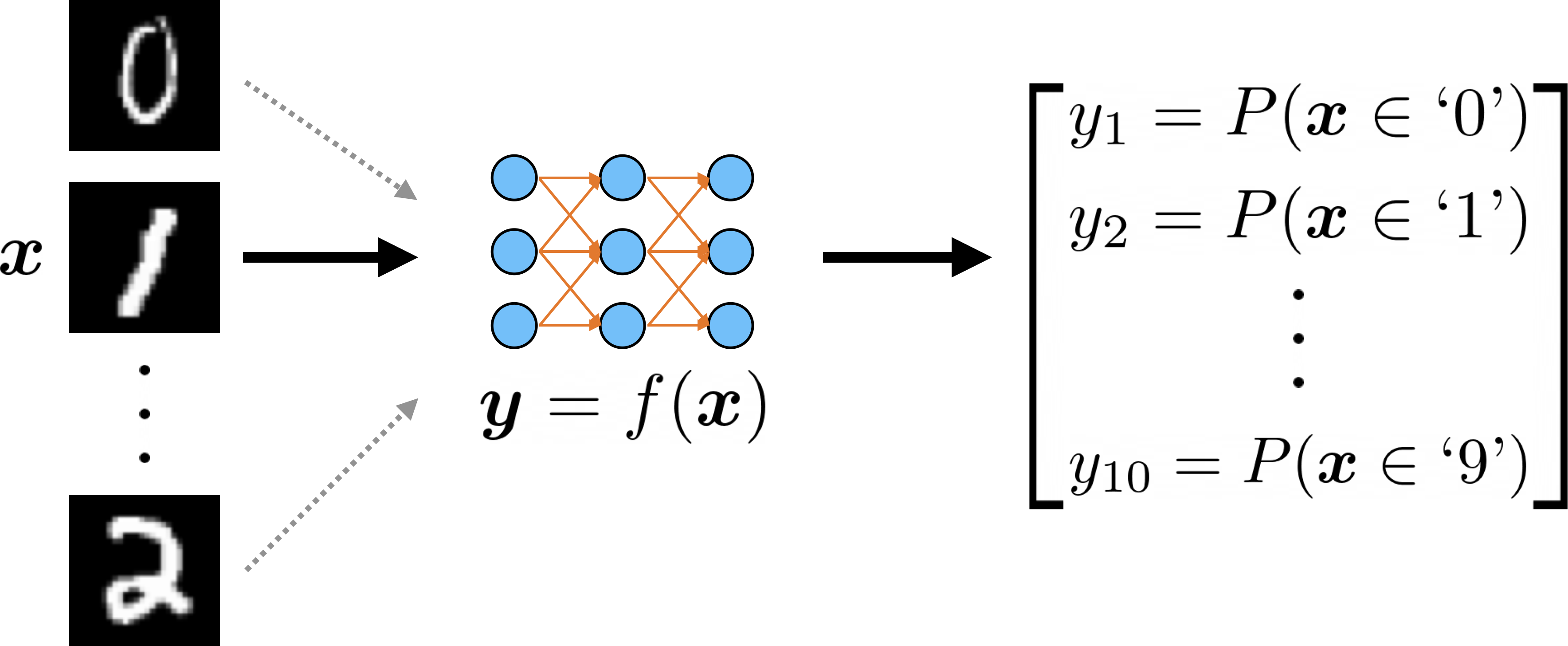}
\caption{Example classification task on the MNIST database of handwritten digits, in which the image of a handwritten digit is given as input and the probability of that digit being from each possible class is provided as output.}
\label{fig:mnist}
\end{figure}

Deep learning has continuously achieved new landmarks in varied areas of artificial intelligence for the past decade. Examples of those areas include predictive tasks in computer vision \citep{Krizhevsky2012,Ciresan2012,Szegedy2015,He2016DeepRL,NoisyStudent}, natural language processing \citep{sutskever2014sequence,ELMo,GPT,BERT}, and speech recognition \citep{Hinton2012,EndToEndSpeech,SpecAugment}. 
The artificial neural networks behind such feats are being used in many applications, 
and there is a growing interest for analytical insights to 
help design such networks and then to leverage the model that they have learned. 
For the most commonly used types of neural networks, some of those results and methods are coming from operations research tools such as polyhedral theory and associated optimization techniques such as Linear Programming~(LP) and Mixed-Integer Linear Programming~(MILP). 
Among other things, 
these connections with mathematical optimization may help us understand what neural networks can represent, how to train them, and 
how to make them more compact.
For example, consider the popular task of classifying images (Figure \ref{fig:mnist}); polyhedral theory and associated optimization techniques may help us answer questions such as the following. How should we train the classifier model? How large should it be? How robust to perturbations is it?

\subsection{What neural networks can model}

We can essentially think of artificial neural networks as functions mapping an input $\vx$ from a given domain to an output $\vy$ for a given application. 
For the classification task in Figure \ref{fig:mnist}, inputs $\vx$ correspond to images from the dataset, and $\vy$ to the associated predicted labels, or probabilities for labels describing the content of those images. 
The basic units of neural networks mimic biological neurons in that they receive inputs from adjacent units, transform those inputs, and may produce an output to subsequent units of the network. 
In other words, every unit is also a function, and in fact 
the output of most units is defined by the composition of a nonlinear function with a linear function. 
The nonlinear function is often denoted as the \emph{activation function} in analogy to how a biological neuron is triggered to send a signal to adjacent neurons when the stimulus caused by the input exceeds a certain activation threshold. Such non-linearity is behind the remarkable expressiveness of neural networks.

This model was pioneered by~\cite{FirstANN}, who considered a thresholding function for activation that is now often denoted as the Linear Threshold Unit~(LTU). 
That activation is also the basis of the classic \emph{perceptron} algorithm by ~\cite{perceptron}, 
which yields a binary classifier of the form 
\begin{equation}
f(\vx) = \left\{ \begin{array}{cl} 1 & \text{if } \vw \cdot \vx + b > 0; \\ 0 & \text{otherwise} \end{array} \right. 
\end{equation}
for an input $\vx \in \mathbb{R}^{n_0}$ and with parameters $\vw \in\mathbb{R}^{n_0}$ and $b \in \mathbb{R}$. Those parameters are chosen by optimizing the predictions for a given task, as discussed below and in Section \ref{sec:training}. The term \emph{single-layer perceptron} is used for a neural network consisting of a set of such units processing the same input in parallel. 
The term \emph{multi-layer perceptron} is used for a generalization of this concept, by which the output of a \emph{layer}---a set of units with the same input---is the input for a subsequent layer. This perceptron terminology has also been loosely applied to neural networks with other activation functions.

\begin{figure}
    \centering
    \includegraphics[width=\textwidth]{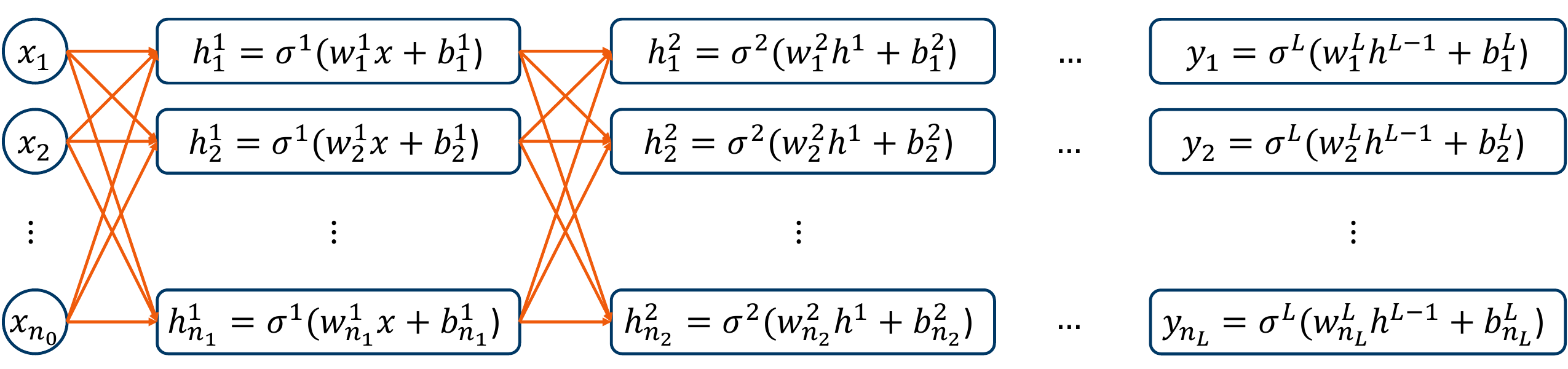}
    \caption{Mapping from $\vx \in \mathbb{R}^{n_0}$ to $\vy \in \mathbb{R}^{n_L}$ through a feedforward neural network with $L$ layers, layer widths $\{ n_l \}_{l \in \sL}$, and activation functions $\{ \sigma_l \}_{l \in \sL}$.}
    \label{fig:ins_outs}
\end{figure}

More generally, neural networks that successively transform inputs through an ordered sequence of layers are also denoted \emph{feedforward networks}. 
The layers that do not produce the final output of the neural network are denoted \emph{hidden layers}. 
For a network with $L$ layers, 
we denote $n_l$ as the number of units---or \emph{width}---of layer $l \in \sL := \{1, 2, \ldots, L\}$ and $h_i^l$ as the output of the $i$-th unit in layer $l$, where $i \in \{1, 2, \ldots, n_l\}$. The output of a unit is given by
\begin{equation} \label{eq:single-neuron}
h_i^l = \sigma^l\left(\vw^l_i \cdot \vh^{l-1} + b^l_i\right),
\end{equation}
where the \emph{weights} $\vw^l_i \in \mathbb{R}^{n_{l-1}}$ and the \emph{bias} $b^l_i \in \mathbb{R}$ are parameters of the unit. 
Those parameters can be aggregated across the layer as the matrix $\mW^l \in \mathbb{R}^{n_l \times n_{l-1}}$ and the vector $\vb^l \in \mathbb{R}^{n_l}$. 
The vector $\vh^{l-1} \in \mathbb{R}^{n_{l-1}}$ represents the aggregated outputs from layer $(l-1)$. 
The activation function $\sigma^l : \mathbb{R} \rightarrow \mathbb{R}$ is applied by the units in layer $l$. 
These definitions implicitly assume that $n_0$ is the size of the network input $\vx \in \mathbb{R}^{n_0}$ and that $\vh^0$ and $\vx$ are the same. 
Figure~\ref{fig:ins_outs} illustrates the operation of a feedforward network as described above. 

\subsection{How neural networks are obtained and evaluated}
In resemblance to how other models for \emph{supervised learning} in machine learning are obtained,  
we can \emph{train} a neural network for a given task by adjusting its behavior with respect to the examples of a \emph{training set} and then evaluate the final trained network on a \emph{test set}. 
Both of these sets consist of inputs for which the correct output $\hat{y}$ is known. 
We can define an objective function to model 
a measure of distance 
between the output $y$ and the correct output $\hat{y}$, which is typically denoted as the \emph{loss function}, and then iteratively update parameters such as $\{ \mW^l \}_{l \in \sL}$ and $\{ \vb^l \}_{l \in \sL}$ to minimize that loss function over the training set. A common objective function is the square error $\|y - \hat{y}\|^2$ summed over the points in the training set. 
The test set contains a separate collection of inputs and their outputs, 
which is used to evaluate the trained neural network with examples that were not seen during training. 
A good performance on the test set may indicate that the trained neural network is able to \emph{generalize} beyond the seen examples, whereas a bad performance may suggest that it \emph{overfits} for the training set. 
Neural networks also have \emph{hyperparameters} that are often chosen manually and do not change during training, such as the \emph{depth} $L$, the widths of the layers $\{ n_l \}_{l \in \sL}$, and the activation functions used in each layer $\{ \sigma^l \}_{l \in \sL}$. 
Different models can be produced by varying the hyperparameters. 
In such a case, a \emph{validation set} disjoint from the training and test sets can be used to compare models with different hyperparameters. 
Whereas the validation set may serve as a benchmark to compare different trained models corresponding to different choices of hyperparameters, 
the test set can only be used to evaluate a single neural network chosen among those evaluated with the validation set. 
The emergent field of \emph{neural architecture search}---recently surveyed by \cite{NAS}---concerns with automatically choosing such hyperparameters. 

One of the key factors for the success of deep learning is that first-order methods for continuous optimization can be effectively applied to train deep networks. 
The interest in neural networks first vanished due to negative results in the Perceptrons book by~\cite{perceptrons}, which showed that single-layer perceptrons cannot represent functions such as the Boolean XOR. 
However, moving to multi-layer perceptrons capable of expressing the Boolean XOR as well as other more expressive models would require a clever training strategy. 
Hence, the interest was regained with papers that popularized the use of \emph{backpropagation}, such as ~\cite{BackPOP1} and ~\cite{BackPOP2}. 
Note that backpropagation was first discussed much earlier in the context of networks by~\cite{BackpropNetworks} and of neural networks explicitly by~\cite{BackpropNN}. 
The backpropagation algorithm calculates the derivative of the loss function with respect to each neural network parameter by applying the chain rule through the units of the neural network, 
which is considerably more efficient than explicitly evaluating the derivative of each network parameter.
Consequently, neural networks are generally trained with gradient descent methods in which the parameters are updated sequentially from the output to the input layer in each step. 
In fact, 
most algorithms for training neural networks are based on Stochastic Gradient Descent~(SGD), 
which is a form of the stochastic approximation through sampling pioneered by~\cite{StochasticApproximation}. 
SGD approximates the partial derivatives of the loss function at each step by using only a subset of the data in order to make the training process more efficient. 
Examples of popular SGD algorithms include momentum \citep{momentum}, Adam \citep{Adam}, and Nesterov Adaptive Gradient \citep{NAG}---the latter inspired by~\cite{Nesterov}. 
Interestingly, 
however, we generally cannot guarantee convergence to a global optimum with gradient descent due to the nonconvexity of the loss function. 
Nevertheless, neural networks trained with adequately parameterized SGD algorithms tend to generalize well. 

\subsection{Why nonlinearity is important in artificial neurons}

{\color{black}The nonlinearity of the activation functions, the resulting layered structure, and/or the loss function yield nonconvex training problems. Likewise, the neural network itself is typically a nonconvex function even after training.}
However, as we will see in Section~\ref{sec:LR}, that same nonlinearity enables the neural network to represent more complex functions as a whole. 
In fact, removing such nonlinearities by using an identity activation function $\sigma^l(u) = u ~ \forall l \in \sL$ would reduce the entire neural network to an affine transformation of the form $f(x) = \mW^{L} ( \mW^{L-1} \left( \ldots \left( \mW^2 \left( \mW^1 x + \vb^1 \right) + \vb^2 \right) + \ldots \right) + \vb^{L-1} ) + \vb^{L}$. 
Hence, a feedforward network without nonlinear activation functions is a {\color{black}(rank-constrained) linear model}.
In contrast, 
neural networks with a single hidden layer of arbitrary width have been long known to be universal function approximators for a broad variety of activation functions \citep{Cybenko1989,funahashi1989approximate,Hornik1989}, as well as for ReLU more recently \citep{yarotsky2017relu}. 
These results have also been extended to the converse case of limited width but arbitrarily large depth \citep{lu2017expressive,hanin2017approximating,park2021width}. 

Although nonlinear activation functions are important for obtaining more complex models, these functions do not need to be overly complex to produce good results. 
In the past, it was common practice to use sigmoid functions for activation \citep{TricksBP}. Those are monotonically increasing functions that approach finite values for arbitrarily large positive and negative inputs, 
such as the standard logistic function $\sigma(u) = \frac{1}{1+e^{-u}}$ and the hyperbolic tangent $\sigma(u) = \tanh(u)$. 
In the present, the most commonly used activation function is the Rectified Linear Unit~(ReLU) $\sigma(u) = \max\{0,u\}$ \citep{CurrentDNN,ReluPop2018}, which was proposed by~\cite{OriginReLU} and first applied to neural networks by~\cite{nair2010rectified}. 
The popularity of ReLU is in part due to experiments by~\cite{nair2010rectified} and~\cite{ReLUGood2} showing that this simpler form of activation yields competitive results. 
Thinking back in terms of the analogy with biological neurons, 
we say that a ReLU is \emph{active} when the output is positive and \emph{inactive} when the output is zero. 
ReLUs have a linear output behavior on the inputs associated with the same ReLUs being active and inactive; this property also holds for other piecewise linear and piecewise constant functions that are used as activation functions in neural networks. 
Table~\ref{tab:activations} lists some of the most commonly used activation functions of that kind. For more comprehensive lists of activation functions, including several other variations based on ReLU, we refer to~\cite{dubey2021activations} and~\cite{tao2022piecewise}.

\begin{table}
\caption{Main piecewise constant and piecewise linear activation functions.}
\label{tab:activations}
\centering
\vspace{2ex}
\begin{tabular}{@{\extracolsep{3pt}}m{0.13\textwidth}cm{0.23\textwidth}}
\textbf{Name} & \textbf{Function} & \textbf{Reference} \\
\cline{1-1}
\cline{2-2}
\cline{3-3}
\noalign{\vskip8pt}
LTU & $\sigma(u) = \left\{ \begin{array}{cl} 1 & \text{if } u > 0 \\ 0 & \text{if } u \leq 0 \end{array}\right.$ & ~\cite{FirstANN} \\
\cline{1-1}
\cline{2-2}
\cline{3-3}
\noalign{\vskip8pt}
ReLU & $\sigma(u) = \max\{0,u\}$ &  ~\cite{OriginReLU,nair2010rectified} \\
\cline{1-1}
\cline{2-2}
\cline{3-3}
\noalign{\vskip8pt}
leaky ReLU & $\begin{array}{c} \sigma(u) = \left\{ \begin{array}{cl} u & \text{if } u > 0 \\ \varepsilon u & \text{if } u \leq 0 \end{array}\right.\\ \text{($\varepsilon$ is small and fixed)} \end{array}$ & ~\cite{leaky} \\
\cline{1-1}
\cline{2-2}
\cline{3-3}
\noalign{\vskip8pt}
parametric ReLU & $\begin{array}{c} \sigma(u) = \left\{ \begin{array}{cl} u & \text{if } u > 0 \\ a u & \text{if } u \leq 0 \end{array}\right. \\ \text{($a$ is a trainable parameter)} \end{array}$ & ~\cite{prelu} \\
\cline{1-1}
\cline{2-2}
\cline{3-3}
\noalign{\vskip8pt}
hard tanh & $\sigma(u) = \left\{ \begin{array}{cl} 1 & \text{if } u > 1 \\ u & \text{if } -1 \leq u \leq 1 \\ -1 & \text{if } u < -1 \end{array}\right.$ & ~\cite{htanh} \\
\cline{1-1}
\cline{2-2}
\cline{3-3}
\noalign{\vskip8pt}
hard sigmoid & $\sigma(u) = \left\{ \begin{array}{cl} 1 & \text{if } u > \frac{1}{2} \\ u + \frac{1}{2} & \text{if } -\frac{1}{2} \leq u \leq \frac{1}{2} \\ 0 & \text{if } u < -\frac{1}{2} \end{array}\right.$ & ~\cite{courbariaux2015bc} \\
\cline{1-1}
\cline{2-2}
\cline{3-3}
\noalign{\vskip8pt}
max pooling & $\begin{array}{c} \sigma(u_1, \ldots, u_k) = \max\{0, u_1, \ldots, u_k\} \\ \text{(each $u_i$ is the output of another neuron)}
\end{array}$ & ~\cite{maxpooling} \\
\cline{1-1}
\cline{2-2}
\cline{3-3}
\noalign{\vskip8pt}
maxout  & $\begin{array}{c} \sigma(u_1, \ldots, u_k) = \max\{u_1, \ldots, u_k\} \\ \text{(each $u_i$ is an affine function)} \end{array}$ & ~\cite{Goodfellow2013} \\
\cline{1-1}
\cline{2-2}
\cline{3-3}
\end{tabular}
\end{table}

\subsection{When deep learning meets polyhedral theory} 
It is commonly accepted in machine learning that a simpler model is preferred if it trains as well as a more complex one, since a simpler model is less likely to overfit. 
Conveniently, 
the successful return of neural networks to relatively simpler activation functions 
prepared the ground 
for deep learning to meet polyhedral theory. 
In other words, 
we are now able to analyze and leverage neural networks through the same lenses and tools that have been successfully used for linear and discrete optimization in operations research for many decades. 
We explain this connection in more detail and 
some lines of research that it has opened up in Section~\ref{sec:poly}. 

\subsection{Scope of this survey and related work}\label{sec:scope}
The interplay between mathematical optimization and machine learning has also been discussed by other recent surveys. \cite{CombOptTour} review the use of machine learning in mathematical optimization, whereas~\cite{OptimizationSurvey} formulate mathematical optimization problems with the main focus of obtaining machine learning models, such as by training neural networks. 
A similar scope has been previously surveyed by~\cite{curtis2017optimization} and~\cite{bottou2018optimization}. 
Our survey complements those by focusing exclusively on neural networks while outlining how linear optimization can be used more broadly in that context, from network training and verification to model embedding and compression, as well as refined through formulation strengthening. 
In addition, we illustrate how polyhedral theory can ground the use of such linear formulations and also provide a more nuanced understanding of the discriminative ability of neural networks.  

The presentation in this survey is {\color{black}mainly} centered on {\color{black}fully-connected} \emph{feedforward rectifier networks} {\color{black}with connections between subsequent layers only}. These are very commonly used networks with only ReLU activations and for which most polyhedral results and applications of linear optimization are known. The focus on a single type of neural network is intended to help the reader capture the intuition behind different developments and understand the nuances involved. 
Despite our focus, there are many variants of interest with fewer or different types of connections that can be interpreted as a special case of fully-connected models. For example, the units of Convolutional Neural Networks~(CNNs or ConvNets) \citep{CNN} have local connectivity:
only a subset of adjacent units defines the output of each unit in the next layer, and the same parameters are used to define the output of different units. In fact, multiple \emph{filters} of parameters can be applied to a set of adjacent units through the output of different units in the next layer. 
CNNs are often applied to identify and aggregate the same local features in different parts of a picture, 
and we can interpret them as a special case of feedforward networks. 
Another common variant, the Residual Network~(ResNet) \citep{He2016DeepRL}, includes \emph{skip connections} that directly connect units in nonadjacent layers. 
Those connections can be emulated by adding units passing their outputs through the intermediary layers. 
Hence, many of the results and applications discussed along the survey are relevant to other variants (e.g., LTU and maxout activations, or those other connectivity patterns); {\color{black}we briefly discuss some of these variations} and provide references to more specific results and applications involving them. 

We also discuss the extent to which other variants remain relevant or can be analyzed through the same lenses. 
For example, \emph{feedback connections} in \emph{recurrent networks} \citep{recurrent1,recurrent2} allow the output of a unit to be used as an input of units in previous layers. 
Recurrent networks such as Long Short-Term Memory~(LSTM) \citep{LSTM} produce outputs that depend on their internal state, and they may consequently process sequential inputs with arbitrary length. 
While feedback connections may not be emulated with a \change{feedforward} network, we discuss in the following paragraph how recurrent networks have been replaced with great success by attention mechanisms, which are implemented with feedforward networks. 
In the realm of variants that remain relevant, it is very common to apply a different type of activation to the output layer of the network, such as the layer-wise softmax function $\sigma : \mathbb{R}^{n_{L}} \rightarrow \mathbb{R}^{n_{L}}$ in which $\sigma(u)_i = e^{u_i}/\sum_{j=1}^{n_{L}} e^{u_j} ~ \forall i \in \{1, \ldots, n_{L}\}$ \citep{softmax}, 
which is used to normalize a multidimensional output as a probability distribution. 
While softmax is not piecewise linear, 
we describe how its output can also be analyzed from a polyhedral perspective.

\paragraph{Other uses of deep learning}
Deep learning is also being used in machine learning beyond the realm of supervised learning. 
In \emph{unsupervised learning}, the focus is on drawing inferences from unlabeled datasets. 
For example, 
Generative Adversarial Networks~(GANs) \citep{GAN} have been used to generate realistic images using a pair of neural networks. 
One of these networks is a \emph{discriminator} trained to identify elements from a dataset and the other is a \emph{generator} aiming to mislead the discriminator with synthetic inputs that could be classified as belonging to the dataset. 

In \emph{reinforcement learning}, the focus is on modeling agents that can interact with their environment through actions and associated rewards. Examples of such applications include neural networks designed for the navigation of self-driving vehicles \citep{gao2020vectornet} and for playing Atari games \citep{Atari}, more contemporary electronic games such as Dota 2 \citep{dota2} and StarCraft II \citep{starcraft2}, and the game of Go \citep{AlphaGo} at levels that are either better or at least comparable to human players.  

A more recent and popular example are generative transformers \citep{GPT}, such as DALL·E 2 \citep{dalle2} producing realistic images from text prompts in mid-2022 
and ChatGPT \citep{ChatGPT} producing realistic dialogues with users in early 2023, 
the latter belonging to the fast-growing family of large language models. 
\change{These models use a new type of architecture, 
which preserves context in a different way from LSTMs: 
they rely on attention heads aimed at scoring the relevance of past states \citep{bahdanau2015translation}, which is the foundation of the transformer architecture \citep{vaswani2017attention}.}

\paragraph{Further reading}
For a historical perspective on neural networks, we recommend \cite{historical}. 
For a recent and broad introduction to the fundamentals of deep learning, we recommend~\cite{zhang2023dive}. 
For other forms of measuring model complexity in neural networks, we refer to~\cite{hu2021complexity}.

\section{The Polyhedral Perspective}\label{sec:poly}

A feedforward rectifier network models a piecewise linear function \citep{arora2018understanding} in which every such piece is a polyhedron \citep{raghu2017expressive}, 
and represents a special case among neural networks modeling piecewise polynomials \citep{balestriero2018spline}. 
Therefore, training a rectifier network is {\color{black}essentially} performing a piecewise linear regression {\color{black}(where the architecture limits the piecewise linear functions considered)},  
and we can potentially interpret such neural networks in terms of what happens in each piece of the function that they model.
However, we are only beginning to answer some of the questions entailed by such a remark. In this survey, we discuss how insights on this subject may help us answer the following questions.
\begin{enumerate}
\item Which piecewise linear functions can or cannot be obtained from training a neural network given its architecture?
\item Which neural networks are more susceptible to adversarial exploitation?
\item Can we integrate the model learned by a neural network into a broader decision-making problem for which we want to find an optimal solution?
\item Is it possible to obtain a smaller neural network that models exactly the same function as another trained  neural network?
\item Can we exploit the polyhedral geometry present in neural networks in the training phase? 
\item Can we efficiently incorporate extra structure when training neural network, such as linear constraints over the weights?
\end{enumerate}

The first question complements the universal approximation results for neural networks. Namely, there is a limit to what functions can be well approximated when limited computational resources are translated into constraints on the depth and width of the layers of neural networks that can be used in practice. 
The functions that can be modeled depend on the particular choice of hyperparameters subject to the computational resources available, and in the long run that may also lead to a more principled approach for the choice of hyperparameters than the current approaches of neural architecture search. 
In Section~\ref{sec:LR}, 
we analyze how a rectifier network partitions the input space into pieces in which it behaves linearly, which we denote as \emph{linear regions}. 
We discuss the geometry of linear regions, the effect of parameters and hyperparameters on the number of linear regions of a neural network, and the extent to which such number of linear regions relates to the accuracy of the network.

The second question relies on formal verification methods to evaluate the robustness of neural networks, which can be approached with mathematical optimization formulations that are also relevant for the third and fourth questions. Such formulations are convenient since a direct inspection of every piece of the function modeled by a neural network is prohibitive given how quickly their number \change{scales} with the size of the network. 
The linear behavior of the network for every choice of active and inactive units implies that we can use linear formulations with binary variables corresponding to the activation of units to model trained neural networks using MILP.  Therefore, we are able to solve a variety of optimization problems over a trained neural network, such as the neural network verification problem, 
identifying the range of outputs for each ReLU of the network, and modeling a trained neural network as part of a larger decision-making problem. 
In Section~\ref{sec:optimizing}, we discuss how to formulate optimization problems over a trained neural network, the applications of such formulations, and the progress toward obtaining stronger formulations that scale more easily with the network size.

The fifth and sixth questions involve the training procedure of a DNN, where linear programming tools have been applied to partially answer them. In Section~\ref{sec:training}, we overview these developments. In terms of the fifth question ---exploiting polyhedrality in training neural networks--- we describe algorithms that use the polyhedral geometry induced by activation sets to solve training problems. We also cover a recently proposed polyhedral construction that can approximately encode multiple training problems at once, showing a strong relationship across training problems that arise from different datasets, for a fixed architecture. Additionally, we review some recent uses of mixed-integer linear programming in the training phase as an alternative to SGD when the weights are required to be integer. Regarding the sixth question ---the incorporation of extra structure when training--- we review multiple approaches that have included techniques related to linear programming within SGD to impose a desirable structure when training, or to find better step-lengths in the execution of SGD.

\section{Piecewise Linear Structure and the Linear Regions of a Neural Network}\label{sec:LR}

Every piece of the piecewise linear function modeled by a neural network is a linear region, 
and ---without loss of generality--- we can think of it as a polyhedron. 
In this section, we define a linear region, discuss what piecewise linear functions can be represented by different architectures, exemplify how linear regions can be so numerous, and what may affect which functions can be represented as well as the number of linear regions in a neural network. 
We also discuss the practical implications of such insights, as well as other related forms of analyzing the ability of a neural network to represent expressive models. 
\change{For the next definition, we remind the reader that a ReLU is said to be \emph{active} when its output is positive.}

\begin{definition}
A linear region corresponds to the set of points from the input space that activate the same units along the neural network, 
and hence can be characterized by the set $\sS^l$ of units that are active in each layer $l \in \sL$.
\end{definition}

\change{ Linear regions are fundamental when analyzing the behavior of a neural network.}
If we restrict the domain of a neural network to a linear region $\sI \subseteq \mathbb{R}^{n_0}$, then the neural network behaves as an affine transformation $\vy_\sI : \sI \rightarrow \mathbb{R}^{n_{L}}$ of the form $\vy_\sI(\vx) = \mT \vx + \vt$ with a matrix $\mT \in \mathbb{R}^{n_{L} \times n_0}$ and a vector $\vt \in \mathbb{R}^{n_{L}}$ that are directly defined by the network parameters and the set of neurons that are activated by any input $\vx \in \sI$. 
For a small perturbation~$\varepsilon$ to some input $\overline{\vx} \in \sI$ such that $\overline{\vx}+\varepsilon \in \sI$, 
the network output for $\bar{x}+\varepsilon$ is given by $\vy_\sI(\overline{\vx}+\varepsilon)$. 
While it is possible that two adjacent regions defined in such way correspond to the same affine transformation, 
thinking of each linear region as having a distinct signature of active units makes it easier to analyze them. 
\newchange{
\cite{LocallyLinear} and \cite{hanin2019deep} have both discussed such distinction, 
and \cite{stargalla2025counting} present additional definitions for the term ``linear region" found in the literature. 
}

The number of linear regions defined by a neural network is one form with which we can measure the complexity of the models that it can represent \citep{DeepArchitectures}. 
Hence, if a more complex model is desired, we may want to design a neural network that can potentially define more linear regions. 
On the one hand, the number of linear regions may grow exponentially \change{with} the depth of a neural network. 
On the other hand, such a number depends on the interplay between network parameters and hyperparameters. 
As we consider how the inputs from adjacent linear regions are evaluated, 
the change to the affine transformation can be characterized in algebraic and geometric terms. 
Understanding such changes may help us grasp how a neural network is capable of telling its inputs apart, including what are the sources of the complexity of the model. 

For neural networks in which the activation function is not piecewise linear, 
\cite{Bianchini2014} have used more elaborate topological measures to compare the expressiveness of shallow and deep neural networks. 
\change{ On that note, a recent study by \cite{ergen2024topological} has also explored topological measures in neural networks with ReLU activations.}  
\cite{hu2020curve} followed a closer approach by producing a linear approximation neural network in which the number of linear regions can be counted.

\subsection{The combinatorial aspect of linear regions}

One of the most striking aspects about analyzing a neural network in terms of its linear regions is how quickly such number grows. 
Early work on this topic by~\cite{pascanu2013on} and \cite{montufar2014on} has drawn two important observations.
First, that it is possible to construct simple deep neural networks with a number of linear regions that grows exponentially in the depth.
Second, that the number of linear regions can be exponential in the number of neurons alone\change{, i.e., when the depth is fixed}. 

\begin{figure}
    \centering
    \includegraphics[width=0.6\textwidth]{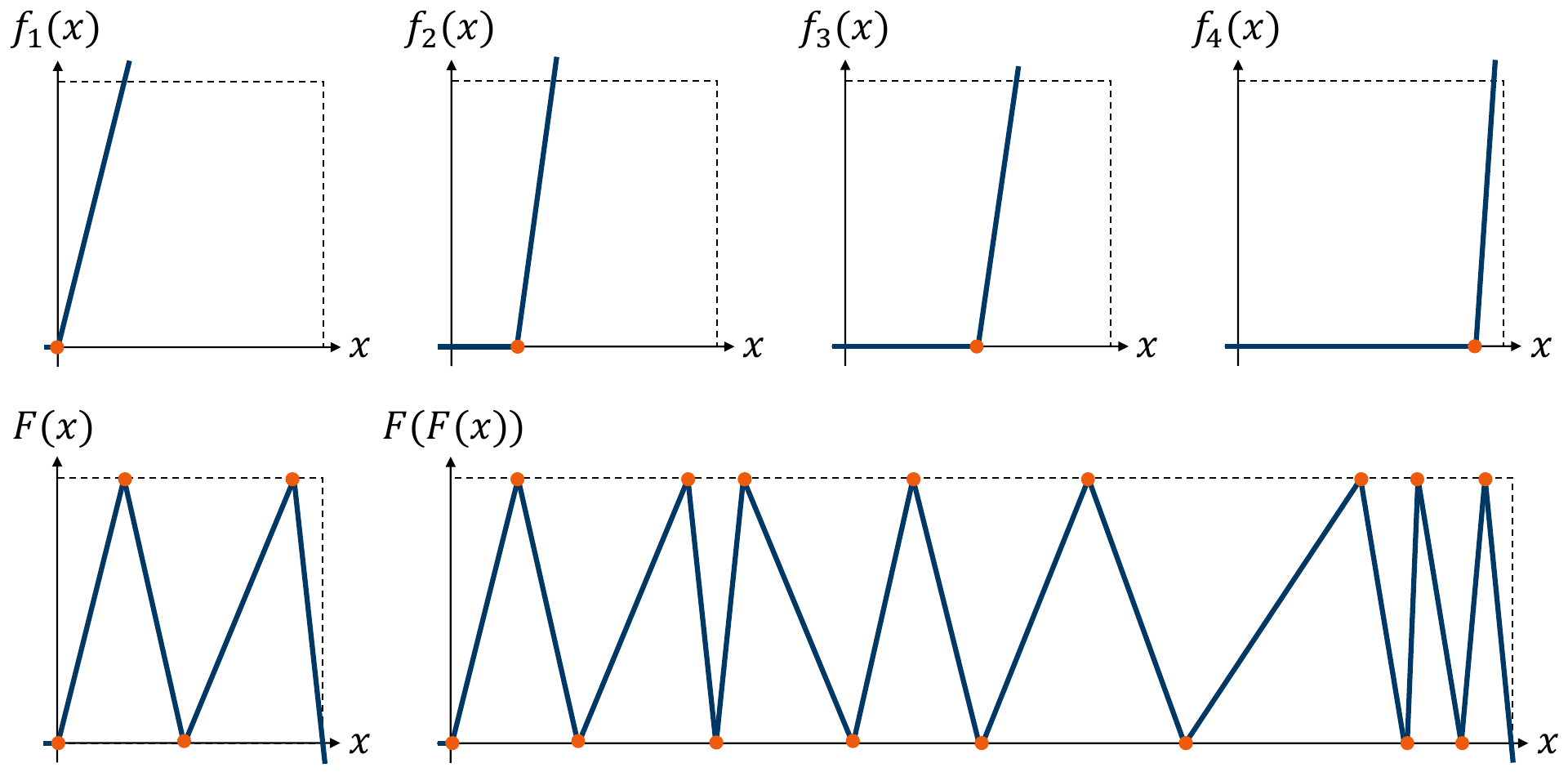}
    \caption{\change{ Output functions} $\{ f_i(x) \}_{i=1}^4$ of the units in the first layer and combined outputs of the first two layers ---$F(x) = f_1(x) - f_2(x) + f_3(x) - f_4(x)$ for the first and $F(F(x))$ for the second--- of a neural network in which the number of linear regions grows exponentially on the depth, as described in Example~\ref{ex:zigzag}.}
    \label{fig:zigzag}
\end{figure}

\begin{figure}
    \centering
    \includegraphics[width=\textwidth]{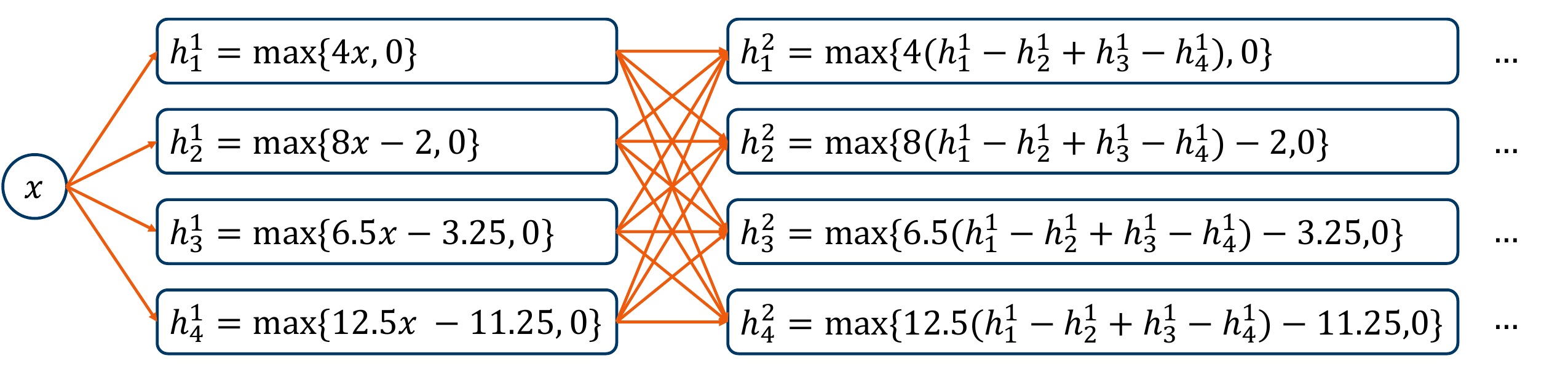}
    \caption{Mapping from the input $x \in [0,1]$ to the intermediary output $\vh^2 \in [0,1]^{4}$ through the first two layers of a neural network in which the number of linear regions \change{grows} exponentially \change{with} the depth, as described in Example~\ref{ex:zigzag}. 
    The parameters of subsequent layers are the same as those in the second layer.}
    \label{fig:architecture_zigzag}
\end{figure}

The first observation comes from analyzing the role of ReLUs in a very simple setting. 
Namely, that of a neural network in which we regard every layer as having a single input in the $[0,1]$ domain, which is produced by combining the outputs of the units from the preceding layer, 
as illustrated by Example~\ref{ex:zigzag}. 
\begin{example}\label{ex:zigzag}
Consider a neural network with input $x$ from the domain $[0,1]$ and layers having 4 neurons with ReLU activation. 
\change{ With the choice of parameters described in the next paragraph, 
we obtain a zigzagging function with 4 slopes in the $[0,1]$ domain (see $F(x)$ in Figure~\ref{fig:zigzag}). 
Each of the slopes of $F(x)$ defines a bijection between segments of the input ---namely, $[0,0.25]$, $[0.25, 0.5]$, $[0.5, 0.9]$, and $[0.9, 1.0]$--- and the image $[0,1]$. 
The effect of repeating such structure in the second layer is that of composing the zigzagging function with itself (see $F(F(x))$ in Figure~\ref{fig:zigzag}), 
with 4 slopes being produced within each of those 4 initial segments. 
Hence, the number of slopes ---and therefore of linear regions--- in the output of such a neural network with $L$ layers of activation functions is $4^L$, 
which implies an exponential growth on depth. 
}

\change{
The choice of parameters for the first and subsequent layers is illustrated in Figure~\ref{fig:architecture_zigzag}. 
For the first layer, the output of the neurons are given by the following functions: 
$f_1(x)=\max\{4x,0\}$, $f_2(x)=\max\{8x-2,0\}$, $f_3(x)=\max\{6.5x-3.25,0\}$, and $f_4(x)=\max\{12.5x-11.25,0\}$ (see first row of charts in Figure~\ref{fig:zigzag})}. 
In other words, $\vh^1_i = f_i(x) ~\forall i \in \{1,2,3,4\}$. 
For the subsequent layers, assume that the outputs coming from the previous layer are combined through the function $F(x)=f_1(x)-f_2(x)+f_3(x)-f_4(x)$ \change{ (see $F(X)$ in Figure~\ref{fig:zigzag})}, 
which substitutes $x$ as the input to the next layer; upon which the same set of functions $\{ f_i(x) \}_{i=1}^4$ defines the output of the next layer. 
In other words, $\vh^l_i = f_i(F(\vh^{l-1})) = f_i(\evh_1^{l-1} - \evh_2^{l-1} + \evh_3^{l-1} - \evh_4^{l-1})~\forall i \in \{1,2,3,4\}, l \in \sL \setminus \{ 1 \}$. 
\change{ Hence, $(F^2)(x) = F(F(x)) = h_1^1-h_2^1+h_3^1-h_4^1$ (see $F(F(x))$ in Figure~\ref{fig:zigzag}). }
\end{example}
In Example~\ref{ex:zigzag}, every neuron changes the slope of the resulting function once it becomes active, 
in which we purposely alternate between positive and negative slopes once the function reaches either 0 or 1, respectively. 
By selecting the network parameters accordingly, 
\cite{montufar2014on} were the first to show that a layer with $n$ ReLUs can be used to create a zigzagging function with $n$ slopes on the $[0,1]$ domain, 
with the image along every slope also corresponding to the interval $[0,1]$. 
Consequently, stacking $L$ of such layers results in a neural network with $n^L$ linear regions.  

The second observation ---that the number of linear regions can grow exponentially in the number of neurons alone--- comes from the interplay between the parts of the input space in which each of the units are active, 
especially in higher-dimensional spaces. 
This is based on some geometric observations that we discuss in Section~\ref{sec:geometry}. 
Even for a \emph{shallow} network ---i.e., the number of layers being $L=1$--- such a number of linear regions may approach $2^n$, 
which corresponds to every possible activation set $\sS \subseteq \{1, \ldots, n\}$ defining a nonempty linear region. 
However, as we discuss later, that is not always the case due to architectural choices such as the number of layers and their width.

\subsection{The algebra of linear regions}
\label{sec:algebraoflinearegions}

Given the activation sets $\{ \sS^l \}_{l \in \sL}$ denoting which neurons are active for each layer of the neural network, 
we can explicitly describe the affine transformation $\vy_\sI(\vx) = \mT \vx + \vt$ associated with the corresponding linear region $\sI$. 
For every activation set $\sS^l$, layer $l$ defines an affine transformation of the form $\Omega^{\sS^l}(\mW^l \vh^{l-1} + \vb^l)$, where $\Omega^{\sS^l}$ is a diagonal $n_l \times n_l$ matrix in which $\Omega^{\sS^l}_{ii} = 1$ if $i \in \sS^l$ and $\Omega^{\sS^l}_{ii} = 0$ otherwise. 
Hence, the matrix $\mT$ and vector $\vt$ are as follows:
\[
\mT = \prod_{l=1}^L \Omega^{\sS^l} \mW^{l},
\qquad
\vt = \sum_{l_1=1}^{L} \left( \prod_{l_2=l_1+1}^{L} \Omega^{\sS^{l_2}} \mW^{l_2} \right) \Omega^{\sS^{l_1}} \vb^{l_1}.
\]
On a side note, 
\cite{takai2021functions} proposed \change{ refining the study of linear regions by identifying linear regions associated with equivalent affine transformations upon permutation of the neurons. 
In other words, that would entail identifying and working with equivalence classes of linear regions and focusing on the number of such classes instead. For simplicity, we assume every linear region to be distinct from other linear regions.}

Each linear region is associated with a polyhedron, 
and we can describe the union of polyhedra $\mathcal{D}$ on $(\vx, \vh^1, \ldots, \vh^L)$, \change{ with $\vx = \vh^0$,} that covers the entire space \change{ $\mathbb{R}^{n_0}$ for the input $\vx$} of the neural network as follows:
\[
\mathcal{D} = 
\bigvee_{\substack{(\sS^1, \ldots, \sS^{L}) \subseteq \\ \{1, \ldots, n_1\}  \\ \times \ldots \times \\ \{1, \ldots, n_{L}\} } }
\left(
\begin{array}{cc}
\vw_i^l \cdot \vh^{l-1} + b_i^l \geq 0 & \forall l \in \sL, i \in \sS^l \\
h_i^l = \vw_i^l \cdot \vh^{l-1} + b_i^l & \forall l \in \sL, i \in \sS^l \\
\vw_i^l \cdot h^{l-1} + b_i^l \leq 0 & \forall l \in \sL, i \notin \sS^l \\
h_i^l = 0 & \forall l \in \sL, i \notin \sS^l 
\end{array}
\right).
\]
Such partitioning entails an overlap between adjacent linear regions when $\vw_i^l \vh^{l-1} + b_i^l = 0$, i.e., at the boundary in which unit $i$ in layer $l$ is active in one region and inactive in another. 
Nevertheless, for any input $\overline{\vx}$ associated with a point at such a boundary between two linear regions $\sI_1$ and $\sI_2$, it holds that $\vy_{\sI_1}(\overline{\vx}) = \vy_{\sI_2}(\overline{\vx})$ even if those affine transformations are not entirely identical since the output of the neural network is continuous. 
More importantly, such overlap implies that each term of $\mathcal{D}$ is defined using only equalities and nonstrict inequalities. 
\change{In other words, each linear region corresponds to a proper closed polyhedra in the extended $(\vx, \vh^1, \ldots, \vh^{L})$ space due to the overlap of boundaries between pairs of adjacent linear regions.}  
Consequently, those linear regions also define polyhedra if projected to the \change{$\vx$ space}, 
since by using Fourier-Motzkin elimination \citep{Fourier,Motzkin} 
we obtain a polyhedral description of the linear region in $\vx$. 
\change{If one of those polyhedra in the $\vx$ space does not have an interior, 
which means that it is not full-dimensional in $\vx$, 
then that polyhedron lies entirely on the boundary of other linear regions.
In such a case, we do not regard it as a proper linear region. } 
\newchange{The interiors of all full-dimensional linear regions are disjoint.} 

\change{Working with linear regions in $(\vx, \vh^1, \ldots, \vh^{L})$ space is convenient since we can explicitly define the linear regions with two linear constraints per neuron, as in the description of $\mathcal{D}$. 
Moreover, since the values for $\vx$ imply the values for $\vh^1, \ldots, \vh^L$, 
as the former is the input of the neural network and the latter are the outputs of each layer, 
then full-dimensional regions in $\vx$ are also those defining $n_0$-dimensional polyhedra in $(\vx, \vh^1, \ldots, \vh^{L})$.}
By looking at the geometry of those linear regions from a different perspective in Section~\ref{sec:geometry} and understanding its impact on the number of \change{full-dimensional} linear regions in Section~\ref{sec:number}, we will see that many terms of $\mathcal{D}$ may actually be \change{either empty or a proper subset of other terms}.

\change{Modeling a neural network as the union of polyhedra
is convenient for modeling optimization problems involving neural networks, 
as we discuss later in Section~\ref{sec:MIPmodels}.
Optimization over the union of polyhedra is the subject of disjunctive programming, 
which has been applied to formulate MILP~\citep{DPBook} and Mixed-Integer Non-Linear Programming~(MINLP) formulations~\citep{gdp1994first,gdp2012survey},
as well as for generating valid inequalities to strengthen those~\citep{CGLP1,CGLP2}.}

\subsection{The geometry of linear regions}\label{sec:geometry}

Another form of looking at the geometry of linear regions is through their transformation along the layers of the neural network. Namely, we can think of the input space as initially being partitioned by the units of the first layer, and then each resulting linear region being further partitioned by the subsequent layers. 
In that sense, we can think of every layer as a particular form of ``slicing'' the input. 
In fact, a layer may slice each linear region that is defined by the preceding layer in a different way due to which neurons are active or not in previous layers.

Let us begin by illustrating how a given layer $l \in \sL$ partitions its input space\change{, or $\vh^{l-1}$ space}. Every neuron $i$ in layer $l$ is associated with an \emph{activation hyperplane} of the form $\vw_i^l \cdot \vh^{l-1} + b_i^l = 0$, which divides the possible inputs of its layer into an open half-space in which the unit is active ($\vw_i^l \cdot \vh^{l-1} + b_i^l > 0$) and a closed half-space in which the unit is inactive ($\vw_i^l \cdot \vh^{l-1} + b_i^l \leq 0$). These hyperplanes define the boundary between adjacent linear regions, and the arrangement of such hyperplanes for a given layer $l \in \sL$ determines how that layer partitions the $\vh^{l-1}$\change{ space}. 
In other words, \change{every input $\vh^{l-1}$} can be located with respect to each of those hyperplanes, which corresponds to the activation set of the linear region to which it belongs. 
However, not every activation set out of the $2^{n_l}$ possible ones maps to a nonempty region of the input space. In the case of Example~\ref{ex:hyperplane_arrangement}, there is no linear region in which the activation set is empty.

\begin{figure}
    \centering
    \includegraphics[width=0.25\textwidth]{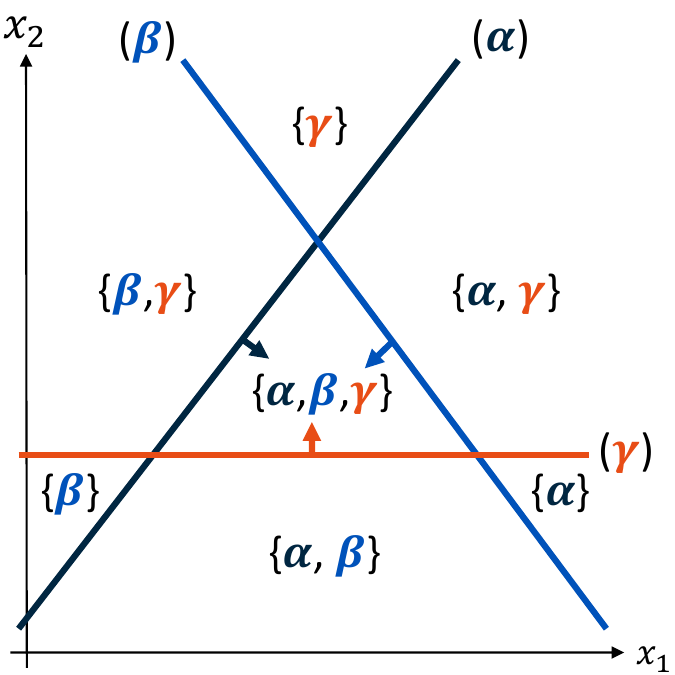}
    \caption{Linear regions defined by the shallow neural network described in Example~\ref{ex:hyperplane_arrangement}. 
    Every line corresponds to the activation hyperplane of a different neuron, which is given by $\alpha$, $\beta$, and $\gamma$ in parentheses. 
    The arrow next to each line points to the half space in which the inputs activate that neuron.
    Every linear region has a subset of $\{\alpha,\beta,\gamma\}$ as its corresponding activation set.}
    \label{fig:hyperplane_arrangement}
\end{figure}

\begin{example}\label{ex:hyperplane_arrangement}
Consider a neural network with domain $\vx \in \mathbb{R}^2$ and a single layer having 3 neurons $\alpha$, $\beta$, and $\gamma$ with outputs given as follows: 
$h^1_{\alpha} = \max\{ 2.3 x_1 - 1.9 x_2 +0.6, 0\}$, $h^1_{\beta} = \max\{ -0.9 x_1 - 0.7 x_2 +4.8, 0\}$, and $h^1_{\gamma} = \max\{ 0 x_1 + 3 x_2 -5, 0\}$. 
These neurons define the activation hyperplanes ($\alpha$) $2.3 x_1 - 1.9 x_2 +0.6 = 0$, ($\beta$) $-0.9 x_1 - 0.7 x_2 +4.8 = 0$, and ($\gamma$) $0 x_1 + 3 x_2 -5 = 0$ in the \change{$\vx$ space}, 
which are illustrated along with the activation sets of the linear regions in Figure~\ref{fig:hyperplane_arrangement}.

Instead of $2^3$ linear regions corresponding to each possible activation set defined by a subset of the neurons in $\{\alpha, \beta, \gamma\}$, the arrangement of such hyperplanes produces 7 linear regions, 
which is in fact the maximum number of 2-dimensional regions that can be defined by drawing 3 lines on a plane. 
\end{example}

The maximum number of full-dimensional regions resulting from a partitioning defined by $n$ hyperplanes depends on the dimension $d$ of the space in which those hyperplanes are defined \citep{Zaslavsky1975}. 
That number never exceeds $\sum\limits_{i=1}^{\min\{d,n\}} \binom{n}{i}$. Such bound only coincides with $2^n$ if $d \geq n$; otherwise, as illustrated in Example~\ref{ex:hyperplane_arrangement}, that number can be smaller. 
As observed by~\cite{hanin2019deep}, 
that bound is $O\left( \frac{n^d}{d!} \right)$ when $n \gg d$.

In fact,
the above bound is all that we need to determine the maximum number of linear regions in shallow networks. 
While not every shallow network may define as many linear regions, it is always possible to put the hyperplanes in what is called a \emph{general position} in order to reach that bound. 
Thus, the maximum number of linear regions defined by a shallow network is \change{$\sum\limits_{i=0}^{\min\{n_0,n_1\}} \binom{n_1}{i}$}.
For the polyhedron associated with each linear region, 
being in general position implies that each vertex lies on exactly $d$ activation hyperplanes.

In the case of deep networks, the partitioning of each linear region by the subsequent layers is based on the output of that linear region. This affects the shape and the number of the linear regions defined by the following layers, which may vary between adjacent linear regions due to which units are active or inactive from one linear region to another, as illustrated in Example~\ref{ex:linear_regions}. 

\begin{figure}
    \centering
    \includegraphics[width=0.24\textwidth]{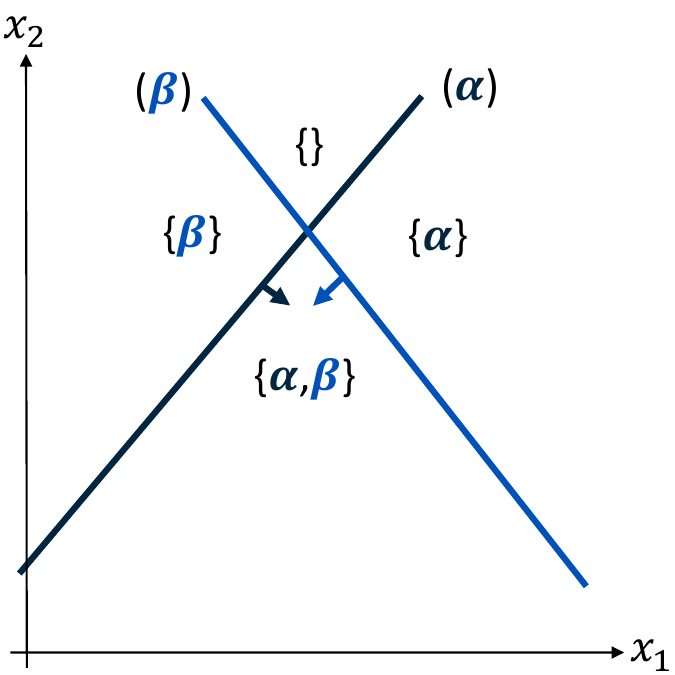}
    \includegraphics[width=0.24\textwidth]{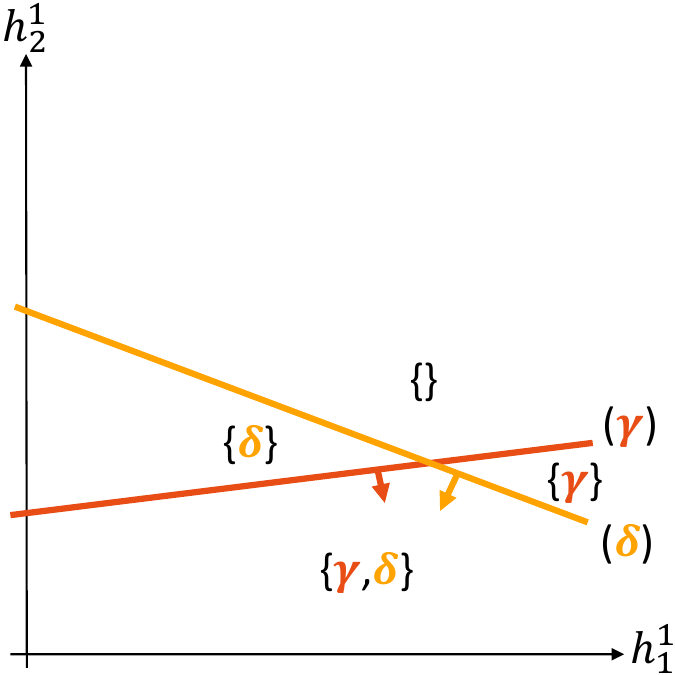} 
    \includegraphics[width=0.24\textwidth]{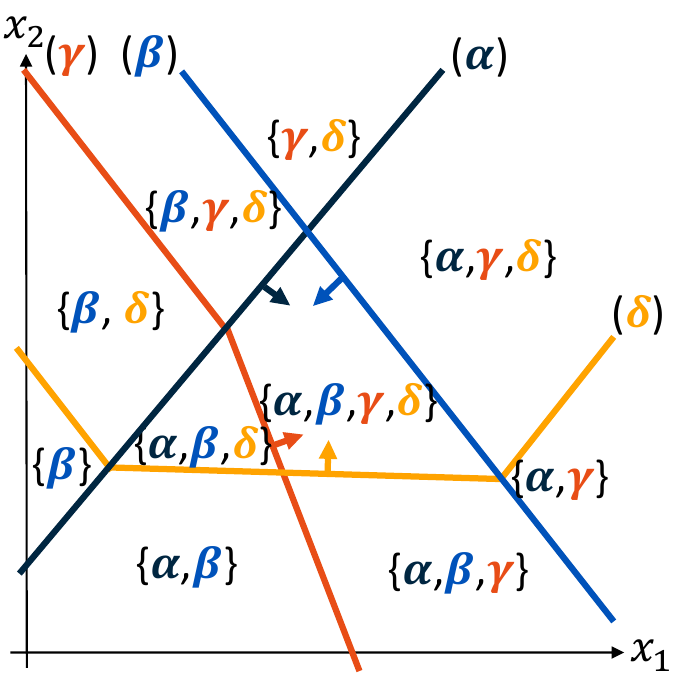}
    \caption{Linear regions defined by the 2 layers of the neural network described in Example~\ref{ex:linear_regions}, 
    following the same notation as in Figure~\ref{fig:hyperplane_arrangement}. The first and second plots show the linear regions and corresponding activation sets defined by the first and the second layers in terms of their input spaces (\change{the $\vx$ and $\vh^1$ spaces}). The third plot shows the linear regions defined by the combination of the 2 layers and the union of their activation sets in terms of the input space of the first layer (\change{the $\vx$ space}).}
    \label{fig:linear_regions}
\end{figure}

\begin{figure}
    \centering
    \includegraphics[width=0.45\textwidth]{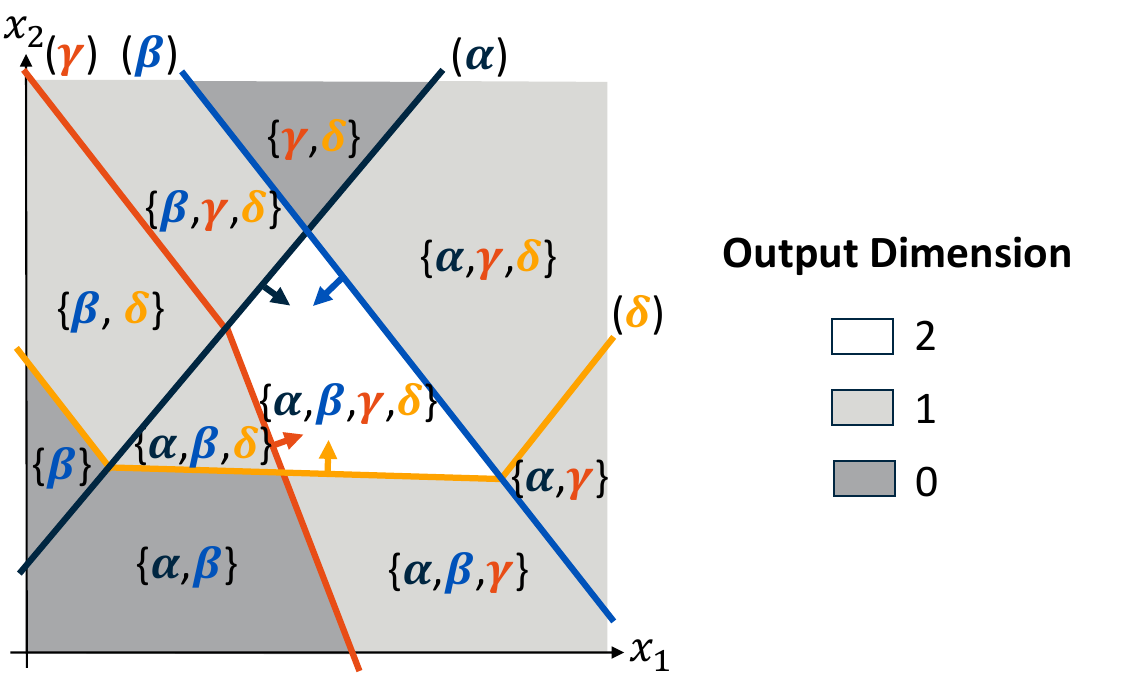}
    \caption{Dimension of the image of the affine function $\vy_{\sI}(\vx)$ associated with each linear region $\sI$ defined by the neural network described in Example~\ref{ex:linear_regions}. The linear regions are the same illustrated in the third plot of Figure~\ref{fig:linear_regions}.}
    \label{fig:dimensions}
\end{figure}

\begin{example}\label{ex:linear_regions}
Consider a neural network with domain $\vx \in \mathbb{R}^2$ and 2 layers having 2 neurons each ---say neurons $\alpha$ and $\beta$ in layer 1, and neurons $\gamma$ and $\delta$ in layer 2--- with outputs given as follows: 
$h^1_{\alpha} = \max\{ 2.3 x_1 - 1.9 x_2 +1.5, 0\}$, $h^1_{\beta} = \max\{ -0.9 x_1 - 0.7 x_2 +5, 0\}$, 
$h^2_{\gamma} = \max\{ 0.4 h^1_1 - 3.1 h^1_2 +4, 0\}$, $h^2_{\delta} = \max\{ -0.6 h^1_1 - 1.6 h^1_2 +5, 0\}$.
These neurons define the activation hyperplanes ($\alpha$) $2.3 x_1 - 1.9 x_2 +1.5 = 0$ and ($\beta$) $-0.9 x_1 - 0.7 x_2 +5 = 0$ in the $\vx$ space 
and the activation hyperplanes ($\gamma$) $0.4 h^1_1 - 3.1 h^1_2 +4 = 0$ and ($\delta$) $-0.6 h^1_1 - 1.6 h^1_2 +5 = 0$ in the $\vh^1$ space, 
which are illustrated along with the activation sets of the linear regions in the first two plots of Figure~\ref{fig:linear_regions}. 
The third plot illustrates the linear regions jointly defined by the two layers in terms of the input space\change{, or $\vx$ space}.  

The third plot is repeated in Figure~\ref{fig:dimensions}, in which the shape of each linear region $\sI$ is filled in accordance to the dimension of the image of $\bar{\vy}_{\sI}(\vx)$ ---the output of the neural network for each linear region $\sI$.
\end{example}

Example~\ref{ex:linear_regions} highlights two important aspects about the structure of linear regions in deep neural networks. 
First, the linear regions defined by a neural network with multiple layers are different because activation hyperplanes after the first layer may look ``bent'' from the input space\change{, or $x$ space}, 
such as with the inflections of hyperplanes $(\gamma)$ and $(\delta)$ in the third plot of Figure~\ref{fig:linear_regions} from one linear region defined by the first layer to another. 
This partitioning of the input space would not be possible with a single layer. 

By comparing side by side the first and the third plots of Figure~\ref{fig:linear_regions}, 
we can see how every linear region of a given layer of a neural network may be partitioned differently by the following layer. 
When defined in terms of the input space\change{, or $\vx$ space}, 
the hyperplanes associated with the second layer differ across the linear regions defined by the first layer 
because each of those linear regions is associated with a different affine transformation from $\vx$ to $\vh^1$. 
Hence, the activation hyperplanes of layer $l$ may break each linear region from layer $l-1$ differently. 
To every linear region defined by the hyperplane arrangement in the $\vh^{l-1}$ space there is a linear transformation $\vh^{l-1} = \Omega^{\sS^{l-1}}(\mW^{l-1} \vh^{l-2} + \vb^{l-1})$ to the points of that linear region based on the set of active neurons $\sS^{l-1}$. 
Consequently, inputs in the $\vh^{l-1}$ space that are associated with different linear regions are transformed differently to the $\vh^l$ space, and therefore the form in which those linear regions are further partitioned by layer $l$ is not the same when seen from the $\vh^{l-1}$ space. 

Second, some combinations of activation sets of multiple layers do not correspond to linear regions even if the activation hyperplanes are in general position with respect to each layer. 
For each layer, the first two plots of Figure~\ref{fig:linear_regions} show that every activation set corresponds to a nonempty region of the layer input. 
However, not every pair of such activation sets would define a nonempty linear region for the neural network. 
For example, the linear region of the first layer associated with the activation set $\sS^1 = \{\}$ defines a linear region in $\vx$ which is always mapped to $\vh^1=0$, 
and thus only corresponds to activation set $\sS^2 = \{ \gamma, \delta \}$ in the second layer because both units are active for such input. 
Thus, no linear region in $\vx$ is associated with only the units in sets $\{\}, \{\gamma\}$, and $\{ \delta \}$ being active ---i.e., there is no linear region such that $\sS^1 \cup \sS^2 = \{\}, \{\gamma\}, \text{or} \{\delta\}$.

More generally, the number of units that are active on each linear region defined by the first layer also imposes a geometric limit to how that linear region can be further partitioned. 
If only one unit is active at a layer, that means that the output of the layer within that linear region has dimension 1, 
and, consequently, the subsequent hyperplane arrangements within that linear region are limited to a 1-dimensional space.
For the network in the example, we thus expect no more than $\sum_{i=0}^1 \binom{2}{i} = 3$ linear regions being defined instead of $2^2 = 4$ when only one unit is active. 
In fact, 
that is precisely the number of subdivisions by the second layer of the linear region defined by activation set $\sS^1 = \{\beta\}$ from the first layer.

For every linear region defined by layer $l$ with an activation set $\sS^l$, the dimension of the output of the corresponding transformation $\Omega_{\sS^l}(\mW^l \vh^{l-1} + \vb^l)$ is at most $|\sS^l|$ since $\text{rank}(\Omega_{\sS^l}) = |\sS^l|$. Hence, the dimension of the output of every linear region defined by a rectifier network is upper bounded by its smallest activation set across all layers. This phenomenon was first identified by~\cite{serra2018bounding} as the \emph{bottleneck effect}. In neural networks with uniform width, this phenomenon leads to a surprising consequence: the number of linear regions with full-dimensional output is at most one. There are also consequences to the maximum number of linear regions that can be defined, as we discuss later.

\subsubsection{The geometry of decision regions}

\change{
Neural networks are often applied to classification problems. 
In those problems, every input $x \in \mathbb{R}^{n_0}$ must be associated with one class among $\mathcal{C}_1, \mathcal{C}_2, \ldots, \mathcal{C}_C$. 
In rectifier networks, that is usually done with the additional layer-wise softmax function $\sigma(\vh^L)_i = e^{\vh^L_i}/\sum_{j=1}^{n_{L}} e^{\vh^L_j} ~ \forall i \in \{1, \ldots, n_{L}\}$, 
in which the width $n_L$ equals to number of classes $C$. 
Every output element $i \in \{1, \ldots, C\}$ represents the predicted probability of the input of being from class $\mathcal{C}_i$, 
and then the largest output element corresponds to the class to which the input is associated. 
With the addition of the softmax layer, 
the neural network no longer corresponds to a piecewise linear function.
However, the region associated with each class corresponds to the union of polyhedra.}

The set of inputs associated with the same class define a \emph{decision region}. 
In rectifier networks coupled with a softmax layer, 
the decision regions can also be defined by polyhedra\change{, as described next}. 
Although the output of the softmax layer is not piecewise linear, 
its largest output corresponds to its largest input. 
Hence, every linear region $\sI$ defined by layers 1 to $L$ is partitioned by the softmax layer into smaller regions 
where $\vh^{L}_i \geq \vh^{L}_j ~ \forall j \neq i$ for each class $\mathcal{C}_i$\change{, 
which may denote as \emph{polyhedral decision subregions}. 
Therefore, 
the decision region for each class $\mathcal{C}_i$ is the union of polyhedral decision subregions in which $\vh^{L}_i \geq \vh^{L}_j ~ \forall j \neq i$.}

Difficulties in modeling functions such as the Boolean XOR in shallow networks are related to limitations on the form of the decision regions, 
which may be limited by the depth of the neural network. 
For example, \cite{makhoul1989twolayer} showed that two layers suffice to obtain disconnected decision regions. 
In fact, we may say further in the typical setting where no hidden layer is wider than the input ---i.e., $n_0 \geq n_l ~ \forall l \in \sL$: 
\cite{pmlr-v80-nguyen18b} showed that at least one layer $l \in \sL$ must be such that $n_l > n_0$ for the network to present disconnected decision regions. 
\change{Under the same setting and for an input size $n_0 \geq 2$,} 
\cite{grigsby2022topology} showed that the decision regions are either empty or unbounded.

\subsection{Piecewise linear representability}\label{sec:representability}

We have seen that it is possible to build neural networks with a very large number of linear regions. However, that does not necessarily mean that any piecewise linear function with a certain number of pieces can be represented by a given architecture. 
A fundamental question of neural network representability is how the different network parameters (primarily width and depth) can affect the piecewise linear functions it can represent. Here, we detail some of the key results in this regard.

One of the most important representability results is that of \cite{arora2018understanding}, which shows the following equivalence between the families of continuous piecewise \change{linear} functions and neural networks.

\begin{theorem}[\cite{arora2018understanding}]
Every $\mathbb{R}^n \to \mathbb{R}$ ReLU DNN represents a piecewise linear function, and every piecewise linear function $\mathbb{R}^n \to \mathbb{R}$ can be represented by a ReLU DNN with at most $\lceil \log_2(n + 1)\rceil$ hidden layers.
\end{theorem}

This result is remarkable in at least two ways: (1) it provides a clean structural characterization of the expressive power of neural networks, and (2) it provides a bound on the depth that is needed for representing an arbitrary piecewise linear function as a neural network.

At the core of the bound proof are the following ingredients. Firstly, a result from \cite{wang2005generalization} that states that every piecewise linear function $f$ can be represented as $f=\sum_{j=1}^p s_j\left(\max _{i \in S_j} \ell_i\right)$, where each $\ell_i$ is an affine linear function, each $S_i$ is a finite set of cardinality at most $n+1$, and each $s_j\in \{-1,+1\}$. Thus, every such $f$ is representable as a linear combination of convex piecewise linear functions, each of which has at most $n+1$ linear regions. Secondly, the authors prove that the maximum of two numbers $\max\{x,y\}$ can be computed with a 2-layer ReLU neural network. Roughly speaking, applying this inductively (along with other operations that are also representable by ReLU networks), yields the logarithmic upper bound.

A natural follow-up question is whether one can do better than this logarithmic bound. \cite{hertrich2021depth} \newchange{conjectured} that this is not true. 
\begin{conjecture}[\cite{hertrich2021depth}]\label{conjecture1}
    Let $\text{ReLU}_n(k)$ be the set of functions $f:\mathbb{R}^n\to \mathbb{R}$ that can be represented by a network with $k$ hidden layers. Then, for every $k\in \{1,\ldots,\lceil \log_2(n+1) \rceil\} $
    \[\text{ReLU}_n(k-1)\subsetneq \text{ReLU}_n(k).\]
\end{conjecture}

To narrow down the path into proving this conjecture,  \cite{hertrich2021depth} show that Conjecture \ref{conjecture1} is \emph{equivalent} to the following.

\begin{conjecture}[\cite{hertrich2021depth}]\label{conjecture2}
    For any $k \in \mathbb{N}$, $n= 2^k$, the function $f(x) = \max\{0, x_1, \ldots , x_n\}$ cannot be represented with $k$ hidden layers.
\end{conjecture}

\newchange{During the review process of this survey, 
Conjectures \ref{conjecture1} and \ref{conjecture2} were disproved by \cite{bakaev2025simplex}.

\begin{theorem}[\cite{bakaev2025simplex}]
Every piecewise linear function $\mathbb{R}^n \to \mathbb{R}$ can be represented by a ReLU DNN with at most $\lceil \log_3(n - 1)\rceil + 1$ hidden layers.
\end{theorem}

Other relevant recent works around these conjectures are
\cite{valerdi2024depth,averkov2025rational,grillo2025depth} and \cite{bakaev2025monotone}.
In what follows, we discuss prior work towards these conjectures, which produced notable results.
For example, \cite{hertrich2021depth} proved the following.}
\begin{theorem}[\cite{hertrich2021depth}]
    Under a technical assumption, there does not exist a 3-layer neural network that represents $f(x) = \max\{0, x_1, x_2, x_3, x_4\}$.
\end{theorem}

Later on, \newchange{the conjecture was shown to be true for the case of networks with integer weights}.

\begin{theorem}[\cite{haase2023lower}]
    Let $\text{ReLU}^\mathbb{Z}_n(k)$ be the set of functions $f:\mathbb{R}^n\to \mathbb{R}$ that can be represented using a network with $k$ hidden layers and integer weights.
    Then, for $n=2^k$, $\max\{0, x_1, \ldots , x_n\}\not\in \text{ReLU}^\mathbb{Z}_n(k)$.
    In particular, for every $k\in \{1,\ldots,\lceil \log_2(n+1) \rceil\} $
    \[\text{ReLU}^\mathbb{Z}_n(k-1)\subsetneq \text{ReLU}^\mathbb{Z}_n(k).\]
\end{theorem}

\change{
The proof of the last theorem relies on the concept of \emph{Newton polytopes}.
In what follows, we highlight some of the key features of the proof by \cite{haase2023lower}.
We remark that parts of this strategy---for instance, the use of Newton polytopes for understanding expressiveness---had already been proposed by \cite{hertrich2021depth} as potential avenues for proving Conjecture \ref{conjecture2},
along with other significant representability results.}

Consider $f$ a positively homogeneous, convex, piecewise linear function (such as the maximum function). Such a function can be written as $f(x) = \max\{a_1^\top x, \ldots, a_p^\top x\}$ for some family of vectors $\{a_i\}_{i=1}^p$; the Newton polytope $P_f$ is defined as $P_f:= \conv\{a_1,\ldots, a_p\}$. Note that if $f(x) = \max\{0, x_1, \ldots , x_n\}$ its Newton polytope is the unit simplex, i.e., $P_f = \Delta_0^n := \conv\{0,e_1, \ldots, e_n\}$ where $e_i$ is the $n$-th dimensional $i$-th unit vector. 

\cite{haase2023lower} \change{leverage a strong} connection between Newton polytopes and ReLU networks: that the $+$ and $\max$ operations over functions (which are the building blocks in ReLU networks) translate well to operations over polytopes, namely, $P_{f+g} = P_f + P_g$ and $P_{\max\{f,g\}} = \conv(P_f \cup P_g)$. Then, to connect the functions representable by networks with $k$ hidden layers and families of polytopes, the authors define the following recursive family.
\begin{align*}
\mathcal{P}'_{k+1} &:= \{\conv(Q_1 \cup Q_2)\, \mid \, Q_1,Q_2\in \mathcal{P}_k\} \\
\mathcal{P}_{k+1} &:= \{Q_1 +\cdots + Q_\ell\, \mid \, Q_1,\ldots, Q_\ell\in \mathcal{P}'_{k+1}\}.
\end{align*}
with $\mathcal{P}_0$ \change{being} the set of lattice points. 
\change{With this, they can show the following.}

\begin{theorem}[\cite{haase2023lower}]
   A positively homogeneous continuous piecewise linear function $f$ can be represented by an integral $k$ hidden-layer neural network if and only if $f = g-h$ for two convex, positively homogeneous continuous piecewise linear
functions $g$ and $h$ with $P_g, P_h \in \mathcal{P}_k$.
\end{theorem}

\change{We remark that a non-integral version of this theorem was originally proved by \cite{hertrich2021depth}, and later carried to the integral case in \cite{haase2023lower}.}

\change{\cite{haase2023lower}} use this theorem to show that, for $n=2^k$, there cannot be two polytopes $P,Q \in \mathcal{P}_k$  such that $P + \Delta_0^n = Q$. Therefore, $\max\{0, x_1, \ldots , x_n\} \not \in \text{ReLU}^\mathbb{Z}_n(k)$. 
The detailed arguments involve analyzing the parity of the volumes of the faces of polytopes in $\mathcal{P}_k$. We refer the reader to \cite{haase2023lower} for those.

\subsection{The number of linear regions}\label{sec:number}

We have previously seen conditions that affect the number of linear regions both positively and negatively. 
We discuss these and other analytical results in Section~\ref{sec:number_analytical}, and then discuss work on counting linear regions in practice in Section~\ref{sec:number_counting}.

\subsubsection{Analytical results}\label{sec:number_analytical}

At least three lines of work on analytical results have brought important insights. 
The first line is based on constructing networks with a large number of linear regions, 
which leads to lower bounds on the maximum number of linear regions. 
The second line is based on showing how the network architecture ---in particular its hyperparameters--- may impact the hyperplane arrangements defined by the layers, 
which leads to upper bounds on the maximum number of linear regions. 
The third line is based on characterizing the parameters of neural networks according to how they are initialized and updated along training, 
which leads to results on the expected number of linear regions for such networks. 

\paragraph{Lower bounds}

The lower bounds on the maximum number of linear regions are obtained through a careful choice of network parameters aimed at increasing the number of linear regions. 
In some cases, they also depend on particular choices of hyperparameters. 
We present them \change{in chronological order, which often coincides with nondecreasing values,} in Table~\ref{tab:lower_bounds}.

The first lower bound was introduced by~\cite{pascanu2013on} and then improved by those authors with a new construction technique in~\cite{montufar2014on}. 
In fact, Example~\ref{ex:zigzag} shows the case in which $n_0=1$ for the technique in~\cite{montufar2014on}. 
While a different construction is proposed by~\cite{Telgarsky2015}, 
subsequent developments in the literature have been based on~\cite{montufar2014on}. 

The lower bound by~\cite{arora2018understanding} is based on a different technique to construct a first wide layer based on zonotopes, 
which is then followed by the same layers as in~\cite{montufar2014on}.  
The first lower bound by~\cite{serra2018bounding} reflects a slight change to the technique used by~\cite{montufar2014on}, 
which in terms of Example~\ref{ex:zigzag} corresponds to using $n$ neurons to define $n+1$ instead of $n$ slopes on $[0,1]$. 
The second lower bound by~\cite{serra2018bounding} extends that of~\cite{arora2018understanding} by changing in the same way the construction of the subsequent layers of the network.

\begin{table}
\caption{Lower bounds on the maximum number of linear regions defined by a neural network.}
\label{tab:lower_bounds}
\centering
\vspace{2ex}
\begin{tabular}{@{\extracolsep{4pt}}cc}
\textbf{Reference} & \textbf{Bound and conditions}  \\
\cline{1-1}
\cline{2-2}
\noalign{\vskip2.5pt}
~\cite{pascanu2013on} & $\left(\prod\limits_{l=1}^{L-1} \left\lfloor \frac{n_l}{n_0} \right\rfloor \right) \sum\limits_{i=0}^{n_0} \binom{n_L}{i}$ \\
\cline{1-1}
\cline{2-2}
\noalign{\vskip4pt}
~\cite{montufar2014on} & $\left(\prod\limits_{l=1}^{L-1} \left\lfloor \frac{n_l}{n_0} \right\rfloor^{n_0}\right) \sum\limits_{i=0}^{n_0} \binom{n_L}{i}$, where $n_l \geq n_0 ~\forall l \in \sL$ \\
\cline{1-1}
\cline{2-2}
\noalign{\vskip4pt}
~\cite{Telgarsky2015} & $2^{\frac{L-3}{2}}$, where $n_i = 1$ for $i$ odd, $n_i = 2$ for $i$ even, and $L-3$ divides by 2 \\
\cline{1-1}
\cline{2-2}
\noalign{\vskip4pt}
~\cite{arora2018understanding} &  $2 \sum\limits_{j=0}^{n_0-1} \binom{m-1}{j} w^{L-1}$, where $2m=n_1$ and $w=n_l ~\forall l \in \sL \setminus \{1\}$ \\
\cline{1-1}
\cline{2-2}
\noalign{\vskip4pt}
~\cite{serra2018bounding} & $\left(\prod\limits_{l=1}^{L-1} \left( \left\lfloor \frac{n_l}{n_0} \right\rfloor + 1 \right)^{n_0} \right) \sum\limits_{i=0}^{n_0} \binom{n_L}{i}$, where $n_l \geq 3 n_0 ~\forall l \in \sL$ \\
\cline{1-1}
\cline{2-2}
\noalign{\vskip4pt}
~\cite{serra2018bounding} & $2 \sum\limits_{j=0}^{n_0-1} \binom{m-1}{j} (w+1)^{L-1}$, where $2m=n_1$ and $w=n_l \geq 3 n_0 ~\forall l \in \sL \setminus \{1\}$\\
\cline{1-1}
\cline{2-2}
\end{tabular}
\end{table}

\paragraph{Upper bounds}

The upper bounds on the maximum number of linear regions are obtained by primarily considering changes to the geometry of the linear regions from one layer to another, as previously outlined and revisited below. 
We present those with a close form 
\change{in chronological order, which often coincides with nonincreasing values,} in Table~\ref{tab:upper_bounds}. 

\cite{pascanu2013on} established the connection between linear regions and hyperplane arrangements, 
which leads to the tight bound for shallow networks based on~\cite{Zaslavsky1975} for activation hyperplanes in general position. 
\cite{montufar2014on} defined the first bound for deep networks based on enumerating all activation sets. 
The subsequent upper bounds extended the result by~\cite{pascanu2013on} to deep networks by considering its successive application through the sequence of layers of the network. 

In the case of \emph{deep} networks, where $L > 1$, 
we need to consider how the linear regions defined up to a given layer of the network can be further partitioned by the next layers. 
We 
start by assuming that every linear region defined by the first $l-1$ layers 
is then subdivided into the maximum possible number of linear regions defined by the activation hyperplanes of layer $l$. 
That leads to the bound 
in~\cite{raghu2017expressive}, 
which is implicit in their proof of an asymptotic bound of $O(n^{n_0 L})$, where $n$ is used as the width of every layer. 
However, there are many ways in which this bound can be refined upon careful examination. 
First, the dimension of the input of layer $l$ ---i.e., the output of layer $l-1$--- within each linear region is never larger than the smallest dimension among layers $1$ to $l$, since for every linear region we have an affine transformation between inputs and outputs of each layer \citep{montufar2017notes}.  
Second, the dimension of the input coming through each linear region is in fact bounded by the smallest number of active units in each of the previous layers \citep{serra2018bounding}.  
This leads to a tight upper bound for $n_0=1$, since it matches the lower bound in ~\cite{serra2018bounding}. 
Finally, 
the activation hyperplane of some units may not partition the linear regions because all possible inputs to the unit are in the same half-space, and in some of those cases the unit may never produce a positive output. 
For the number $k$ of active units in a given layer $l$, we can use the network parameters to calculate the maximum number of units that can be active in the next layer, $\mathcal{A}_l(k)$, as well as the number of units that can be active or inactive for different inputs, $\mathcal{I}_{l}(k)$ \citep{serra2020empirical}. 

\cite{hinz2019framework} observed that the upper bound by \cite{serra2018bounding} can be tightened by explicitly computing a recursive histogram of linear regions on the layers of the neural network according to the dimension of their image subspace. However, the resulting bound is not explicitly defined in terms of the network hyperparameters, and hence cannot be included on the table. This work is further extended in~\cite{hinz2021histograms} by also allowing a composition of bounds on subnetworks instead of only on the sequence of layers. 
Another extension of the framework from~\cite{hinz2019framework} by \cite{yutong2020framework} highlights that residual connections prevent the bottleneck effect in ResNets, 
by which reason such networks tend to have more linear regions. 

\cite{cai2023pruning} proposed a separate recursive bound based on \cite{serra2018bounding} to account for the sparsity of the weight matrices, 
which illustrates how pruning connections may affect the maximum number of linear regions.

The results above have also been extended to other architectures. 
In some cases, results on other types of activations are also part of the papers previously mentioned. 
\change{First, \cite{montufar2014on} and \cite{serra2018bounding} also present upper bounds for \emph{maxout} networks. 
In a more recent development, \cite{montufar2022maxout} present considerably tighter upper bounds for the number of linear regions in maxout networks with rank $k=3$ or greater. 
Second,} \cite{raghu2017expressive} present an upper bound for networks using \emph{hard tanh} activation. 
In other cases, the ideas discussed above have been adapted for sparser networks with parameter sharing: 
\cite{xiong2020cnn} present upper and lower bounds for convolutional networks, 
which are shown to asymptotically define more linear regions per parameter than rectifier networks with the same input size and number of layers. 
\cite{chen2022gcn} present upper and lower bounds for graph convolutional networks. 
\cite{matoba2022maxpooling} discuss the expresiveness of the maxpooling layers typically used in convolutional neural networks through their equivalence to a sequence of rectifier layers. 
Moreover, \cite{goujon2022role} present results for recently proposed activation functions, 
such as DeepSpline~\citep{agostineli2015spline,unser2019representer,bohra2020learning} and GroupSort~\citep{anil2019groupsort}. 

Some of the results above were also revisited through the lenses of tropical algebra, 
in which every linear region corresponds to a tropical hypersurface \citep{zhang2018tropical,charisopoulos2018tropical,maragos2021tropical,montufar2022maxout}.

In a sense, we can think of the discussion in Section~\ref{sec:representability} as a 
converse line of work exploring the minimum dimensions of a neural network capable of representing a given piecewise linear function. 
\change{Starting from the function to be represented rather than from the architecture, 
\cite{he2020finite} 
and 
\cite{hertrich2021depth} 
propose constructions with depth that is logarithmic on the input at the cost of a large width.
By relaxing the depth to be only logarithmic on the number of linear regions, 
\cite{chen2022bounds} propose a construction with smaller depth. 
More recently, 
\cite{brandenburg2025decomposition} propose a construction aimed at balancing depth and width.}  
On a related note, 
\cite{MPC} show that linear time-invariant systems in model predictive control can be exactly expressed by rectifier networks and provide bounds on the width and number of layers necessary for a given system, whereas \cite{ferlez2020aren} describe an algorithm for producing architectures that can be parameterized as an optimal model predictive control strategy.

\begin{table}
\caption{Upper bounds on the maximum number of linear regions defined by a neural network.}
\label{tab:upper_bounds}
\centering
\vspace{2ex}
\begin{tabular}{@{\extracolsep{4pt}}cc}
\textbf{Reference} & \textbf{Bound and conditions}  \\
\cline{1-1}
\cline{2-2}
\noalign{\vskip2.5pt}
~\cite{pascanu2013on} & $\sum\limits_{i=0}^{n_0} \binom{n_1}{n_0}$ for shallow networks, $n_1 \geq n_0$ \\
\cline{1-1}
\cline{2-2}
\noalign{\vskip4pt}
~\cite{montufar2014on} & $2^{\sum\limits_{l=1}^{L} n_l}$ \\
\cline{1-1}
\cline{2-2}
\noalign{\vskip4pt}
~\cite{raghu2017expressive} & $\prod\limits_{l=1}^{L} \sum\limits_{j=0}^{n_{l-1}} \binom{n_l}{j}$ \\
\cline{1-1}
\cline{2-2}
\noalign{\vskip4pt}
~\cite{montufar2017notes} & $\prod\limits_{l=1}^{L} \sum\limits_{j=0}^{d_l} \binom{n_l}{j}$, $d_l = \min\{n_0, n_1, \ldots, n_{l-1}\}$ \\
\cline{1-1}
\cline{2-2}
\noalign{\vskip4pt}
~\cite{serra2018bounding} & $\begin{array}{r}\sum\limits_{(j_1,\ldots,j_L) \in J} \prod\limits_{l=1}^L \binom{n_l}{j_l}, J = \{(j_1, \ldots, j_L) \in \mathbb{Z}^L: 0 \leq j_l \leq d_l ~\forall l \in \sL \},\\ d_l = \min\{n_0, n_1 - j_1, \ldots, n_{l-1} - j_{l-1}, n_l\}\ \end{array}$ \\
\cline{1-1}
\cline{2-2}
\noalign{\vskip4pt}
~\cite{serra2020empirical}  & $\begin{array}{r}\sum\limits_{(j_1,\ldots,j_L) \in J} \prod\limits_{l=1}^L \binom{\mathcal{I}_l(k_{l-1})}{j_l}, J = \{(j_1, \ldots, j_L) \in \mathbb{Z}^L: 0 \leq j_l \leq d_l, \\ d_l = \min\{n_0, k_1, \ldots, k_{l-1}, \mathcal{I}_l(k_{l-1})\},\\ k_0 = n_0, k_l = \mathcal{A}_{l}(k_{l-1}) - j_{l-1} ~\forall l \in \sL \}\end{array}$\\
\cline{1-1}
\cline{2-2}
\end{tabular}
\end{table}

\paragraph{Expected number}

The third analytical approach has been the evaluation of the expected number of linear regions. 
In a pair of papers, Hanin and Rolnick studied the number of linear regions based on how the network parameters are typically initialized.  
In the first paper \citep{hanin2019complexity}, 
they show that the average number of linear regions along 1-dimensional subspaces of the input grows linearly with respect to the number of neurons, irrespective of the network depth. 
In the second paper \citep{hanin2019deep}, 
they show that the average number of linear regions in higher-dimensional subspaces of the input also grows similarly in deep and shallow networks. 
For $N = \sum_{i=1}^L n_i$ as the total number of \change{neurons}, 
the expected number of linear regions is  $O(2^N)$ if $N \leq n_0$ and $O\left(\frac{(T N)^{n_0}}{n_0!}\right)$ otherwise, 
where $T > 0$ is a constant based on the network parameters. 
Moreover, some of their experiments suggest that the number of linear regions in shallow networks is slightly greater. 
According to the authors, 
these bounds reflect the fact that the family of functions that can be represented by neural networks in the way that they are typically initialized is considerably smaller. 
They further argue that training as currently performed is unlikely to expand the family of functions much further, as illustrated by their experiments. 
Similar results on the expected number of linear regions for maxout networks are presented by~\cite{tseran2021expected}, 
and an application of the results above results to data manifolds is explored by~\cite{tiwari2022manifolds}. 
Additional results for specific architectures of rectifier networks are conjectured by~\cite{wang2022estimation}, although without proof.

\subsubsection{Counting linear regions}\label{sec:number_counting}

Counting the actual number of linear regions of a given network has been a more challenging topic to explore. 
\cite{serra2018bounding} have shown that the linear regions of a trained network can be enumerated as the solutions of an MILP formulation, 
which has been slightly corrected in~\cite{cai2023pruning}\footnote{The MILP formulation of neural networks is discussed in Section~\ref{sec:optimizing}.}. 
However, MILP solutions are generally counted one by one \citep{danna2007multiple}, with exception of special cases \citep{serra2020nearoptimal} and small subproblems \citep{serra2020enumerative}, which makes this approach impractical for large neural networks. 
\newchange{In fact, 
\cite{stargalla2025counting} have shown that counting linear regions, as we define them in this survey, is \#P-complete.} 
\cite{serra2020empirical} have shown that approximate model counting methods, which are commonly used to count the number of feasible assignments in propositional satisfiability, can be easily adapted to solution counting in MILP, which leads to an order-of-magnitude speedup in comparison with exact counting. 
This type of approach is particularly suitable for obtaining probabilistic lower bounds, which can complement the analytical upper bounds for the maximum number of linear regions.  
In \cite{craighero2020compositional} and \cite{craighero2020understanding}, a directed acyclic graph is used to model the sets of active neurons on each layer and how they connected with those in subsequent layers. 
\cite{yang2020reachability} describe a method for decomposing the input space of rectifier networks into their linear regions by representing each linear region in terms of its face lattice, upon which the splitting operations corresponding to the transformations performed by each layer can be implemented. As the number of linear regions \change{grows}, these splitting operations can be processed in parallel. \cite{yang2021reachability} extend that method to convolutional neural networks. 
Moreover, \cite{wang2022estimation} describes an algorithm for enumerating linear regions that counts adjacent linear regions with same corresponding affine function as a single linear region. 

Another approach is to enumerate the linear regions in subspaces, which limits their number and reduces the complexity of the task. 
This idea was first explored by \cite{novak2018sensitivity} for measuring the complexity of a neural network in terms of the number of transitions along a single line. 
\cite{hanin2019complexity,hanin2019deep} also use this method with a bounded line segment or rectangle as a single set representing the input and then sequentially partitioning it. 
If this first set is intersected by the activation hyperplane of a neuron in the first layer, 
then we replace this set by two sets corresponding to the parts of the input space in which that neuron is active and not. 
Once those sets are further subdivided by all activation hyperplanes associated with the neurons in the first layer, 
the process can be continued with the neurons in the following layers. 
This method is used to count the number of linear regions along subspaces of the input with dimension 1 in \cite{hanin2019complexity} and dimension 2 in \cite{hanin2019deep}. 
A generalized version for counting the number of linear regions in affine subspaces spanned by a set of samples using an MILP formulation is presented in \cite{cai2023pruning}. 
An approximate approach for counting the number of linear regions along a line by computing the closest activation hyperplane in each layer is presented in \cite{gamba2022equal}.

Other approaches have obtained lower bounds on the number of linear regions of a trained network by limiting the enumeration or considering exclusively the inputs from the dataset. 
In \cite{xiong2020cnn}, the number of linear regions is estimated by sampling points from the input space and enumerating all activation patterns identified through this process. 
In \cite{cohan2022evolution}, the counting is restricted to the linear regions found between consecutive states of a neural network modeling a reinforcement learning policy.

\subsection{Applications and insights}

Thinking about neural networks in terms of linear regions led to a variety of applications. 
In turn, that inspired further studies on the structure and properties of linear regions under different settings. 
We organize the literature about applications and insights around some central themes in the subsections below.

\subsubsection{The number of linear regions}

From our discussion, the number of linear regions emerges as a potential proxy for the complexity of neural networks, 
which has been studied by some authors and exploited empirically by others. 
\cite{novak2018sensitivity} observed that the number of transitions between linear regions in 1-dimensional subspaces correlates with generalization. 
\cite{hu2020distillation} used bounds on the number of linear regions as proxy to model the capacity of a neural network used for learning through distillation, in which a smaller network is trained based on the outputs of another network. 
\cite{chen2021nas} and \cite{chen2021nas2} present one of the first approaches to training-free neural architectural search through the analysis of network properties. One of the two metrics that they have shown to be effective for that purpose is the number of linear regions associated with a sample of inputs from the training set on randomly initialized networks. 
\cite{biau2021wgans} observed that obtaining a discriminator network for Wasserstein GANs~\citep{arjovsky2017wgan} that correctly approximates the Wasserstein distance entails that such a discriminator network has a growing number of linear regions as the complexity of the data distribution increases. 
\cite{park2021unsupervised} maximized the number of linear regions in unsupervised learning in order to produce more expressive encodings for downstream tasks using simpler classifiers. 
In neural networks modeling reinforcement learning policies, 
\cite{cohan2022evolution} observed that the number of transitions between linear regions in inputs corresponding to consecutive states increases by 50\% with training while the number of repeated linear regions decreases. 
\cite{cai2023pruning} proposed a method for pruning different proportions of parameters from each layer by \change{maximizing} the bound on the number of linear regions, 
which leads to better accuracy than uniform pruning across layers. 
On a related note, 
\cite{liang2021brelu} proposed a new variant of the ReLU activation function for dividing the input space into a greater number of linear regions. 

The number of linear regions also inspired further theoretical work. 
\cite{amrami2021benefit} presented an argument for the benefit of depth in neural networks based on the number of linear regions for correctly separating samples associated with different classes. 
\cite{liu2021approximation} studied upper and lower bounds on the optimal approximation error of a convex univariate function based on the number of linear regions of a rectifier network. 
\cite{henriksen2022repairing} used the maximum number of linear regions as a metric for capacity that may limit repairing incorrect classifications in a neural network.

\subsubsection{The shape of linear regions}

Some studies aimed at understanding what affects the shape of linear regions in practice, including how to train neural networks in such a way to induce certain shapes in the linear regions. 
\cite{empirical2020iclr} observed that multiple training techniques may lead to similar accuracy, but very different shapes for the linear regions. 
For example, batch normalization~\citep{ioffe2015batchnorm} and dropout~\citep{srivastava2014dropout} lead to more linear regions. 
While batch normalization breaks the space in regions with uniform size, more orthogonal norms, and more gradient variability across adjacent regions; dropout produces more linear regions around decision boundaries, norms are more parallel, and data points less likely to be in the region containing the decision boundary. 
\cite{croce2019max} and \cite{LocallyLinear} applied regularization to the loss function to push the boundary of each linear region away from points in the training set that it contains, as long as those points are correctly classified. They show that this form of regularization improves the robustness of the neural network while making the linear regions larger.
In fact, \cite{zhu2020local} observed that the boundaries of the linear regions move away from the training data; 
and \cite{HashEncoders} conjectured that the linear regions near training samples becomes smaller through training, or that conversely the activation patterns are denser around the training samples. 
\cite{gamba2020biased} presented an empirical study on the angles between activation hyperplanes defined by convolutional layers, and observed that their cosines tend to be similar and more negative with depth after training.

The geometry of linear regions also led to other theoretical and algorithmic advances. 
Theoretically, \cite{phuong2020equivalence} proved that architectures with nonincreasing layer widths have unique parameters ---upon permutation and scaling--- for representing certain functions. 
In other words, some pairs of neural networks are only equivalent if their parameters only differ by permutation and multiplication. 
\cite{grigsby2023symmetries} showed that equivalences other than by permutation are less likely to occur with greater input size and width, but more likely with greater depth.
Algorithmically, 
\cite{ReverseEngineering} proposed a procedure to reconstruct a neural network by evaluating several inputs in order to determine regions of the input space 
for which the output of the neural network can be defined by an affine function ---and thus consist of a single linear region. 
Depending on how the shape changes between adjacent linear regions, 
the boundaries of the linear regions are replicated with neurons in the first hidden layer or in subsequent layers of the reconstructed neural network. 
\cite{masden2022combinatorial} provided theoretical results and an algorithm for characterizing the face lattice of the polyhedron associated with each linear region.

\subsubsection{Activation patterns and the discrimination of inputs}

Another common theme is understanding how inputs from the training and test sets are distributed among the linear regions, and what can be inferred from the encoding of the activation patterns associated with the linear regions.
\cite{gopinath2019property} noted that many properties of neural networks, including the classes of different inputs, are associated with activation patterns ---and thus with their linear regions. 
Several works~\citep{HashEncoders,sattelberg2020locally,tropex2021iclr} observed that each training sample is typically located in a different linear region when the neural network is sufficiently expressive; whereas
\cite{HashEncoders} noted that simple machine learning algorithms can be applied using the activation patterns as features, and 
\cite{sattelberg2020locally} noted that there is some similarity between activation patterns of different neural networks under affine mapping, meaning that the training of these neural networks lead to similar models.
\cite{chaudhry2020continual} exploited the idea of continual learning with different tasks being encoded in disjoint subspaces, 
which thus corresponds to a disjoint set of activation sets on each layer being associated with classifications for each of those tasks. 
Based on their approach for enumerating linear regions, 
\cite{craighero2020compositional} and \cite{craighero2020understanding} have found that inputs from larger linear regions are often correctly classified by the neural network, that inputs from smaller linear regions are often incorrectly classified, and that the number of distinct activations sets reduces along the layers of the neural network. 
\cite{gamba2022equal} also discussed the issue of some linear regions being smaller and thus less likely to occur in practice. 
Moreover, they propose a measurement for the similarity of the affine functions associated with linear regions along a line 
and observed that the linear regions tend to be less similar to one another when the network is trained with incorrectly classified labels. 

\change{Another important topic on the discrimination of inputs is determining if the function modeled by a neural network is injective and surjective, as recently studied by 
\cite{puthawala2022injective}
and \cite{froese2024injectivesurjective}.}
 
\subsubsection{Function approximation}

Because of the linear behavior of the output within each linear region, 
we can approximate the output of the neural network based on the output of its linear regions. 
\cite{chu2018pwnn} and \cite{sudjianto2020unwrapping} produced linear models based on this local behavior; whereas 
\cite{glass2021relumot} observed that we can interpret neural networks as equivalent to local linear model trees~\citep{nelles2000lolimot}, in which a distinct linear model is used at each leaf of a decision tree, and provided a method to produce such models from neural networks.
\cite{tropex2021iclr} described how to extract the linear regions associated with the inputs from the training set as means to approximate the output of the inputs from the test set. 
\cite{robinson2019dissecting} presented another approach for explicitly representing the function modeled by a neural network through the enumeration of its linear regions. 
On a related note, 
\cite{chaudhry2020continual} used the assumption of training samples remaining within the same linear region during gradient descent to simplify the analysis of backpropagation.

This topic also relates to the broad literature on neural networks as universal function approximators, to which the concept of linear regions helps articulating ideal conditions. 
As observed by \cite{mhaskar2020approximation}, the optimal number of linear regions in a neural network ---or, correspondingly, of pieces of the piecewise linear function modeled by it--- depends on the function being approximated.
In addition, linear regions were also used explicitly to build function approximations. 
\cite{kumar2019equivalent} have shown that rectifier networks can \change{be} approximated to arbitrary precision with two hidden layers, the largest of which having a neuron corresponding to each different activation pattern of the original network; an exact counterpart of this result \change{based on taking the limit of arbitrarily large parameters} was later presented by~\cite{villani2023shallow}. 
\newchange{A more recent development by \cite{zanotti2025depth} have shown that two hidden layers suffice to represent a piecewise linear function with a number of neurons linear on the number of linear regions.} 
\cite{fan2020quasiequivalence} described the transformation between sufficiently wide and deep networks while arguing that the fundamental measure of complexity should be counting simplices within linear regions. 
In subsequent work, \cite{fan2023simple} empirically observed that linear regions tend to have a small number of higher-dimensional faces, or facets.
\change{More recently, 
\cite{safran2024max} studied the number of neurons necessary to approximate the maximum function.}

Another line of study aimed at understanding the expressiveness and approximability of neural networks in terms of their number of parameters, 
in particular when the number of linear regions is greater than the number of parameters \citep{fractals2019neurips,fractals2020ieee,daubechies2022approximation}. 
They all discuss how the composition the modeled functions tend to present the self-similarity property of fractal distributions, 
which is one reason why they have so many linear regions. 
\cite{keup2022origami} interpreted the connection between linear regions in different parts of the input space in terms of how paper origamis are constructed: by ``folding'' the data for separability. 

Another related topic is computing the Lipschitz constant $\rho$ of the function $f(x)$ modeled by the neural network, 
the smallest $\rho$ for which $\| f(x') - f(x) \| \leq \rho \| x' - x \|$ for any two inputs $x$ and $x'$. 
Note that the first derivative of the output of a linear region is constant, 
which is leveraged by~\cite{hwang2020unrectifying} to evaluate the stability of the network by computing the constant across linear regions by changing the activation pattern. 
Interestingly, \cite{zhou2019lipschitz} showed that the constant grows similarly to the number of linear regions: polynomial in width and exponential in depth. 
A smaller constant limits the susceptibility of the network to adversarial examples~\citep{huster2018limitations}, which are discussed in Section~\ref{sec:optimizing}, 
and also lead to smaller bias variance \citep{loukas2021training}. 
While calculating the exact Lipschitz constant is NP-hard and encourages approximations \citep{scaman2018lipschitz,combettes2019certificates}, 
the exact constant can be computed using MILP \citep{jordan2020exact}. 
Notably, 
many studies have focused on relaxations such as linear programming \citep{zou2019cnnlp}, semidefinite programming \citep{fazlyab2019sdp,chen2020sdp}, and polynomial optimization \citep{latorre2020polynomial}. 
An alternative approach is to use more sophisticated activation functions for limiting the value of the constant \citep{anil2019groupsort,aziznejad2020controledlipschitz}. 

\change{In most of the discussion so far, 
we have seen the benefits of using neural networks for compactly representing---and computing---piecewise linear functions and piecewise linear approximations of arbitrary functions. 
We would like to acknowledge an interesting line for further work 
posited by one of the anonymous reviewers of this survey, 
and which would go in the opposite direction: 
what functions cannot be efficiently represented by a neural network, 
even if they can be efficiently computed by other means?}
\newchange{
For example, \cite{hertrich2024virtual} establish lower bounds on the size of certain types of neural networks that solve linear programming problems. 
}

\subsubsection{Optimizing over linear regions}

As an alternative to optimizing over neural networks as described next in Section~\ref{sec:optimizing}, 
a number of approaches have resorted to techniques that are equivalent to systematically enumerating or traversing linear regions and optimizing over them \citep{croce2018gcpr,croce2020ijcv,khedr2020verification,vincent2021icra,xu2022advml}.
Notably, \cite{vincent2021icra} and \cite{xu2022advml} are mindful of the facet-defining inequalities associated with a linear region, 
which are the ones to change when moving toward an adjacent linear region. 
On a related note, 
\cite{seck2021lp} alternates between gradient steps and solving a linear programming model within the current linear region.

\section{Optimizing Over a Trained Neural Network}\label{sec:optimizing}
In Section \ref{sec:training} we will see how polyhedral-based methods can be used to \emph{train} a neural network. In this section, we will focus on how polyhedral-based methods can be used to do something with a neural network \emph{after it has been trained.}
Specifically, after the network architecture and all parameters have been fixed, a neural network $f$ is merely a function. If each activation function $\sigma$ used to describe the network is piecewise linear (as is the case with those presented in Table~\ref{tab:activations}), $f$ is also a piecewise linear function. Therefore, any optimization problem \change{involving} $f$ in some way will be a piecewise linear optimization problem. For example, in the simple case where the output of $f$ is univariate, the optimization problem 
\[
    \min_{x \in \mathcal{X}} f(x)
\]
is a piecewise linear optimization problem. 
As discussed in Section \ref{sec:LR}, this problem can have an enormous number of ``pieces'' (linear regions) when $f$ is a neural network; solving this \change{problem} thus heavily depends on the size and structure of the neural network $f$. For example, the training procedure by which $f$ is obtained can greatly influence the performance of optimization strategies \citep{tjeng2017evaluating,xiao2018training}. 

In this section, we first explore situations in which \change{one} might want to optimize over a trained neural network in this manner. We will then survey available methods for solving this \change{problem (either exactly or deriving bounds)} using polyhedral-based methods. 
We conclude with a brief view of future directions. 

\subsection{Applications of optimization over trained networks}
Applications where \change{one} might want to optimize over a trained neural network $f$ broadly fall into two categories: those where $f$ is the ``true'' object of interest, and those where $f$ is a convenient proxy for some unknown, underlying behavior.

\subsubsection{Neural network verification} \label{sec:verification}
\begin{figure}
    \centering
    \begin{tikzpicture}
        \node[] at (0,0){\includegraphics[width=2.5cm]{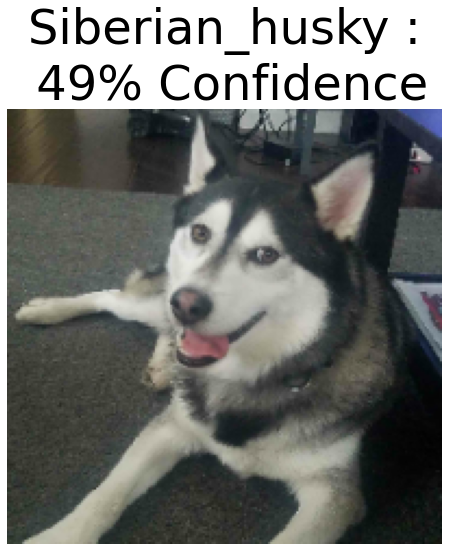}};
        \node[] at (1.9,-0.30){\Large \bf{+}};
        \node[] at (3.8,-0.26){\includegraphics[width=2.5cm]{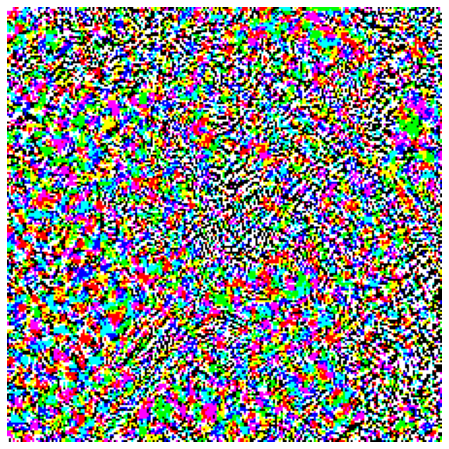}};
        \node[] at (5.7,-0.30){\Large \bf{=}};
        \node[] at (7.6,0){\includegraphics[width=2.5cm]{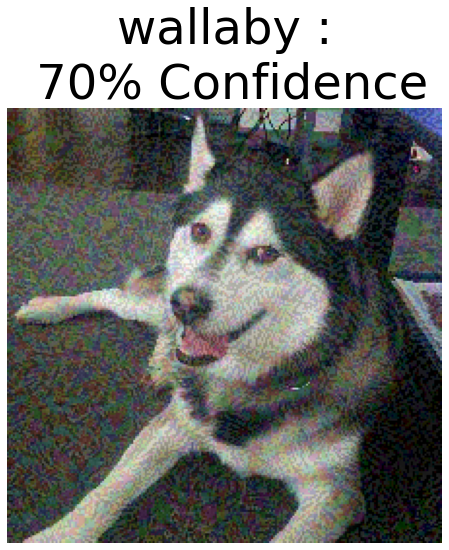}};
        \node[] at (3.8,1.2){$\times (\epsilon = 0.15)$};
    \end{tikzpicture}
    \caption{Example of adversarial attack on MobileNetV2 \citep{sandler2018mobilenetv2}. The original image taken by one of the survey authors is classified as `siberian\_husky,' but is re-classified as `wallaby' with a small (in an $\ell_\infty$-norm sense) targeted attack.}
    \label{fig:adversarial}
\end{figure}
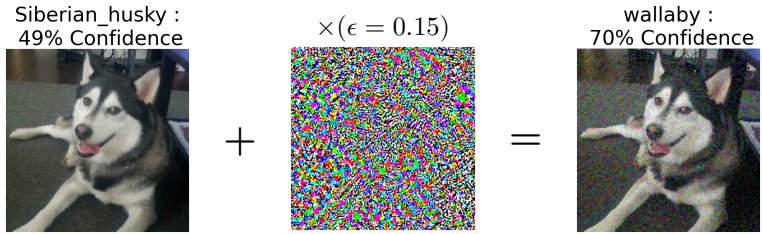

Neural network verification is a burgeoning field of study in deep learning. 
Starting in the early 2000s, researchers began to recognize the importance of rigorously verifying the behavior of neural networks, mainly in aviation-related applications \citep{schumann2003verification,zakrzewski2001verification}. 
More recently, the seminal works of \cite{szegedy2014intriguing} and \cite{goodfellow2015explaining} observed that neural networks are unusually susceptible to \emph{adversarial attacks}. These are small, targeted perturbations that can drastically affect the output of the network; as shown in Figure \ref{fig:adversarial}, even powerful models such as MobileNetV2 \citep{sandler2018mobilenetv2} are susceptible. The existence and prevalence of adversarial attacks in deep neural networks has raised justifiable concerns about the deployment of these models in mission-critical systems such as autonomous vehicles \citep{deng2020analysis}, aviation \citep{kouvaros2021formal}, or medical systems \citep{finlayson2019adversarial}. One fascinating empirical work by  \cite{eykholt2018robust} showed the susceptibility of standard image classification networks that might be used in self-driving vehicles to a very analogue form of attacks: black/white stickers, placed in a careful way, could confuse these models enough that they would mis-classify road signs (e.g., mistaking stop signs for ``speed limit 80'' signs).

Neural network verification seeks to prove (or disprove) a given input-output relationship, i.e., $x \in \mathcal{X} \Rightarrow y \in \mathcal{Y}$, that gives some indication of model robustness. 
Methods for verifying this relationship are classified as being sound and/or complete. 
A method that is \textit{sound} will only certify the relationship if it is indeed true (no false positives), while a method that is \textit{complete} will (i) always return an answer and (ii) only disprove the relationship if it is false (no false negatives). 
An early set of papers \citep{FischettiMIP,LomuscioMIP,tjeng2017evaluating} recognized that MILP provides an avenue for verification that is both sound and complete, given that $\mathcal{X}$ and $f(x)$ are both linear, or piecewise linear. 
We refer the readers to recent reviews \citep{huang2020survey,leofante2018automated,li2022sok,liu2021algorithms} for a more comprehensive treatment of the landscape of verification methods, including MILP- and LP-based technologies.

\begin{example}
Consider a classification network $f : [0,1]^{n_0} \to \mathbb{R}^d$ where the $j$-th output, $f(x)_j$, corresponds to the probability that input $x$ is of class $j$.\footnote{In actuality, we will instead typically work with the outputs corresponding to ``logits'', or unnormalized probabilities. These are typically fed into a softmax layer that normalizes them to correspond to a probability distribution over the classes. This nonlinear softmax transformation is not piecewise linear, but, thankfully, it can be omitted in the context of the verification task w.l.o.g.} Then consider a labeled image $\hat{x}$ known to be of class $i$, and a ``target'' adversarial class $k \neq i$. 
Then verifying local robustness of the prediction corresponds to checking $x \in \{ x: ||x-\hat{x}|| \leq \epsilon \} \Rightarrow y=f(x) \in \{ y: y_i \geq y_k \}$, where $\epsilon > 0$ is a constant which prescribes the radius around which $\hat{x}$ we will search for an adversarial example. 

This verification task can be formulated as an optimization problem of the form:
\begin{equation} \label{eq:verification}
\begin{aligned} 
    \max_{x \in [0,1]^{n_0}} \quad& f(x)_k - f(x)_i \\
    \text{s.t.}& ||x - \hat{x}|| \leq \epsilon.
\end{aligned}
\end{equation}Any feasible solution $x$ to this problem with positive cost is an adversarial example: it is very ``close'' to \change{$\hat{x}$} which has true label $i$, yet the network believes it is more likely to be of class $k$.\footnote{Alternative objectives may allow us to strengthen this statement to say that the network \emph{will} classify $x$ to be of class $k$. However, this will require a more complex reformulation to model this problem via MILP, so we omit it for simplicity.} If, on the other hand, it is proven that the optimal objective value is negative, this proves that $f$ is robust (at least in the neighborhood around $\hat{x}$). 
Note that the verification problem can terminate once the sign of the optimal objective value is determined, but solving the problem returns an optimal adversarial example. 
\end{example}

The objective function of \eqref{eq:verification} models the desired input-output relationship, $x \in \mathcal{X} \Rightarrow y \in \mathcal{Y}$, while the constraints model the domain $\mathcal{X}$. 
The domain $\mathcal{X}$ is typically a box or hyperrectangular domain. Extensions to this are described in Section \ref{sec:domains}.
Some explanation-focused verification applications define the input-output relationship in a derivative sense, e.g., $x \in \mathcal{X} \Rightarrow \partial y / \partial x \in \mathcal{Y}'$ \citep{wicker2022robust}. 
As the derivative of the ReLU function is also piecewise linear, this class of problems can also be modeled in MILP. For example, in the context of fairness and explainability, \cite{liu2020monotonic} and \cite{jordan2020exact} used MILP to certify monotonicity and to compute local Lipschitz constants, respectively.

Although in this survey we focus on optimization over trained neural networks, it is important to note that polyhedral theory underlies numerous strategies for neural network verification. 
For example, SAT and SMT (Satisfiability Modulo Theories) solvers designed for Boolean satisfiability problems (and more general problems for the case of SMT) can also be used to search through activation patterns for a neural network \citep{pulina2010abstraction}, resulting in tools that are sound and complete, such as Planet \citep{ehlers2017formal} and Reluplex \citep{katz2017reluplex}. 
\cite{bunel2018unified} presented a unified view to compare MILP and SMT formulations, as well as the relaxations that result from these formulations (we will revisit this in Section~\ref{sec:relaxations}).  
On the other hand, strategies such as ExactReach~\citep{xiang2017reachable} exploit polyhedral theory to compute reachable sets: given an input set to a ReLU function defined as a union of polytopes, the output reachable set is also a union of polytopes. 
Other methods over-approximate the reachable set to improve scalability, e.g., for vision models \citep{yang2021reachability}, resulting in methods that are sound, but not complete.

\subsubsection{Neural network as proxy}
Another situation in which \change{one} may want to solve an optimization problem containing trained neural networks is when \change{one} would like to optimize some other, unknown function for which historical input/output data exist. 
A similar situation arises when solving an optimization problem where (some of) the constraints are overly complicated, but samples from the constraints can be queried and used to train a simpler \textit{surrogate} model. 
In these cases, \change{a neural network can be trained} to approximate this underlying, unknown or complicated function. Then, since the neural network is known, the resulting optimization problem is a deterministic piecewise linear problem. 
Note that we focus here on using a neural network as a surrogate; neural networks can additionally learn other components of an optimization problem, e.g., uncertainty sets for robust optimization \citep{goerigk2023data} \newchange{or second-stage value functions for two-stage stochastic programming~\citep{alcantara2025quantile,patel2022neur2sp}}. 

Several software tools have been developed for this class of problems. For the case of constraint learning, \texttt{JANOS} \citep{bergman2022janos} and \texttt{OptiCL} \citep{maragno2021mixed} both provide functionality for learning a ReLU neural network to approximate a constraint based on data and embedding the learned neural network in MILP. 
The \texttt{reluMIP} package \citep{reluMIP} has also been introduced to handle the latter embedding step. 
More generally, \texttt{OMLT} \citep{ceccon2022omlt} translates neural networks to \texttt{pyomo} optimization blocks, including \newchange{various architectures (e.g,. CNNs, GNNs),} MILP formulations, and activation functions. 
\newchange{More recently, \texttt{PySCIPOpt-ML} \citep{turner2025pyscipopt} and \texttt{MathOptAI.jl} \citep{dowson2025mathoptai} provide similar interfaces for \texttt{SCIP} and \texttt{JuMP}, respectively.}  
Finally, recent developments in \texttt{gurobipy}\footnote{\url{https://github.com/Gurobi/gurobi-machinelearning}} enable directly parsing in trained neural networks.

Applications of this paradigm can be envisioned in a number of domain areas. This approach is common in deep reinforcement learning, where neural networks are used to approximate an unknown ``$Q$-function,'' or the long-term cost of taking a particular action at a particular state. In $Q$-learning, this $Q$-function is optimized iteratively to produce new candidate policies, which are then evaluated (typically via simulation) to produce new training data for future iterations. 
Optimization over the learned $Q$-function must be relatively fast in control applications, and several practical methods have been proposed. 
When the action space is discrete, the $Q$-function neural network is trained with one output value for each possible action, simplifying optimization to evaluating the model and selecting the largest output. 
Continuous action spaces require the $Q$ network be optimized over \citep{burtea2023safe,delarue2020reinforcement,ryu2020caql}, or an ``actor'' network can be trained to learn the optimal actions \citep{lillicrap2015continuous}. 
In a related vein, ReLU neural networks can be used as a process model for optimal scheduling or control \citep{wu2020scalable}.

Chemical engineering also presents applications where surrogate models have proven beneficial for optimization, as is the subject of recent reviews \citep{bhosekar2018advances,mcbride2019overview,tsay2019110th}.  
In particular, ReLU neural networks can be seamlessly embedded in larger MILP problems such as flow networks and reservoir control where the other constraints are also mixed-integer linear \citep{grimstad2019surrogate,Planning,yang2022modeling,liu2025diffconvex}. 
\newchange{Other recent uses for ReLU neural networks within MILP models are for representing the uncertainty in response variables in stochastic optimization~\citep{alcantara2023misp} and as a proxy for solving routing problems associated with a facility location problem~\citep{kaleem2024locationrouting}.}
Focusing on control applications where the neural network is embedded in a MILP that must be solved repeatedly, \cite{katz2020integrating} use multiparametric programming to learn the solution map of the resulting MILP itself, which is also piecewise affine. 
An emerging area of research uses verification tools to reason about neural networks used as controllers, e.g., see \citet{ARCH_COMP_20}. These applications involve optimization formulations combining the neural network with constraints defining the controlled system. 
For example, verification can be used to bound the reachable set \citep{sidrane2022overt,sosnin2024scaling} (alongside piecewise linear bounds on the dynamical system) or the maximum error against a baseline controller \citep{schwan2022stability}.

Finally, applications for optimization over neural networks arise in machine learning applications. 
MILP formulations can be used to compress neural networks \citep{serra2020lossless,serra2021compression,elaraby2020importance}, which consequently result in more tractable surrogate models \citep{kody2022modeling}. The main idea is to use MILP to identify \textit{stable} nodes, i.e., nodes that are always on or off over an input domain, which can then be algebraically eliminated. 
Optimization has also been employed in techniques for feature selection, based on identifying strongest input nodes \citep{sildir2022mixed,zhao2023model}. 
In the context of Bayesian optimization, \cite{volpp2020meta} use reinforcement learning to meta-learn acquisition functions parameterized as neural networks; selecting ensuing query points then requires optimization over the trained acquisition function. 
Later work modeled both the acquisition function and feasible region in black-box optimization as neural networks \citep{papalexopoulos22constrained}. 
In that work, exploration and exploitation are balanced via Thompson sampling and training multiple neural networks from a random parameter initialization. 
\newchange{Another application is the generation of counterfactual explanations for the neural network~\citep{kanamori2021counterfactual}.}

\paragraph{A word of caution}
Standard supervised learning algorithms aim to learn a function which fits the underlying function according to some distribution under which the data is generated. However, optimizing a function corresponds to evaluating it at a single point. This means that we may end up with a model that well-approximates the underlying function in distribution, but for which the pointwise minimizer is a poor approximation of the true function. This phenomenon is referred to as the ``Optimizer's curse'' \citep{smith2006}.

\subsection{\change{Preliminaries}}
\change{ We now mathematically outline our focus on single neuron formulations and detail the scope of this survey.}

\subsubsection{Single neuron relaxations}

For the following subsections, consider the $i$-th neuron in the $l$-th layer of a neural network, endowed with a ReLU activation function, whose behavior is governed by \eqref{eq:single-neuron}. Presume that \change{the input domain} of interest $\mathcal{D}^{l-1} \subset \mathbb{R}^{n_l}$ is a bounded region. Further, since $\mathcal{D}^{l-1}$ is bounded, presume that finite bounds are known on each input component, i.e. that vectors $L^{l-1},U^{l-1} \in \mathbb{R}^{n_l}$ are known such that $\mathcal{D}^{l-1} \subseteq [L^{l-1},U^{l-1}] \subset \mathbb{R}^{n_l}$. We can then write the \emph{graph} of the neuron, which couples together the input and the output of the nonlinear ReLU activation function:
\begin{align*}
    \gr &= \Set{ (\vh^{l-1},h^l_i) \in \mathcal{D}^{l-1} \times \mathbb{R} | h^l_i = 0 \geq \vw^l_i \vh^{l-1} + b^l_i} \\
    &\cup \Set{ (\vh^{l-1},h^l_i) \in \mathcal{D}^{l-1} \times \mathbb{R} | h^l_i = \vw^l_i \vh^{l-1} + b^l_i \geq 0}.
\end{align*}
This is a disjunctive representation for $\gr$ in terms of two polyhedral alternatives.
We assume that every included neuron exhibits this disjunction, i.e., every neuron can be on or off depending on the model input. 
This assumption of \emph{strict activity} implies that \change{the sign of the preactivation term $(\vw^l_i \vh^{l-1} + b^l_i)$ is unknown}, noting that neurons not satisfying this property can be exactly pruned from the model \citep{serra2020lossless}. 
\change{Specifically, if a neuron is strictly inactive, it can be directly removed; if a neuron is strictly active, it can be replaced with a skip connection, or the weights of downstream connections updated accordingly.}

We observe that, given this (or any) formulation for each individual unit, it is straightforward to construct a formulation for the entire network. For example, if we take $X^l_i = \Set{(\vh^{l-1},h^l_i,z^l_i) | \eqref{eqn:relu-big-m} }$ for each layer $l$ and each unit $i$, we can construct a MILP formulation for the graph of the entire network, $\Set{(x,f(x)) : x \in \mathcal{D}^0}$:
\[
    (\vh^{l-1},h^l_i,z^l_i) \in X^l_i \quad \forall l \in \sL, i \in \llbracket n_l \rrbracket.
 \]  

This also generalizes in a straightforward manner to more complex feedforward network architectures (e.g. convolutions, or sparse or skip connections), though we omit the explicit description for notational simplicity. 
\newchange{One interesting line of work here explores formulations specific to graph neural networks (GNNs), using MILP techniques to simultaneously handle the GNN and the combinatorial input search space over graphs~\citep{mcdonald2024mixed,zhang2023optimizing}.}

\subsubsection{Beyond the scope of this survey}

The effectiveness of the single-neuron formulations described above is bounded by the tightness of the optimal univariate formulation; this property is known as the ``single-neuron barrier'' \citep{salman2019convex}. 
This has motivated research in convex relaxations that jointly account for multiple neurons within a layer \citep{singh2019beyond}. 
Nevertheless, the analysis of polyhedral formulations for multiple neurons simultaneously quickly becomes intractable, and is beyond the scope of this survey. Instead, we point the interested reader to the recent survey by \cite{roth2021primer}, and highlight a few approaches taken in the literature. Multi-neuron analysis has been used to: improve bounds tightening schemes \citep{rossig2021advances}, prune linearizable neurons \citep{botoeva2020efficient}, design dual decompositions \citep{ferrari2022complete}, and generate strengthening inequalities \citep{serra2020empirical, zhou2024scalable}. 
Similarly, we do not review formulations for ensembles of ReLU networks, though MILP formulations have been proposed \citep{wang2021ensemble,wang2023optimizing}. 

Additionally, recent works have exploited polyhedral structure to develop sampling-based strategies, which can be used to \change{warm-start MILP (i.e., primal heuristics)} or accelerate local search in verification \citep{perakis2022optimizing,tong2024optimization,wu2022efficient}. 
\newchange{Other recent approaches include the use of surrogate models based on pruned neural networks~\citep{pham2025surrogate} and dualizing the constraints associated with ReLU in the optimization model~\citep{liu2025diffconvex}.}
\cite{lombardi2017empirical} computationally compare MILP against local search and constraint programming approaches.  
\newchange{\cite{plate2025analysis} analyze the impact of a variety of factors on the runtimes of MILP with embedded neural networks.} 
In a related vein, \cite{cheon2022outer} examines local solutions and proposes an outer approximation method to improve gradient-based optimization. 
Finally, following \cite{raghunathan2018semidefinite}, a large body of work has presented optimization-based methods for verification that use semidefinite programming concepts \citep{dathathri2020enabling,fazlyab2020safety,newton2021exploiting}.
Notably, \cite{batten2021efficient} showed how combining semidefinite and MILP formulations can produce a new formulation that is tighter than both. This was later extended with reformulation-linearization technique, or RLT, cuts \citep{lan2022tight}. 
While related to linear programming and other methods based on convex relaxations, this stream of work is beyond the scope of this survey. 
We refer the reader to \cite{zhang2020tightness} for a discussion of the tightness of these formulations. 
\newchange{More recently, \cite{karia2025kan} studied the embedding of Kolmogorov-Arnold networks~\citep{liu2025kan}, which may have activation functions requiring a nonlinear optimization formulation.}

\subsection{Exact models using mixed-integer programming}
\label{sec:MIPmodels}

Mixed-integer programming offers a powerful algorithmic framework for \emph{exactly} modeling nonconvex piecewise linear functions. The Operations Research community has a \change{long and storied} history of developing MILP-based methods for piecewise linear optimization, with research spanning decades \citep{croxton2003comparison,dantzig1960significance,geissler2012using,huchette2022nonconvex,lee2001polyhedral,misener2012global,padberg2000approximating,vielma2010mixed}.
However, many of these techniques are specialized for low-dimensional or separable piecewise linear functions. While a reasonable assumption in many OR problems, this is not the case when modeling neurons in a neural network. Therefore, the standard approach in the literature is to apply general-purpose MILP formulation techniques to model neural networks.

\paragraph{Connection to Boolean satisfiability}
Some SMT-based methods such as Reluplex \citep{katz2017reluplex} and Planet \citep{ehlers2017formal} effectively construct branching technologies similar to MILP solvers. 
Indeed, \texttt{Marabou} \citep{katz2019marabou} builds on Reluplex, and a recent extension \texttt{MarabouOpt} can optimize over trained neural networks \citep{strong2021global}. 
The authors also outline general procedures to extend verification solvers to optimization. 
Our focus in this review is on more general MILP formulations, or those that can be incorporated into off-the-shelf MILP solvers with relative ease. 
\cite{bunel2020branch,bunel2018unified} provide a more comprehensive discussion of similarities and differences to SMT. 

\subsubsection{The big-$M$ formulation} 

The big-$M$ method is a standard technique used to formulate logic and disjunctive constraints using mixed-integer programming \citep{bonami2015mathematical,vielma2015mixed}. Big-$M$ formulations are typically very simple to reason about and implement, and are quite compact, though their convex relaxations can often be quite poor, leading to weak dual bounds and (often) slow convergence when passed to a mixed-integer programming solver. Since $\gr$ is a disjunctive set, the big-$M$ technique can be applied to produce the following formulation:
\begin{subequations} \label{eqn:relu-big-m}
\begin{align}
    h^l_i &\geq \vw^l_i \vh^{l-1} + b^l_i \label{eqn:relu-big-m-1} \\
    h^l_i &\leq \left(\vw^l_i \vh^{l-1} + b^l_i\right) - M^{l}_{i,-}(1-z) \\
    h^l_i &\leq M^{l}_{i,+}z \\
    (\vh^{l-1},h^l_i) &\in [L^{l-1},U^{l-1}] \times \mathbb{R}_{\geq 0} \label{eqn:relu-big-m-4} \\
    z^l_i &\in \{0,1\}. \label{eqn:relu-big-m-5}
\end{align}
\end{subequations}
Here, $M^l_{i,-}$ and $M^l_{i,+}$ are data which must satisfy the inequalities
\begin{align*}
    M^l_{i,-} &\leq \min_{\vh^{l-1} \in \mathcal{D}^{l-1}} \vw^l_i \vh^{l-1} + b^l_i \\
    M^l_{i,+} &\geq \max_{\vh^{l-1} \in \mathcal{D}^{l-1}} \vw^l_i \vh^{l-1} + b^l_i.
\end{align*}
This big-$M$ formulation for ReLU-based networks has been used extensively in the literature \citep{bunel2018unified,Cheng2017,DuttaMIP,FischettiMIP,kumar2019equivalent,LomuscioMIP,serra2020empirical,serra2018bounding,tjeng2017evaluating,xiao2018training}.

The big-$M$ formulation is compact, with one binary variable and $\mathcal{O}(1)$ general inequality constraints for each neuron. 
Applied for each unit in the network, this leads to a MILP formulation with $\mathcal{O}(\sum_{l \in \sL} n_l) = \mathcal{O}(Ln_{\max} )$ binary variables and general inequality constraints, where $n_{\max} = \max_{l \in \sL} n_L$. However, it has been observed \citep{anderson2019strong,anderson2020strong} that this big-$M$ formulation is not strong in the sense that its LP relaxation does not, in general, capture the convex hull of the graph of a given unit; see Figure~\ref{fig:relu-neuron} for an illustration. In fact, this LP relaxation can be arbitrarily bad \citep[Example 2]{anderson2019strong}, even in fixed input dimension.
As MILP solvers often bound the objective function between the best feasible point and its tightest optimal continuous relaxation, a weak formulation can negatively impact performance, often substantially. 

It is worth dwelling on where this lack of strength comes from. If the input $\vh^{l-1}$ is one dimensional, the big-$M$ formulation is \emph{locally} ideal \citep{vielma2015mixed}: the extreme points of the linear programming relaxation (\ref{eqn:relu-big-m-1}-\ref{eqn:relu-big-m-4}) naturally satisfy the integrality constraints \eqref{eqn:relu-big-m-5}. However, this fails to hold in the general multivariate input case. To see why, observe that the bounds on the input variables $\vh^{l-1}$ are only coupled with the logic involving the binary variable $z$ only in an aggregated sense, through the coefficients $M^l_{i,-}$ and $M^l_{i,+}$. In other words, the ``shape'' of the pre-activation domain is not incorporated directly into the big-$M$ formulation. 
Furthermore, the strength of this formulation highly depends on the big-$M$ coefficients. 
These coefficients can be obtained using techniques ranging from basic interval arithmetic to optimization-based bounds tightening. 
\cite{grimstad2019surrogate} show how constraints external to the neural network can yield tighter bounds via optimization- or feasibility-based bounds tightening. 
\cite{rossig2021advances} compare several methods for deriving bounds and further develop optimization-based bounds tightening based on pairwise dependencies between variables. 
\newchange{\cite{zhou2024scalable} present a rolling-horizon decomposition for bounds tightening based on solving a subproblem for each layer of the network.} 

\tikzstyle{yzx} = [
  x={(1.2*.9625cm, 1.2*.9625cm)},
  y={(1.2*2.5cm, 0cm)},
  z={(0cm, 1.2*3cm)},
]

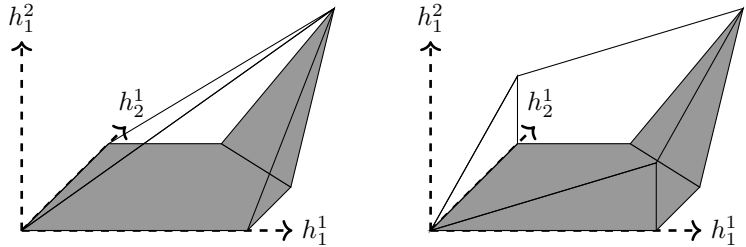
\begin{figure}[t]
    \centering
    \begin{tikzpicture}[yzx]
        \draw [->, dashed, line width=1] (0,0,0) -- (1.2,0,0);
        \draw [->, dashed, line width=1] (0,0,0) -- (0,1.2,0);
        \draw [->, dashed, line width=1] (0,0,0) -- (0,0,0.7);
        \node[above right] at (1.2,-0.075,0) {$h^1_2$};
        \node[right] at (0,1.2,0) {$h^1_1$};
        \node[above] at (0,0,.7) {$h^2_1$};
        \coordinate (LL) at (0,0,0);
        \coordinate (UL) at (1,0,0);
        \coordinate (LU) at (0,1,0);
        \coordinate (UU) at (1,1,0.5);
        \coordinate (UM) at (1,0.5,0);
        \coordinate (MU) at (0.5,1,0);
        
        \draw [fill=gray!80] (LL) -- (UL) -- (UM) -- (MU) -- (LU) -- cycle;
        \draw [fill=gray!80] (UU) -- (UM) -- (MU) -- cycle;
        \draw (LL) -- (UL) -- (UU) -- cycle;
        \draw (LL) -- (LU) -- (UU) -- cycle;
    \end{tikzpicture}  \hspace{2em}
    \begin{tikzpicture}[yzx]
        \draw [->, dashed, line width=1] (0,0,0) -- (1.2,0,0);
        \draw [->, dashed, line width=1] (0,0,0) -- (0,1.2,0);
        \draw [->, dashed, line width=1] (0,0,0) -- (0,0,0.7);
        \node[above right] at (1.2,-0.075,0) {$h^1_2$};
        \node[right] at (0,1.2,0) {$h^1_1$};
        \node[above] at (0,0,.7) {$h^2_1$};
        \coordinate (LL) at (0,0,0);
        \coordinate (UL) at (1,0,0);
        \coordinate (LU) at (0,1,0);
        \coordinate (UU) at (1,1,0.5);
        \coordinate (UM) at (1,0.5,0);
        \coordinate (MU) at (0.5,1,0);
        \coordinate (E1) at (1,0,1/4);
        \coordinate (E2) at (0,1,1/4);

        \draw [fill=gray!80] (LL) -- (UL) -- (UM) -- (MU) -- (LU) -- cycle;
        \draw [fill=gray!80] (UU) -- (UM) -- (MU) -- cycle;
        \draw (LL) -- (E2) -- (UU) -- (E1) -- cycle;
        \draw (LL) -- (UL) -- (E1) -- cycle;
        \draw (LL) -- (LU) -- (E2) -- cycle;
        
    \end{tikzpicture}
    \caption{\textbf{Left:} The convex hull of a ReLU neuron \eqref{eqn:balas-relu}, and \textbf{Right:} the convex relaxation offered by the big-$M$ formulation \eqref{eqn:relu-big-m} Adapted from Anderson et al. \cite{anderson2020strong,anderson2019strong}}
    \label{fig:relu-neuron}
\end{figure}

\subsubsection{A stronger extended formulation} 
A much stronger MILP formulation can be constructed through a classical method, the extended formulation for disjunctions \citep{balas1998disjunctive,jeroslow1984modelling}. This formulation for a given ReLU neuron takes the following form~\citep[Section 2.2]{anderson2019strong}:
\begin{subequations} \label{eqn:balas-relu}
\begin{align}
    (\vh^{l-1},h^l_i) &= (x^+,y^+) + (x^-,y^-) \label{eqn:balas-relu-1} \\
    y^- &= 0 \geq \vw^{l}_i x^- + b^l_i(1-z) \\
    y^+ &= \vw^{l}_i x^+ + b^l_iz \geq 0 \\
    L^{l-1}(1-z) &\leq x^- \leq U^{l-1}(1-z) \\
    L^{l-1}z &\leq x^+ \leq U^{l-1}z \label{eqn:balas-relu-5} \\
    z &\in \{0,1\}. \label{eqn:balas-relu-6}
\end{align}
\end{subequations}
This formulation requires one binary variable and $\mathcal{O}(n_{l-1})$ general linear constraints and auxiliary continuous variables. It is also locally ideal, i.e., as strong as possible. While the number of variables and constraints for an individual unit seems quite tame, applying this formulation for unit in a network leads to a formulation with $\mathcal{O}(n_0 + \sum_{l \in \sL} n_ln_{l-1}) = \mathcal{O}(|\sL|n_{\max}^2)$ continuous variables and linear constraints. Moreover, while the formulation for \emph{an individual unit} is locally ideal, the composition of many locally ideal formulations will, in general, fail to be ideal itself. 
Consider that, while each node can be modeled as a two-part disjunction, the full network requires exponentially many disjuncts, each corresponding to one activation pattern. 

Despite its strength and relatively modest increase in size relative to the big-$M$ formulation \eqref{eqn:relu-big-m}, it has been empirically observed that this formulation often performs worse than expected \citep{anderson2019strong,vielma2019small}, both in the verification setting and more broadly.

\subsubsection{A class of intermediate formulations}
The previous sections observed that the big-$M$ formulation \eqref{eqn:relu-big-m} is compact, but may offer a weak convex relaxation, while the extended formulation \eqref{eqn:balas-relu} offers the tightest possible convex relaxation for an individual unit, at the expense of a much larger formulation.
\cite{kronqvist2022psplit,kronqvist2021steps} present a strategy for obtaining formulations intermediate to \eqref{eqn:relu-big-m} and \eqref{eqn:balas-relu} in terms of both size and strength. 
The key idea is to partition $\vw_i^l \vh^{l-1}$ into a number of aggregated variables, $\vw_i^l \vh^{l-1} = \sum_{p=1}^P \hat{x}_p$. 
Each auxiliary variable $\hat{x}_p$ is defined as a sum of a subset of the $j$-th weighted inputs $\hat{x}_p = \sum_{j \in \mathbb{S}_p} w_{i,j}^l h_j^{l-1}$, with $\mathbb{S}_1, ..., \mathbb{S}_P$ partitioning $\{1,...,n_{l-1}\}$.  
This technique can be applied to the ReLU function, giving the convex hull over the directions defined by $\hat{x}_p$ \citep{tsay2021partition}: 
\begin{subequations} \label{eq:relu-partition}
\begin{align}
\left( \sum_{j \in \mathbb{S}_p} w_{i,j}^l h_j^{l-1},h^l_i \right) &= (\hat{x}_p^+,y^+) + (\hat{x}_p^-,y^-) \label{eq:Pextended1} \\
    y^- &= 0 \geq \sum_p \hat{x}_p^- + b^l_i(1-z) \\
    y^+ &= \sum_p \hat{x}_p^+ + b^l_iz \geq 0 \\
    \hat{\boldsymbol{M}}_{i,-}^l(1-z) &\leq \hat{x}^- \leq \hat{\boldsymbol{M}}_{i,+}^l(1-z) \\
    \hat{\boldsymbol{M}}_{i,-}^l z &\leq \hat{x}^+ \leq \hat{\boldsymbol{M}}_{i,+}^l z \\
    z &\in \{0,1\}. \label{eq:Pextended_1end}
\end{align}
\end{subequations}

Here, the $p$-th elements of $\hat{\boldsymbol{M}}_{i,-}^l$ and $\hat{\boldsymbol{M}}_{i,+}^l$ must satisfy the inequalities
\begin{align*}
    \hat{M}^l_{i,-,p} &\leq \min_{\vh^{l-1} \in \mathcal{D}^{l-1}} \sum_{j \in \mathbb{S}_p} w_{i,j}^l h_j^{l-1} \\
    \hat{M}^l_{i,+,p} &\geq \max_{\vh^{l-1} \in \mathcal{D}^{l-1}} \sum_{j \in \mathbb{S}_p} w_{i,j}^l h_j^{l-1}.
\end{align*}

These coefficients can be derived using techniques analogous to those for the big-$M$ formulation (note that tighter bounds may be derived by considering $\hat{x}^-$ and $\hat{x}^+$ separately). 
Observe that when $P=1$, we recover the same tightness as the big-$M$ formulation \eqref{eqn:relu-big-m}, as, intuitively, the formulation is built over a single ``direction'' corresponding to $\vw_i^l \vh^{l-1}$. Conversely, when $P=n_{l-1}$, we recover the tightness of the extended formulation \eqref{eqn:balas-relu}, as each direction corresponds to a single element of $\vh^{l-1}$.
\cite{tsay2021partition} study partitioning strategies and show that intermediate values of $P$ result in formulations that can outperform the two extremes, by balancing formulation size and strength.

\subsubsection{Cutting plane methods: Trading variables for inequalities} \label{sec:strengthening-inequalities}
The extended formulation \eqref{eqn:balas-relu} achieves its strength through the introduction of auxiliary continuous variables. However, it is possible to produce a formulation of equal strength by projecting out these auxiliary variables, leaving an ideal formulation in the ``original'' $(\vh^{l-1},h^l_i,z)$ variable space. While in general this projection may be difficult computationally, for the simple structure of a single ReLU neuron it is possible to characterize in closed form. The formulation is given by \cite{anderson2020strong,anderson2019strong} as
\begin{subequations} \label{eqn:relu-ideal}
\begin{align}
    h^l_i &\geq \vw_i^l \vh^{l-1} + b^l_i \label{eqn:relu-ideal-1} \\
    h^l_i &\leq \sum_{j \in J} w^{l}_{i,j}(h^{l-1}_i - \breve{L}^{l}_j(1-z))  + \left(b + \sum_{j \not\in J} w^{l}_{i,j}\breve{U}_j \right)z \quad \forall J \subseteq \llbracket n_{l-1} \rrbracket \label{eqn:relu-ideal-2} \\
    (\vh^{l-1},h^l_i) &\in \mathcal{D}^{l-1} \times \mathbb{R}_{\geq 0} \label{eqn:relu-ideal-3} \\
    z^l_i &\in \{0,1\},
\end{align}
\end{subequations}
where notationally, for each $j \in \llbracket n_{l-1} \rrbracket$, we take 
\begin{align*}
\breve{L}^{l-1}_j &= \begin{cases} L^{l-1}_j & w^{l}_{i,j} \geq 0 \\ U^{l-1}_j & w^{l}_{i,j} < 0 \end{cases} \quad \mathrm{and} \quad 
\breve{U}^{l-1}_j = \begin{cases} U^{l-1}_j & w^{l}_{i,j} \geq 0 \\ L^{l-1}_j & w^{l}_{i,j} < 0 \end{cases} 
\end{align*}

We note a few points of interest about this formulation. First, it is ideal, and so recovers the convex hull of a ReLU activation function, coupled with its preactivation affine function and bounds on each of the inputs to that affine function. Second, it can be shown that, under very mild conditions, each of the exponentially many constraints in \eqref{eqn:relu-ideal-2} are necessary to ensure this property; none are redundant and can be removed without affecting the relaxation quality. Third, note that by selecting only those constraints in \eqref{eqn:relu-ideal-2} corresponding to $J = \emptyset$ and $J = \llbracket n_{l-1} \rrbracket$, we recover the big-$M$ formulation \eqref{eqn:relu-big-m} in the case where $\mathcal{D}^{l-1} = [L^{l-1},U^{l-1}]$. This suggests a practical approach for using this large family of inequalities: Start with the big-$M$ formulation, and then dynamically generate violated inequalities from \eqref{eqn:relu-ideal-2} as-needed in a cutting plane procedure. As shown by \cite{anderson2020strong}, this separation problem is separable in the input variables, and hence can be completed in $\mathcal{O}(n_{l-1})$ time. 

The cutting plane strategy is in general compatible with weaker formulations, such as relaxation-based verification \citep{zhang2022general} and formulations from the class \eqref{eq:relu-partition}. 
In fact, \cite{tsay2021partition} show that the intermediate formulations in \eqref{eq:relu-partition} effectively pre-select a number of inequalities from \eqref{eqn:relu-ideal-2}, in terms of their continuous relaxations. 
While adding these constraints results in a tighter continuous relaxation, the added constraints can eventually significantly increase the model size. Practical implementations may therefore only perform cut generation at a limited number of branch-and-bound search nodes \citep{depalma2021scaling,tsay2021partition}. 

\paragraph{A subtlety when using \eqref{eqn:relu-ideal}}
This third point above raises a subtlety discussed in the literature~\citep[Appendix F]{depalma2021scaling}. Often, additional structural information is known about $\mathcal{D}^{l-1}$ beyond bounds on the variables. In this case, it is typically possible to derive tighter values for the big-$M$ coefficients. In this case, when using a separation-based approach it is preferable to initialize the formulation with these tightened big-$M$ constraints, and then proceed with the cutting plane approach as normal from there.

\subsection{Scaling further: Convex relaxations and linear programming}
\label{sec:relaxations}

The above demonstrate MILP as a powerful framework for exactly modeling complex, nonconvex trained neural networks, but standard solvers are often not sufficiently scalable to adequately handle large-scale networks. A natural approach to increase the scalability, then, is to \emph{relax} the network in some manner, and then apply convex optimization methods. For the verification problem discussed in Section~\ref{sec:verification}, this yields what is known as an \emph{incomplete verifier}: any certification of robustness provided can be trusted (no false positives), but there may be robust instances that the method cannot prove are \change{indeed robust} (some false negatives). 
In other words, \change{using only a relaxation produces a verifier that is sound, but not complete}. 

While a variety of methods exist for accomplishing this, in this section, we briefly outline techniques relevant to polyhedral theory. 
In particular, we focus on some techniques for building convex polyhedral relaxations. 
The most natural convex relaxation for a MILP formulation is its linear programming (LP) relaxation, constructed by dropping any integrality constraints. For example, the LP relaxation of \eqref{eqn:relu-big-m} is given by the system (\ref{eqn:relu-big-m-1}-\ref{eqn:relu-big-m-4}). This is a compact linear programming relaxation for a ReLU-based network, and is the basis for methods due to \cite{bunel2020lagrangian} and \cite{ehlers2017formal}.

\subsubsection{Projecting the big-$M$ and ideal MILP formulations}
This section examines projections of the linear relaxations of formulations \eqref{eqn:relu-big-m} and \eqref{eqn:relu-ideal}. \\

\textbf{(Projecting the big-$M$).}
Note that the LP relaxation given by (\ref{eqn:relu-big-m-1}--\ref{eqn:relu-big-m-4}) maintains the variables $z^l_i$ in the formulation, though they are no longer required to satisfy integrality. Since these variables are ``auxiliary'' and are no longer necessary to encode the nonconvexity of the problem, they can be projected out without altering the quality of the convex relaxation. Doing this yields what is commonly known as the ``triangle'' or ``$\Delta$'' relaxation \citep{salman2019convex}:
\begin{subequations} \label{eqn:triangle-relaxation}
\begin{align}
    h^l_i &\geq \vw^l_i \vh^{l-1} + b^l_i \\
    h^l_i &\leq \frac{M^l_{i,+}}{M^l_{i,+} - M^l_{i,-}}(\vw^l_i \vh^{l-1} + b^l_i) \\
    (\vh^{l-1},h^l_i) &\in [L^{l-1},U^{l-1}] \times \mathbb{R}_{\geq 0}.
\end{align}
\end{subequations}

While the LP relaxation \eqref{eqn:triangle-relaxation} for an individual unit is compact, modern neural network architectures regularly comprise millions of units. The resulting LP relaxation for the entire network may then require millions of variables and constraints. Additionally, unless special precautions are taken, many of these constraints will be relatively dense. All this quickly leads to \change{LPs} that are \change{beyond the reach} of modern off-the-shelf LP solvers. As a result, researchers have explored alternative schemes for scaling LP-based methods to these larger networks.
\cite{salman2019convex} present a framework for LP-based methods (LP solvers, propagation, dual methods), which we review in the following subsections. However, they do not account for the ideal formulation developed in later works \citep{anderson2020strong,depalma2021scaling}.

\textbf{(Projecting the ideal).}
Figure~\ref{fig:relu-neuron} shows that the triangle (big-$M$) relaxation fails to recover the convex hull of the ReLU activation function and the multivariate input to the affine pre-activation function.  
We can similarly project the LP relaxation of the ideal formulation \eqref{eqn:relu-ideal} into the space of input/output variables \citep{anderson2020strong}, yielding a description for the convex hull of $\{ (\vh^{l-1},h^l_i) | L^{l-1} \leq \vh^{l-1} \leq U^{l-1}, \: h^l_i = \sigma(\vw^{l}_i \vh^{l-1} + b^l_i) \}$:
\begin{subequations}
\begin{align}
    h^l_i &\geq \vw^l_i \vh^{l-1} + b^l_i \\
    h^l_i &\leq \sum_{k \in I} w_{i,k}^l (x_k - \breve{L}_k) + \frac{\ell(I)}{\breve{U}_h - \breve{L}_h}(x_h - \breve{L}_h) \quad \forall (I,h) \in \mathcal{J} \label{eqn:projected-ideal-relu-2} \\
    (\vh^{l-1},h^l_i) &\in [L^{l-1},U^{l-1}] \times \mathbb{R}_{\geq 0},
\end{align}
\end{subequations}
where $l(I) \coloneqq \sum_{k \in I} w^l_{i,k}\breve{L}_k + \sum_{k \not\in I} w^l_{i,k} \breve{U}_k + b^l_i$ and
\begin{multline*}
    \mathcal{J} \coloneqq \{ (I,h) \in 2^{\llbracket n_{l-1} \rrbracket} \times \llbracket n_{l-1} \rrbracket | l(I) \geq 0, \: l(I \cup \{h\} < 0, \: w^l_{i,k} \neq 0 \forall k \in I \}.
\end{multline*}
\cite{anderson2020strong} also show that the inequalities \eqref{eqn:projected-ideal-relu-2} can be separated over in $\mathcal{O}(n_{l-1})$ time. 
Interestingly, in contrast to \eqref{eqn:relu-ideal}, the number of facet-defining inequalities depends heavily on the affine function. While in the worst case the number of inequalities will grow exponentially in the input dimension, there exist instances where the convex hull can be fully described with only $\mathcal{O}(n_{l-1})$ inequalities.

\subsubsection{Dual decomposition methods}

A first approach for greater scalability for LP-based methods is decomposition, a standard technique in the large-scale optimization community. Indeed, the cutting plane approach of Section~\ref{sec:strengthening-inequalities} can be viewed as a decomposition method operating in the original variable space. However, the method is initialized with the big-$M$ formulation for each neuron, and hence this initial model will be of size roughly equal to the size of the network. Therefore, it should be understood to use decomposition for providing a tighter verification bound, rather than for providing greater scalability to larger networks.

In contrast, dual decomposition can be used to scale inexact verification methods to larger networks. Such methods maintain dual feasible solutions throughout the algorithm, meaning that upon termination they yield valid dual bounds on the verification instance, and hence serve as incomplete verifiers.

\cite{wong2018provable,wong2018scaling} use as their starting point the triangle relaxation \eqref{eqn:triangle-relaxation} for each neuron, and then take the standard LP dual of the (relaxed) verification problem. Alternatively, \cite{dvijotham2018dual} propose a Lagrangian-based approach for decomposing the original nonlinear formulation of the problem \eqref{eq:verification}. Crucially, since the complicating constraints coupling the layers in the network are imposed as objective penalties instead of ``hard'' constraints, the optimization problem (given fixed dual variables) decomposes along each layer and the subproblems induced by the separability can be solved in closed form. This approach dualizes separately the equations characterizing the pre-activation and post-activation functions:
\begin{align*}
    \max_{\mu,\lambda}\quad \min_{\vh,\hat{\vh}} \quad& \left( \mW^L {\vh}^{L-1} + \vb^L \right) + \sum_{k=1}^{L-1} \Big( \mu_k^T(\hat{\vh}^k - \mW^k \vh^{k-1} - \vb^k)  + \lambda_k^T(\vh^k - \sigma(\hat{\vh}^k) \Big) \\
    \text{s.t.} \quad& L^k \leq \hat{\vh}^k \leq U^k \quad \forall k \in \llbracket n - 1 \rrbracket \\
    & \sigma(L^k) \leq \vh^k \leq \sigma(U^k) \quad \forall k \in \llbracket n - 1 \rrbracket.
\end{align*}

Here, the $\hat{\vh}$ variables track the pre-activation values for the neurons in the network. 
The dual variables $\mu_k^T$ correspond to the equality constraints defining the pre-activation values, $\hat{\vh}^k = \mW^k \vh^{k-1} + \vb^k$. 
Likewise, the dual variables $\lambda_k^T$ correspond to enforcing the ReLU activation function, $\vh^k = \sigma(\hat{\vh}^k) = \mathrm{max} (0, \hat{\vh}^k)$. 
Any feasible solution for the neural network is feasible for this dualized problem, making the multiplier terms for $\mu_k^T$ and $\lambda_k^T$ zero. 
Thus, the inner problem gives a lower bound for the original problem---a property known as \emph{weak duality}. 
The outer (dual) problem optimizing over the Lagrangian multipliers then seeks to maximize this lower bound, i.e., to give the tightest possible lower bound. This can be solved using a subgradient-based method, or learned along with the model parameters in a ``predictor-verifier'' approach \citep{dvijotham2018training}. 

On the other hand, this approach can be combined with other relaxation-based methods. 
The Lagrangian decomposition can be applied to dualize only the coupling constraints between layers, and a convex relaxation used for the activation function \citep{bunel2020lagrangian}: 
\begin{align*}
    \max_\lambda \quad \min_{\vh,\hat{\vh}} \quad& \left( \mW^L {\vh}^{L-1} + \vb^L \right) + \sum_{k=1}^{L-1} \left( \lambda_k^T(\vh^k - \sigma(\hat{\vh}^k) \right) \\
    \text{s.t.} \quad
    & L^k \leq \hat{\vh}^k \leq U^k \quad \forall k \in \llbracket n - 1 \rrbracket \\
    & \hat{\vh}^{k} = \mW^{k} \vh^{k-1} + b^k \quad \forall k \in \llbracket n - 1 \rrbracket \\
    & \vh^k \geq 0 \quad \forall k \in \llbracket n - 1 \rrbracket \\
    & \vh_i^k \geq \hat{\vh}_i^k \quad \forall k \in \llbracket n - 1 \rrbracket, \forall i \in \llbracket n_k \rrbracket \\
    & \vh^k_i \leq \frac{U^k_i(\hat{\vh}^k_i - L^k_i)}{U^k_i - L^k_i} \quad \forall k \in \llbracket n - 1 \rrbracket, \forall i \in \llbracket n_k \rrbracket.
\end{align*}
Note that the final three constraints apply the big-$M$/triangle relaxation \eqref{eqn:triangle-relaxation} to each ReLU activation function. 
The dual problem can then be solved via subgradient-based methods, proximal algorithms, or, more recently, a projected gradient descent method applied to a nonconvex reformulation of the problem \citep{bunel2020efficient}.

More recently, \cite{depalma2021scaling} presented a dual decomposition approach based on \eqref{eqn:relu-ideal}. 
However, creating a dual formulation from the exponential number of constraints produces an exponential number of dual variables. The authors therefore propose to maintain an ``active set'' of dual variables to keep the problem sparse. A selection algorithm (e.g., selecting entries that maximize an estimated super-gradient) can then be used to append the active set. 
Similar to the above discussion on cut generation, the frequency of appending the active set should be chosen strategically.

\subsubsection{Fourier-Motzkin elimination and propagation algorithms}
\newcommand{\defeq}{\vcentcolon=}
Alternatively, one can project out \emph{all} of the decision variables. For example, in order to solve the linear programming problem $\min_{x \in \mathcal{X}} c \cdot x$, we can augment the problem with a new decision variable to $\min_{(x,y) \in \Gamma} y$ for $\Gamma \defeq \Set{(x,y) \in \mathcal{X} \times \mathbb{R} : y = c \cdot x}$, and project out the $x$ variables. The transformed problem is the a trivial univariate optimization problem: $\min_{y \in \operatorname{Proj}_y(\Gamma)} y$.

Of course, the complexity of the approach described is hidden in the projection step, or building $\operatorname{Proj}_y(\Gamma)$. The most well-known algorithm for computing projections of linear inequality systems is Fourier-Motzkin elimination, described by \cite{dantzig1973fourier}, which is notorious for its practical inefficiency. 
The process effectively comprises replacing variables from a set of inequalities with all possible implied inequalities, which can produce many unnecessary constraints.  
However, it turns out that neural network verification problems are well-structured in such a way that Fourier-Motzkin elimination can be performed very efficiently: for instance, by imposing one inequality upper bounding and one inequality lower bounding each ReLU function. 
Note that while Section \ref{sec:algebraoflinearegions} describes the use of Fourier-Motzkin elimination to obtain \emph{exact} input-output relationships in linear regions of neural networks, here we are interested in obtaining linear \emph{bounds} for a nonlinear function.

\begin{figure}
    \centering
    \includegraphics[width=\textwidth]{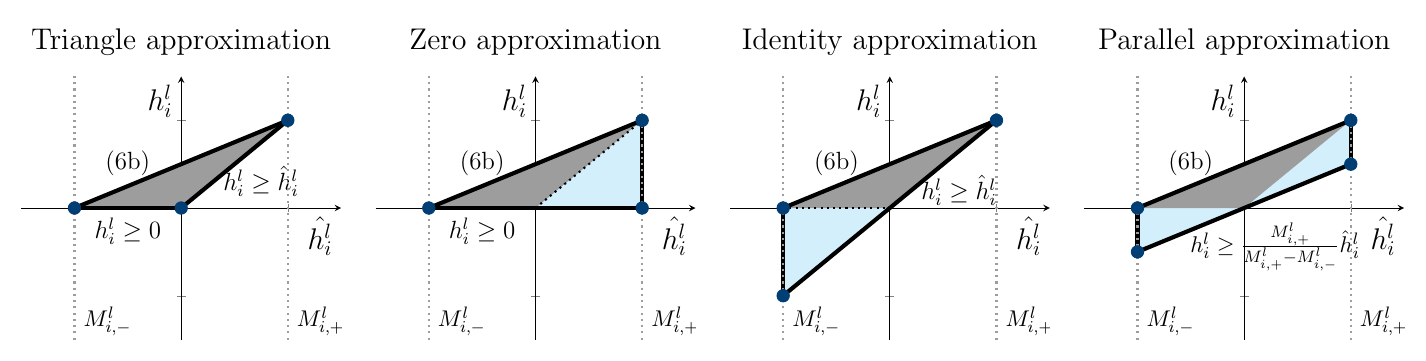}
    \caption{Convex approximations for the ReLU function commonly used by propagation algorithms, given as a function of the preactivation function $\hat{h_i^l}$. The ReLU applies $h_i^l = \max(0,\hat{h_i^l} )$.}
    \label{fig:convexapproximations}
\end{figure}

In fact, this general approach was independently developed in the verification community.
While MILP research has focused on formulations tighter than the big-M, such as \eqref{eqn:balas-relu} and \eqref{eq:relu-partition}, the verification community often prefers greater scalability at the price of weaker convex relaxations.
The continuous relaxation of the big-M is equivalent to the triangle relaxation \eqref{eqn:triangle-relaxation}: the optimal convex relaxation for a single input, or in terms of the aggregated pre-activation function, as shown in Figure \ref{fig:convexapproximations}.  
However, the lower bound involves two linear constraints, which is not used in several propagation-based verification tools owing to scalability or compatibility. 

Such tools use methods such as abstract transformers to propagate polyhedral bounds, i.e., \textit{zonotopes}, through the layers of a  neural network. 
DeepZ \citep{singh2018fast}, Fast-Lin \citep{weng2018towards}, and Neurify \citep{wang2018efficient} employ a parallel approximation, with the latter also implementing a branch-and-bound procedure towards completeness. 
Subsequently, DeepPoly \citep{singh2019abstract}  and CROWN \citep{zhang2018efficient} select between the zero and identity approximations by minimizing over-approximation area. 
OSIP \citep{hashemi2021osip} selects between the three approximations using optimization: approximations for a layer are select jointly to minimise bounds for the following layer. 
These technologies are also compatible with interval bounds, propagating box domains \citep{mirman2018differentiable}. 
Interestingly, bounds on neural network weights can also be propagated using similar methods, allowing reachability analysis of Bayesian neural networks \citep{wicker2020probabilistic} \newchange{and reachability analysis of the weights themselves during training steps~\citep{sosnin2024certified,wickercertification}}. 

\cite{tjandraatmadja2020convex} provide an interpretation of these propagation techniques through the lens of Fourier-Motzkin elimination. 
Consider the problem of propagating bounds through a ReLU neural network: for a node $h_i^l = \mathrm{max} \{0, \hat{h_i^l} \}$, convex bounds for $h_i^l$ can be obtained given bounds for $\hat{h}_i^l$ (Figure \ref{fig:convexapproximations}). Assuming the inputs are outputs of ReLU activations in the previous layer, $\hat{h}_i^l = \vw_i^l \vh^{l-1} + b_i^l$. 
Computing an upper bound can then be expressed as:
\begin{equation*}
    \begin{aligned}
    \max_{\vh^{l-1}} \quad& \vw_i^l \vh^{l-1} + b_i^l \\
    \text{s.t.} \quad& \mathcal{L}_k(\vh^{l-2}) \leq h^{l-1}_k \leq \mathcal{G}_k(\vh^{l-2}) \forall k \in \{ 1,...,n_{l-1} \}.
    \end{aligned}
\end{equation*}

As the objective function is linear, the solution of this problem can be computed by propagation without explicit optimization. 
For each element in $\vh^{l-1}$, we only need to consider the associated objective coefficient in $\vw_i^l$ to determine whether $\mathcal{L}_k(\vh^{l-2}) \leq h^{l-1}_k$ or $h^{l-1}_k \leq \mathcal{G}_k(\vh^{l-2})$ will be the active inequality at the optimal solution. 
We can thus replace $h^{l-1}_k$ with $\mathcal{L}_k(\vh^{l-2})$ or $\mathcal{G}_k(\vh^{l-2})$ accordingly. 
This projection is mathematically equivalent to applying Fourier-Motzkin elimination, while avoiding redundant inequalities resulting from the `non-selected' bounding function. 
Repeating this procedure for each layer results in a convex relaxation for the outputs that only involves the input variables.
We naturally observe the desirability of simple lower bounds $\mathcal{L}_k(\vh^{l-2})$: imposing two-part lower bounds in each layer would increase the number of propagated constraints in an exponential manner, similar to Fourier-Motzkin elimination.

\paragraph{A path towards completeness}
Given an input-output bound, the reachable set can be refined by splitting the input space \citep{henriksen2021deepsplit,rubies2019fast}---a strategy similar to spatial branch and bound. In other words, completeness is achieved by branching in the input space, rather than activation patterns: this strategy is especially effective when the input space is low-dimensional \citep{strong2022zope}. For example, ReluVal \citep{wang2018formal} propagates symbolic intervals and implements splitting procedures on the input domain. As the interval extension of ReLU is Lipschitz continuous, the method converges to arbitrary accuracy in a finite number of splits.

\subsection{Generalizing the single neuron model}

\subsubsection{Extending to other domains} \label{sec:domains}
In general, we will expect that the effective input domain $\mathcal{D}^{l-1}$ for a given unit may be quite complex. For the first layer ($l=1$) this may derive from explicitly stated constraints on the inputs of the networks, while for later layers this will typically derive from the complex nonlinear transformations applied by the preceding layers. 
For example, in the context of surrogate models \cite{yang2022modeling} propose bounding the input to the convex hull of the training data set, while other works \citep{schweidtmann2022obey,shi2022careful} propose machine learning-inspired techniques for learning the trust region implied by the training data. 
In effect, these methods assume a trained model is locally accurate around training data, which is a property similar to that which verification seeks to prove. 

Nevertheless, most research focuses on hyperrectangular input domains, largely motivated by practical considerations: i) there are efficient, well-studied methods for computing valid (though not necessarily optimally tight) variable bounds, ii) characterizing the exact effective domain may be computationally impractical, and iii) and the hyperrectangular structure makes analysis simpler for complex formulations like those presented in Section~\ref{sec:strengthening-inequalities}. 
We note that \cite{jordan2019neurips} use polyhedral analyses to perform verification over arbitrary (including non-polyhedral) norms, by fitting a $p$-norm ball in the decision region and checking adjacent linear regions. 
On the other hand, robust optimization can be employed to find $p$-norm adversarial regions (rather than verifying robustness), as opposed to a single point adversary \citep{maragno2023finding}. 

\cite{anderson2020strong} present two closely related frameworks for constructing ideal and hereditarily sharp formulations for ReLU units with arbitrary polyhedral input domains. This characterization is derived from Lagrangian duality, and requires an infinite number of constraints (intuitively, one for each choice of dual multipliers). Nonetheless, separation can still be done over this infinite family of inequalities via a subgradient-based algorithm; this approach will be tractable if optimization over $\mathcal{D}^{l-1}$ is tractable. 
Many propagation algorithms are also fully compatible with arbitrary polyhedral input domains, as the projected problem (i.e., a linear input-output relaxation) remains an LP. 
\cite{singh2021overcoming} show that simplex input domains can actually be beneficial, creating tighter formulations by propagating constraints on the inputs through the network layers.
Similarly, optimization-based bound tightening problems based on solving LPs can embed constraints defining polyhedral input domains.

In certain cases, additional structural information about the input domain can be used to reduce this semi-infinite description to a finite one. For example, this can be done when $\mathcal{D}^{l-1}$ is a Cartesian product of unit simplices \citep{anderson2020strong} (note that this generalizes the box domain case, wherein each simplex is one-dimensional). This particular structure is particularly useful for modeling input domains with combinatorial constraints. For example, a network trained to predict binding propensity of a given length-$n$ DNA sequence is naturally modeled via an input domain that is the product of $n$ 4-dimensional simplices--one simplex for each letter in the sequence, each of which is selected from an alphabet of length 4.

\subsubsection{Extending to other activation functions}
The big-$M$ formulation technique can be \change{directly applied to} any piecewise linear activation function. While much of the literature focuses on the ReLU due to its widespread popularity, models for other activation functions have been explored in the literature. For example, multiple papers \citep[Appendix K]{serra2018bounding} \citep[Appendix A.2]{tjeng2017evaluating} present a big-$M$ formulation for the maxout activation function. 
Adapting a formulation from \cite{anderson2020strong} \citep[Proposition 10]{anderson2020strong}, a formulation for the maxout unit is
\begin{align*}{2}
    y^l_i &\leq u_j(\vh^{l-1}) + M^l_{i,j}(1-z_j) \quad \forall j \in \llbracket k \rrbracket \\
    y^l_i &\geq u_j(\vh^{l-1}) \quad \forall j \in \llbracket k \rrbracket \\
    \sum_{j=1}^k z_j &= 1 \\
    (\vh^{l-1},v^l_i,z) &\in \mathcal{D}^{l-1} \times \mathbb{R} \times \{0,1\}^k,
\end{align*}
where each $M^l_{i,j}$ is selected such that
\[
    M^l_{i,j} \geq \max_{\tilde{\vh} \in \mathcal{D}^{l-1}} u_j(\tilde{\vh}).
\]

We can observe that the big-$M$ formulation can also handle other discontinuous activation functions, such as binary/sign activations \citep{han2021single} or more general quantized activations \citep{nguyen2022}. 
Nevertheless, the binary activation function naturally lends itself towards Boolean satisfiability, and most work therefore focuses on alternative methods such as SAT \citep{cheng2018verification,jia2020efficient,narodytska2018verifying}. 

While this survey focuses on neural networks with piecewise linear activation functions, we note that recent research has also studied smooth activation functions with a similar aim. 
For example, optimization over smooth activation functions can be handled by piecewise linear approximation and conversion to MILP \citep{sildir2022mixed}. 
Researchers have also studied convex/concave bounds for nonlinear activation functions, which can then be embedded in spatial branch-and-bound procedures \citep{schweidtmann2019deterministic,wilhelm2022convex}. 
\newchange{Notably, \citet{carrasco2024tightening} take a similar approach to the extended formulations described in Section~\ref{sec:MIPmodels} to derive tight formulations for convex and S-shaped activation functions; the relaxations were later implemented for the hyperbolic tangent activation function in a branch-and-bound solver by \citet{wittedeterministic}.}
In contrast to MILP formulations for ReLU neural networks, these problems are typically nonlinear programs that must be solved via spatial branch and bound.

Propagation methods \citep{singh2018fast,zhang2018efficient} can also naturally handle general activation functions: given convex polytopic bounds for an activation function, these tools can propagate them through network layers using the same techniques. 
For example, Fastened CROWN \citep{lyu2020fastened} employs a set of search heuristics to quickly select linear upper and lower bounds on ReLU, sigmoid, and hyperbolic tangent activation functions. 
\change{More recent work combines linear bounds with tailored branching rules for networks with general nonlinear functions~\citep{shi2024neural}.} 
Tighter polyhedral bounds can be employed, such as piecewise linear upper and lower bounds \citep{benussi2022individual}. 
\subsubsection{Extending to adversarial training} \label{sec:adversarialtraining}
As described in Section \ref{sec:intro}, the \emph{training} of neural networks seeks to minimise a measure of distance between the output $y$ and the correct output $\hat{y}$.
For instance, if this distance is prescribed as a loss function $\ERMfunction(y,\hat{y})$, this corresponds to solving the \emph{training} optimization problem: 
\begin{equation} \label{eq:training}
\underset{\{ \mW^l \}_{l \in \sL}, \{ \vb^l \}_{l \in \sL}}{\mathrm{min}} \ERMfunction(y,\hat{y}).
\end{equation}

Further details about the training problem and solution methods are described in Section~\ref{sec:training}. 
\change{An interesting direction here is to produce networks with bounded global~\citep{cisse2017parseval,leino2021globally} or local~\citep{huang2021training} Lipschitz constants---a paradigm with some connections to polyhedral theory. }
Nevertheless, given the focus of this section on optimization, we aim to briefly outline how verification techniques can be embedded in training. 
Specifically, solutions or bounds to the verification problem (Section \ref{sec:verification}) provide a metric of how robust a trained neural network is to perturbations. 
\newchange{These bounds can be used directly to bound the reachable set of network weights during training, given some perturbations~\citep{sosnin2024certified,wickercertification}.} 
Alternatively, verification-based metrics can be embedded in the training problem to obtain a more robust network during training, often resulting in a bilevel training problem. For instance, the verification problem \eqref{eq:verification} can be embedded as a lower-level problem, giving the robust optimization problem:
\begin{align*}
\underset{\{ \mW^l \}_{l \in \sL}, \{ \vb^l \}_{l \in \sL}}{\mathrm{min}} \quad \underset{||x - \hat{x}|| \leq \epsilon}{\mathrm{max}} \quad \ERMfunction(y=f(x),\hat{y}).
\end{align*}
Solving these problems generally involves either bilevel optimization, or computing an adversarial solution/bound at each training step, conceptually similar to the robust cutting plane approach. 
\cite{madry2018towards} proposed this formulation and solved the nonconvex inner problem using gradient descent, thereby losing a formal certification of robustness. 
These approaches may also benefit from reformulation strategies, such as by taking the dual of the inner problem and using any feasible solution as a bound \citep{wong2018provable} \newchange{or by using projection-based approaches \citep{francobaldi2025smle}}. 
The resulting models are not only more robust, but several works have also found it to be empirically easier to verify robustness in them  \citep{mirman2018differentiable,wong2018provable}. 

Alternatively, robustness can be induced by designing an additional penalty term for the training loss function, in a similar vein to regularization. For example: 
\begin{equation*} 
\underset{\{ \mW^l \}_{l \in \sL}, \{ \vb^l \}_{l \in \sL}}{\mathrm{min}} \kappa \ERMfunction(y,\hat{y}) + (1-\kappa) \ERMfunction_\mathrm{robust}(\cdot).
\end{equation*}

Additionally, if these robustness penalties are differentiable, they can be embedded into standard gradient descent based optimization algorithms \citep{dvijotham2018dual,mirman2018differentiable}. 
In the above formulation, the parameter $\kappa$ controls the relative weighting between fitting the training data and satisfying some robustness criterion, and its value can be scheduled during training, e.g., to first focus on model accuracy \citep{gowal2018effectiveness}.
In these cases, over-approximation of the reachable set is less problematic, as it merely produces a model \emph{more} robust than required. 
Nevertheless, \cite{balunovic2020adversarial} improve relaxation tightness by searching for adversarial examples in the ``latent'' space between hidden layers, reducing the number of propagation steps. \cite{zhang2020towards} provide an implementation that that tightens relaxations by also propagating bounds backwards through the network.

\section{Linear Programming and Polyhedral Theory in Training}\label{sec:training}

In the previous sections, we have almost exclusively focused on tasks involving neural networks that have already been constructed, i.e., we have assumed that the training step has already concluded (with the exception of Section \ref{sec:adversarialtraining}). In this section, we focus on the training phase, whose goal is to construct a neural network that can represent the relationship between the input and output of a given set of data points.

Let us consider a set of points, or sample, $(\sampleinput_i,\sampleoutput_i)_{i=1}^\samplesize$, and assume that these points are related via a function $\functolearn$, i.e., $\functolearn(\sampleinput_i) = \sampleoutput_i$ $i=1,\ldots, \samplesize$. In the training phase, we look for a function $\change{\tilde{f}}$ in a pre-defined class (e.g. neural networks with a specific architecture) \change{such that $\change{\tilde{f}}(\sampleinput_i) \approx \functolearn(\sampleinput_i)$}. Typically, this is done by solving an Empirical Risk Minimization problem
\begin{equation}
    \label{eq:ERMproblem}
    \min_{\change{\tilde{f}} \in F} \frac{1}{\samplesize} \sum_{i=1}^\samplesize \loss(\change{\tilde{f}}(\sampleinput_i), \sampleoutput_i)
\end{equation}
where $\loss$ is a loss function and $F$ is the class of functions we are restricted to. We usually assume the class $F$ is parametrized by $(\weights,\biases)\in \parameterset$ (the network weights and biases), so we are further assuming that there exists a function $\NNfunction{\cdot}{\cdot}{\cdot}$ (the network architecture) such that
\[\forall \change{\tilde{f}} \in F,\, \exists (\weights,\biases)\in \parameterset,\, \change{\tilde{f}}(\vx) = \NNfunction{\vx}{\weights}{\biases},\]
and thus, the optimization is performed over the space of parameters. In many cases, $\parameterset = \mathbb{R}^\parameternumber$---the parameters are unrestricted real numbers---but we will see some cases when a different parameter space can be used.

As mentioned in the introduction, nowadays, most of the practically successful \emph{training} algorithms for neural networks, i.e., that solve or approximate \eqref{eq:ERMproblem}, are based on Stochastic Gradient Descent (SGD). From a fundamental perspective, optimization problem \eqref{eq:ERMproblem} is typically a \emph{non-convex, unconstrained} problem that needs to be solved efficiently and where finding a \emph{local minimum} is sufficient. Thus, it is not too surprising that linear programming appears to be an unsuitable tool in this phase, in general. Nonetheless, there are some notable and surprising exceptions to this, which we review here.

Linear programming played an interesting role in training neural networks before SGD became the predominant training method and provided an efficient approach for constructing neural networks with 1 hidden layer in the 90s. This approach has some common points in its polyhedral approach with the first known algorithm that can solve \eqref{eq:ERMproblem} to provable optimality for a 1-hidden-layer ReLU neural network, which was proposed in 2018. Recently, a stream of work has exploited similar polyhedral structures to obtain convex optimization reformulations of regularized training problems of ReLU networks. Linear programming tools have also been used within SGD-type methods in order to compute optimal \emph{step-sizes} in the optimization of \eqref{eq:ERMproblem} or to strictly enforce structure in $\parameterset$. From a different perspective, a \emph{data-independent} polytope was used to describe approximately all training problems that can arise from an uncertainty set. Additionally, a back-propagation-like algorithm for training neural networks, which solves mixed-integer linear problems in each layer, was proposed as an alternative to SGD. Furthermore, when the neural network weights are required to be discrete, the applicability of SGD is impaired, and mixed-integer linear models have been proposed to tackle the corresponding training problems.

In what follows, we review these roles of (mixed-integer) linear programming and polyhedral theory within training contexts. We refer the reader to the book by \cite{goodfellow2016deep}, and the surveys by \cite{curtis2017optimization}, \cite{bottou2018optimization}, and \cite{wright2018optimization} for in-depth descriptions and analyses of the most commonly used training methods for neural networks.

We remark that solving the training problem to global optimality for ReLU neural networks is computationally complex. Even in architectures with just one hidden node, the problem is NP-hard \citep{dey2020approximation,goel2021tight}. Also see \cite{blum1992training,boob2022complexity,chen2022learning,froese2022computational,froese2023training} for other hardness results. Furthermore, it has been recently shown that training ReLU networks is $\exists \mathbb{R}$-complete \citep{abrahamsen2021training,bertschinger2022training}, which implies that it is likely that the problem of optimally training ReLU neural networks is not even in NP. 
Therefore, it is not strange to see that some of the methods we review below, even when they are solving hard problems as sub-routines (like mixed-integer linear problems), either make some non-trivial assumptions or relax some requirements. For example, boundedness and/or integrality of the weights, architecture restrictions such as the output dimension, or not having optimality guarantees. 

It is worth mentioning that, in contrast, for LTUs, exact exponential-time training algorithms are known for much more general architectures than in the ReLU case \citep{khalife2022neural,ergen2023globally,khalife2023neural}. These are out of scope of this survey, though we will provide a high-level overview of some of them, as they share some similarities to approaches designed for ReLU networks.

\subsection{Training neural networks with a single hidden layer}

Following the well-known XOR limitation of the perceptron \citep{perceptrons}, a natural interest arose in the development of training algorithms that could handle at least one hidden layer. In this section, we review training algorithms that can successfully minimize the training error in a one-hidden-layer setting and rely on polyhedral approaches.

\subsubsection{Problem setting and solution scheme}

Suppose we have a sample of size $\samplesize$ $(\sampleinput_i,\sampleoutput_i)_{i=1}^\samplesize$ where $\sampleinput_i \in \R^n$ and $\sampleoutput_i\in \mathbb{R}$. In a training phase, we would like to find a neural network function $\NNfunction{\cdot}{\cdot}{\cdot}$ that represents in the best possible way the relation $\NNfunction{\sampleinput_i}{\weights}{\biases} = \sampleoutput_i$.

Note that when a neural network $f$ has only one hidden layer, its behavior is almost completely determined by the sign of each component of the vector
\(\weights^1 x - \biases^1.\)
These are the cases of ReLU activations $\sigma(z) = \max\{0,z\}$ and LTU activations $\sigma(z) = \text{sgn}(z)$. The training algorithms we show here heavily exploit this observation and construct $(\weights^1,\biases^1)$ by embedding in this phase a \emph{hyperplane partition} problem based on the sample $(\sampleinput_i,\sampleoutput_i)_{i=1}^\samplesize$.
While the focus of this survey is mainly devoted to ReLU activations, we also discuss some selected cases with LTU activations as they share some similar ideas.

\subsubsection{LTU activations and variable number of nodes}

One stream of work dedicated to developing training algorithms for one-hidden-layer networks concerned the use of \emph{backpropagation} \citep{BackPOP1,BackPOP2,BackpropNN}. In the early 90s, an alternative family of methods was proposed, which was heavily based on linear programs (see, e.g., \cite{bennett1990neural,bennett1992robust,roy1993polynomial,mukhopadhyay1993polynomial}). These approaches can construct a 1-hidden-layer network without the need for an \emph{a-priori} number of nodes in the network. We illustrate the high-level idea of these next, based on the survey by \cite{mangasarian1993mathematical}.\\

Suppose that $\sampleoutput_i\in \{-1,1\}$, thus the NN we construct will be a classifier. The training phase can be tackled via the construction of a \emph{polyhedral partition} of $\mathbb{R}^n$ such that no two points (or few) $\sampleinput_i$ and $\sampleinput_j$ such that $\sampleoutput_i \neq \sampleoutput_j$ lie in the same element of the partition. To achieve this, the following approach presented by \cite{bennett1992robust} can be followed. Let $\setofones = \{i \in [\samplesize]\, :\, \sampleoutput_i = 1\}$ and $N = [\samplesize]\setminus \setofones$, and consider the following optimization problem
\begin{subequations}\label{eq:separationLP}
\begin{align}
    \min_{\vw,\change{\theta},y,z} \quad & \frac{1}{|\setofones|} \sum_{i\in \setofones} y_i + \frac{1}{|N|} \sum_{i\in N} z_i \\
    & \vw^\top \sampleinput_i - \change{\theta} + y \geq 1 && \forall i\in \setofones \\
    & - \vw^\top \sampleinput_i + \change{\theta} + z \geq 1 && \forall i \in N \\
    & y,z\geq 0.
\end{align}
\end{subequations}
This LP aims at finding a hyperplane $\vw^\top x  = \change{\theta}$ separating the data according to their value of $\sampleoutput_i$. Since the data may be \change{inseparable}, the LP is minimizing the following classification error
\[\frac{1}{|\setofones|} \sum_{i\in \setofones} (-\vw^\top \sampleinput_i + \change{\theta} + 1)_+ + \frac{1}{|N|} \sum_{i\in N} ( \vw^\top \sampleinput_i - \change{\theta} + 1 )_+
\]
The LP \eqref{eq:separationLP} is a linear reformulation of the latter minimization problem, where the auxiliary values $y,z$ take the value of each element in the sum.

Once the LP \eqref{eq:separationLP} is solved, we obtain \change{a halfspace} classifying our data points. In order to obtain a richer classification and lower error, we can iterate the procedure by means of the Multi-Surface Method Tree (MSMT, see \cite{bennett1992decision}), which solves a sequence of LPs as \eqref{eq:separationLP} in order to produce a polyhedral partition of $\mathbb{R}^n$. Let us illustrate how \change{one variant of} this procedure works in a simplified case. Assume that solving \eqref{eq:separationLP} results in a vector $\vw_1$ \change{and suppose that from it we can compute $a_1,b_1\in\mathbb{R}$ with $b_1 > a_1$ such that}
\[\{i \,: \, \vw_1^\top \sampleinput_i \geq b_1 \} \subseteq \setofones\quad \land \quad \{i \,: \, \vw_1^\top \sampleinput_i \leq a_1 \} \subseteq N, \]
%
\change{with at least one of the sets $\{(\sampleinput_i, \sampleoutput_i) \,: \, \vw_1^\top \sampleinput_i \geq b_1 \}$ or $\{(\sampleinput_i, \sampleoutput_i) \,: \, \vw_1^\top \sampleinput_i \leq a_1 \}$ being non-empty}\footnote{\change{For the case when this is not possible, see \cite{mangasarian1993mathematical} and references therein.}}. \change{In this case, we can remove these ``correctly classified'' points} from the data set. We then redefine \eqref{eq:separationLP} accordingly, in order to obtain a new vector $\vw_2$ and \change{then compute} scalars $b_2,a_2$ that would be used to classify within the region $\{x\in \R^n \,: \,  a_1 < \vw_1^\top \vx < b_1 \}$.

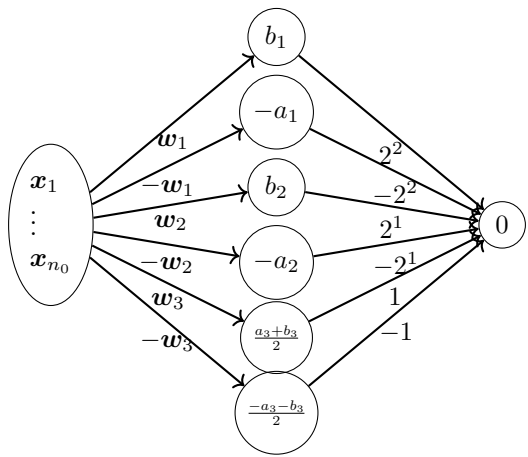
\begin{figure}
    \centering
\begin{tikzpicture}[auto]
    \node[draw, ellipse,align=left] (input) at (0,0) {$\vx_1$\\$\vdots$\\$\vx_\inputdimension$};
     \node[draw, circle] (h1) at (3,2.5) {$b_1$};
     \node[draw, circle] (h2) at (3,1.5) {$-a_1$};
    \node[draw, circle] (h3) at (3,0.5) {$b_2$};
     \node[draw, circle] (h4) at (3,-0.5) {$-a_2$};
     \node[draw, circle,scale=0.75] (h5) at (3,-1.6) {$\frac{a_3+b_3}{2}$};
     \node[draw, circle,scale=0.75] (h6) at (3,-2.8) {$\frac{-a_3-b_3}{2}$};
      \node[draw, circle] (output) at (6,0) {$0$};
     \draw[->,thick] (input) edge node[below]{$\vw_1$} (h1) 
     (input) edge node[below]{$-\vw_1$} (h2) 
     (input) edge node[below]{$\vw_2$} (h3) 
     (input) edge node[below]{$-\vw_2$} (h4)
     (input) edge node[below]{$\vw_3$} (h5) 
     (input) edge node[below]{$-\vw_3$} (h6)
     (h1) edge node[below]{$2^2$} (output)
     (h2) edge node[below]{$-2^2$} (output)
     (h3) edge node[below]{$2^1$} (output)
     (h4) edge node[below]{$-2^1$} (output)
     (h5) edge node[below]{$1$} (output)
     (h6) edge node[below]{$-1$} (output); 
\end{tikzpicture}
\caption{Illustration of Neural Network with LTU activations using MSMT. Inside each node of the hidden layer, we show the thresholds used in each LTU activation.}
    \label{fig:HyperplaneToNN}
\end{figure}

This procedure can be iterated, and the polyhedral partition of $\R^n$ induced by the resulting hyperplanes can be easily transformed into a Neural Network with 1 hidden layer and LTU activations (see \cite{bennett1990neural} for details). We illustrate this transformation with the following example: suppose that after 3 iterations, we have the following regions, with the arrow indicating to which class each region is associated:
\begin{subequations}\label{eq:hyperplanestoNN}
\begin{align}
&\{\vx\in \R^n \,: \, \vw_1^\top \vx \geq b_1 \} \to Y, \label{eq:firstYclass}\\
&\{\vx\in \R^n \,: \, \vw_1^\top \vx \leq a_1 \} \to N , \label{eq:firstNclass}\\
&\{\vx\in \R^n \,: \, a_1 < \vw_1^\top \vx < b_1,\,  \vw_2^\top \vx \geq b_2  \}\to Y, \label{eq:secondYclass}\\
&\{\vx\in \R^n \,: \, a_1 < \vw_1^\top \vx < b_1,\,  \vw_2^\top \vx \leq a_2  \}\to N, \label{eq:secondNclass}\\
&\{\vx\in \R^n \,: \, a_1 < \vw_1^\top \vx < b_1,\,  a_2 < \vw_2^\top \vx < b_2,\, \vw_3^\top x \geq (a_3 + b_3)/2 \} \to Y, \label{eq:thirdYclass}\\
&\{\vx\in \R^n \,: \, a_1 < \vw_1^\top \vx < b_1,\,  a_2 < \vw_2^\top \vx < b_2,\,  \vw_3^\top x < (a_3 + b_3)/2 \} \to N \label{eq:thirdNclass}.
\end{align}
\end{subequations}
Since regions \eqref{eq:thirdYclass} and \eqref{eq:thirdNclass} are the last defined by the algorithm (under some stopping criterion), they both use $(a_3 + b_3)/2$ in order to obtain a well-defined partition of $\mathbb{R}^n$. In Figure \ref{fig:HyperplaneToNN} we show a one-hidden-layer neural network that represents such a classifier. The structure of the neural network can be easily extended to handle more regions.
For other details, we refer the reader to \cite{bennett1992decision}, and for variants and extensions, see \cite{mangasarian1993mathematical} and references therein.

Some key features of this procedure are the following:
\begin{itemize}
    \item Each solution of \eqref{eq:separationLP}, i.e., each new hyperplane, can be represented as a new node in the hidden layer of the resulting neural network.
    \item The addition of a new hyperplane comes with a reduction in the current loss; this can be iterated until a target loss is met.
    \item Thanks to the universal approximation theorem \citep{Hornik1989}, with enough nodes in the hidden layer, we can always obtain a neural network $f$ with zero classification error; although, this can lead to over-fitting.\\
\end{itemize}

The work of \cite{roy1993polynomial} and \cite{mukhopadhyay1993polynomial} follow a related idea, although the classifiers \change{that} are built are quadratic functions. To illustrate the approach, we use the same set-up for \eqref{eq:separationLP}. The approach in \cite{roy1993polynomial} and \cite{mukhopadhyay1993polynomial} aims at finding a function
\[f_{\mV,\vw,b}(x) = \vx^\top \mV \vx + \vw^\top \vx + b \]
such that
\[f_{\mV,\vw,b}(\sampleinput_i) \geq 0 \Longleftrightarrow i\in \setofones.\]
Since this may not be possible, the authors propose solving the following LP
\begin{subequations}\label{eq:separationLP2}
\begin{align}
    \min_{\change{\mV,\vw},b,\epsilon} \quad & \epsilon \\
    & \sampleinput_i^\top \mV \sampleinput_i + \vw^\top\sampleinput_i + b  \geq \epsilon && \forall i\in Y \\
    & \sampleinput_i^\top \mV \sampleinput_i + \vw^\top\sampleinput_i + b  \leq -\epsilon && \forall i\not\in Y \\
    & \epsilon \geq \epsilon_0
\end{align}
\end{subequations}
for some fixed tolerance $\epsilon_0 >0$. When this LP is infeasible, the class $\setofones$ is partitioned into $\setofones_1$ and $\setofones_2$, and an LP as \eqref{eq:separationLP2} is solved for both $Y_1$ and $Y_2$. The algorithm then follows iteratively (see below for comments on these iterations). In the end, the algorithm will construct $k$ quadratic functions $f_1, \ldots, f_k$, which the authors call ``masking functions'', that will classify an input $\vx$ in the class $\setofones$ if and only if
\(\exists i\in [k],\, f_i(\vx) \geq 0.\)

In order to represent the resulting classifier as a single-layer neural network, the authors proceed in a similar manner to a linear classifier; the input layer of the resulting neural network not only includes each entry of $\vx$, but also the bilinear terms $\vx\vx^\top$. Using this input, the classifier built by \eqref{eq:separationLP2} can be thought of as a linear classifier (much like a polynomial regression can be cast as a linear regression).

As a last comment on the work of \cite{roy1993polynomial} and  \cite{mukhopadhyay1993polynomial}, the authors' algorithm does not iterate in a straightforward fashion. They add clustering iterations alternating with the steps described above in order to (a) identify outliers and remove them from the training set, and  (b) subdivide the training data when \eqref{eq:separationLP2} is infeasible. These additions allow them to obtain a polynomial-time algorithm.\\

The methods described in this section are able to produce a neural network with arbitrary quality; however, there is no guarantee on the size of the resulting neural network. When the size of the network is fixed, the story changes, which is the case we will describe next.

\subsubsection{Fixed number of nodes and ReLU activations}
\label{sec:arora}
As mentioned at the beginning of this section, training a neural network
is a complex optimization problem in general, with some results indicating that the problem is likely to not even be in NP \citep{abrahamsen2021training,bertschinger2022training}.

Nonetheless, by restricting the network architecture sufficiently and allowing exponential running times, exact algorithms can be conceived. An important step in the construction of such algorithms was taken by  \cite{arora2018understanding}. In this work, the authors studied the training problem in detail, providing the first 
optimization algorithm capable of solving the training problem to provable optimality for a fixed network architecture with one hidden layer and with an output dimension of 1. As we anticipated, this algorithm shares some similarities \change{with the previous approaches discussed above}.

Let us \change{now consider} a ReLU activation. Also, we no longer assume $\sampleoutput_i \in \{-1,1\}$, but we keep the output dimension as 1. The problem considered by \cite{arora2018understanding} reads
\begin{equation}
    \label{eq:ERMproblemArora}
    \min_{\weights,\biases} \frac{1}{\samplesize} \sum_{i=1}^\samplesize \loss (\weights^2(\sigma(\weights^1\sampleinput_i + \biases^1)), \sampleoutput_i),
\end{equation}
with $\loss: \R \times \R \to \R$ a convex loss. Note that this problem, even if $\loss$ is convex, is a non-convex optimization problem.

\begin{theorem}[\cite{arora2018understanding}]
Let $\layerwidth_1$ be the number of nodes in the hidden layer. There exists an algorithm to find a global optimum of \eqref{eq:ERMproblemArora} in time $O(2^{\layerwidth_1} \samplesize^{\inputdimension \cdot \layerwidth_1} \text{poly}(\samplesize,\inputdimension,\layerwidth_1))$.
\end{theorem}

Roughly speaking, the algorithm works by noting that one can assume the weights in $\weights^2$ are in $\{-1,1\}$, since $\sigma$ is positively-homogeneous. Thus, problem \eqref{eq:ERMproblemArora} can be restated as
\begin{equation}
    \label{eq:ERMproblemArora-simplified}
    \min_{\weights^1,\biases^1, s} \frac{1}{\samplesize} \sum_{i=1}^\samplesize \loss (s (\sigma(\weights^1\sampleinput_i + \biases^1)), \sampleoutput_i)
\end{equation}
where $s\in \{-1,1\}^{\layerwidth_1}$. In order to handle the non-linearity, \cite{arora2018understanding} ``guess'' the values of $s$ and the sign of each component of $\weights^1\sampleinput_i + \biases^1$. Enforcing a sign for each component of $\weights^1\sampleinput_i + \biases^1$ is similar to the approach discussed in the previous section: it fixes how the input part of the data $(\sampleinput_i)_{i=1}^\samplesize$ is partitioned in polyhedral regions by a number of hyperplanes. The difference is that, in this case, the number of hyperplanes to be used is assumed to be fixed.

Using the hyperplane arrangement theorem (see e.g. \cite[Proposition 6.1.1]{matousek2002lectures}), there are at most $\samplesize^{\inputdimension \layerwidth_1}$ ways of fixing the signs of $\weights^1 \sampleinput_i + \biases^1$. Additionally, there are at most $2^{\layerwidth_1}$ possible vectors in $\{-1,1\}^{\layerwidth_1}$. Once these components are fixed, \eqref{eq:ERMproblemArora-simplified} can be solved as an optimization problem with a convex objective function and a polyhedral feasible region imposing the desired signs in $\weights^1\sampleinput_i + \biases^1$. This results in the $O(2^{\layerwidth_1} \samplesize^{\inputdimension \layerwidth_1} \text{poly}(\samplesize,\inputdimension,\layerwidth_1))$ running time.
This algorithm was recently generalized to concave loss functions by \cite{froese2022computational}.

\cite{dey2020approximation} developed a polynomial-time approximation algorithm in this setting for the case of $\layerwidth_1 = 1$ (i.e., one ReLU neuron) and square loss. This approximation algorithm has a better performance when the input dimension is much larger than the sample size, i.e., $\inputdimension \gg \samplesize$. The approach by \cite{dey2020approximation} also relies on fixing the signs of $\weights^1\sampleinput_i + \biases^1$, and then solving multiple convex optimization problems, \change{but employs a different strategy than} \cite{arora2018understanding}; in particular, \cite{dey2020approximation} only explore a polynomial number of the possible ``fixings'', which yields the approximation. 

We note that the result by \cite{arora2018understanding} shows that the training problem on their architecture is in NP. This is in contrast to  \cite{bertschinger2022training}, who show that training a neural network with one hidden layer is likely to not be in NP. The big difference lies in the assumption on the output dimension: in the case of \cite{bertschinger2022training}, the output dimension is two. It is quite remarkable that such a sharp complexity gap is produced by a small change in the output dimension.

\subsubsection{An exact training algorithm for arbitrary LTU architectures}

Recently, \cite{khalife2023neural} presented a new algorithm akin to that in \cite{arora2018understanding}, capable of solving the training problem to global optimality for any fixed LTU architecture with a convex loss function $\loss$. In this case, no assumption on the network's depth is made. The algorithm runs in polynomial time on the sample size $\samplesize$ when the architecture is fixed.

We will not describe this approach in detail, as it heavily relies on the structure given by LTU activations, which is intricate and beyond the scope of this survey.
Although we note that it shares some high-level similarities to the algorithm of \cite{arora2018understanding} for ReLU activations, such as ``guessing" the behavior of the neurons' activity and then solving multiple convex optimization problems. However, the structural and algorithmic details are considerably different.

It is important to note that this result reveals the big gap between what is known for LTU versus ReLU activations in terms of their training problems. In the case of the former, there is an exact algorithm for arbitrary architectures; in the case of the latter, the known results are much more restricted and strong computational limitations exist.

\subsection{Convex reformulations in regularized training problems}

For the case when the training problem is regularized, the following stream of work has developed several convex reformulations of it.
\cite{pilanci2020neural} presented the first convex reformulation of a training problem for the case with one hidden layer and one-dimensional outputs. As the approach described in Section \ref{sec:arora}, this reformulation uses hyperplane arrangements according to the activation patterns of the ReLU units, but instead of using them algorithmically directly, they use them to find their convex reformulations. 
This framework was further extended to CNNs by \cite{ergen2021implicit}.
Higher-dimensional outputs in neural networks with one hidden layer were considered in \cite{pmlr-v108-ergen20a,ergen2021convex,sahiner2021vectoroutput}.
Recently, \cite{sahiner2024scaling} proposed improving tractability of these convex problems using the Burer-Monteiro Factorization \citep{burer2003nonlinear}.
This convex optimization perspective was also applied in Batch Normalization by \cite{ergen2022demystifying}. 

These approaches provide polynomial-time algorithms when some parameters (e.g., the input dimension $\inputdimension$) are considered constant. We note that this \change{seems to contradict the results of \cite{froese2023training}, which showed that training remains NP-Hard even
when the input dimension is fixed, there is only one hidden layer, and the output is one-dimensional. The subtlety relies on the regularization: the results discussed in this section include a regularizing term, whereas \cite{froese2023training} does not}. We explain below where the regularizing term plays an important role.
Training via convex optimization was further developed to handle deeper regularized neural networks in \cite{ergen2021global,ergen2021path,ergen2021revealing}.

In what follows, we review the convex reformulation in \cite{pilanci2020neural} (one hidden layer and one-dimensional output) to illustrate some of the base strategies behind these approaches. See \cite{ergen2023convex} for an extended version of this work.
We refer the reader to the previously mentioned articles for the most recent and intricate developments, as well as numerical experiments.\\

As before, let $\layerwidth_1$ be the number of nodes in the hidden layer. Let us consider the following regularized training problem; to simplify the discussion, we omit biases.
\begin{align}
    \min_{\weights} \quad & \frac{1}{2} \left\|\sum_{j=1}^{\layerwidth_1} \weights^2_j\sigma(\tilde{\mX}\weights^1_j) - \sampleoutput \right\|^2  + 
    \frac{\beta}{2}\sum_{j=1}^{\layerwidth_1}( \| \weights^1_j\|^2 + (\weights^2_j)^2)  \label{eq:ERMproblemCVX}
\end{align}
Here, $\beta > 0$, $\tilde{\mX}$ is a matrix whose $i$-th row is $\sampleinput_i$ and $\weights^1_j$ is the vector of weights going into neuron $j$. Thus, $\tilde{\mX}\weights^1_j$ is a vector whose $i$-th component is the input to neuron $j$ when evaluating the network on $\sampleinput_i$. $\weights^2_j$ is a scalar: it is the weight on the arc from neuron $j$ to the output neuron (one-dimensional). Note that there is a slight notation overload: $(\weights^2_j)^2$ is the square of the scalar $\weights^2_j$. However, we will quickly remove this (\change{potentially} confusing) term.

Problem \eqref{eq:ERMproblemCVX} is a regularized version of \eqref{eq:ERMproblemArora} when $\ell$ is the squared difference. We modified its presentation to match the structure in \cite{pilanci2020neural}. The authors first prove that \eqref{eq:ERMproblemCVX} is equivalent to 
\begin{align*}
\min_{\|\weights^1_j\| \leq 1} \min_{\weights^2_j} \quad & \frac{1}{2} \left\|\sum_{j=1}^{\layerwidth_1} \weights^2_j\sigma(\tilde{\mX}\weights^1_j) - \sampleoutput \right\|^2  + \beta\sum_{j=1}^{\layerwidth_1} | \weights^2_j|
\end{align*}
Then, through a series of reformulations and duality arguments, the authors first show that this problem is lower bounded by
\begin{subequations} \label{eq:ERMproblemCVX2}
\begin{align}
    \max & \quad - \frac{1}{2} \left\|v - \sampleoutput \right\|^2 + 
    \frac{1}{2}\|\sampleoutput\|^2\\
    \mbox{s.t} & \quad |v^\top \sigma(\tilde{\mX} u )| \leq \beta \qquad \qquad \forall u, \, \|u\|\leq 1\\
    &\quad v\in \mathbb{R}^\samplesize
\end{align}
\end{subequations}
Problem \eqref{eq:ERMproblemCVX2} has $\samplesize$ variables and infinitely many constraints.
The authors show that this lower bound is tight when the number of neurons in the hidden layer is large enough; specifically, they require $\layerwidth_1 \geq m^*$, where $m^* \in \{1,\ldots, \samplesize\}$ is defined as the number of Dirac deltas in an optimal solution of a dual of \eqref{eq:ERMproblemCVX2} (see \cite{pilanci2020neural,ergen2023convex} for details). 

Regarding the presence of infinitely many constraints, the authors address this by considering all possible patterns of signs of $\tilde{\mX} u$ (similarly to \cite{arora2018understanding}, as discussed in Section \ref{sec:arora}). For each fixed sign pattern (hyperplane arrangement), they apply a duality argument which allows them to recast the constraint $\max_{u\in \mathcal{B}}|v^\top \sigma(\tilde{\mX} u )|  \leq \beta$ as a finite collection of second-order cone constraints with $\beta$ on the right-hand side.

Finally, using that $\beta > 0$, they show that the reformulated problem satisfies Slater's condition, and thus from strong duality they obtain the following convex optimization problem, which has the same objective value as \eqref{eq:ERMproblemCVX2}.
\begin{subequations} \label{eq:finalcvx}
    \begin{align}
    \min&\quad \frac{1}{2} \left\|\sum_{j=1}^{P} M_i \tilde{\mX}(v_i - w_i) - \sampleoutput \right\|^2  + \beta \sum_{j=1}^{P}( \| v_i\| + \|w_i\|)\\
    \mbox{s.t} & \quad  (2M_i - I_\samplesize) \tilde{\mX} v_i \geq 0  && \forall i\in [P]\\
    & \quad  (2M_i - I_\samplesize) \tilde{\mX} w_i \geq 0 && \forall i\in [P] \\
    & \quad  v_i \in \mathbb{R}^\inputdimension  &&\forall i\in [P]\\
    & \quad w_i \in \mathbb{R}^\inputdimension && \forall i\in [P]
    \end{align}
\end{subequations}
Here, $I_\samplesize$ is the $\samplesize \times \samplesize$ identity matrix, $P$ is the number of possible activation patterns for $\tilde{\mX}$, and each  $M_i$ is a $\samplesize \times \samplesize$ binary diagonal matrix whose diagonal indicates the $i$-th possible sign pattern of $\tilde{\mX} u$. This means that $(M_{i})_{j,j}$ is 1 if and only if $\sampleinput_j^\top u \geq 0$ in the $i$-th sign pattern of $\tilde{\mX} u$. Moreover, the authors provide a formula to recover a solution of \eqref{eq:ERMproblemCVX} from a solution of \eqref{eq:finalcvx}.

Using that $P\leq 2r(e(\samplesize-1)/r)^r$, where $r=\mbox{rank}(\tilde{\mX})$, the authors note that the formulation \eqref{eq:finalcvx} yields a training algorithm with complexity $O(\inputdimension^3 r^3 (\samplesize/r)^{3r})$. Note that if one fixes $r$, the resulting algorithm runs in polynomial time. In particular, fixing $\inputdimension$ fixes the rank of $\tilde{\mX}$ and results in a polynomial time algorithm as well. In contrast, the algorithm by \cite{arora2018understanding} discussed in Section \ref{sec:arora} remains exponential even after fixing $\inputdimension$. Moreover,  \cite{froese2023training} showed that the training problem is NP-Hard even for fixed $\inputdimension$.
This apparent contradiction is explained by two key components of the convex reformulation: the regularization term and the presence of a ``large enough'' number of hidden neurons. This facilitates the exponential improvement of the training algorithm with respect to \cite{arora2018understanding}.

\subsection{Frank-Wolfe in DNN training algorithms}

Another stream of work that has included components of linear programming in DNN training involves the Frank-Wolfe method. We briefly describe this method in the non-stochastic version next. In this section, we omit the biases $\biases$ to simplify the notation.\\

Gradient descent (and its variants) is designed for problems of the form
\begin{equation}\label{eq:unconstrained}
\min_{\weights \in \mathbb{R}^\parameternumber} \ERMfunction(\weights)
\end{equation}
and it is based on iterations of the form
\begin{equation}\label{eq:gradientupdate}
\weights(i+1) = \weights(i) - \alpha_i \nabla \ERMfunction(\weights(i)) 
\end{equation}
where $\alpha_i$ is known as the \emph{learning rate}. In the stochastic versions, $\nabla \ERMfunction(\weights(i))$ is replaced by a stochastic gradient.
In this setting, these algorithms would find a local minimum, which is global when $\ERMfunction$ is convex. 

In the presence of constraints $\weights \in \parameterset$, however, this strategy may not work directly. A regularizing term is typically used in the objective function instead of a constraint, that ``encourages'' $\weights \in \parameterset$ but does not enforce it. If we strictly require that $\weights \in \parameterset \neq \mathbb{R}^n$, and $\parameterset$ is a convex set, one could modify \eqref{eq:gradientupdate} to
\begin{equation}\label{eq:projectionphi}
\weights(i+1) = \text{proj}_{\parameterset} \left( \weights(i) - \alpha_i \nabla \ERMfunction(\weights(i)) \right).    
\end{equation}
and thus ensure that all iterates $\weights(i) \in \parameterset$. Unfortunately, a projection is a costly routine. An alternative to this projection is the Frank-Wolfe method \citep{frank1956algorithm}. Here, a direction $\vd_i$ is computed via the following linear-objective convex optimization problem
\begin{equation}
    \label{eq:FrankWolfeLP}
    \vd_i \in \arg\min_{\vd \in \parameterset} \vv_i^\top \vd
\end{equation}
where normally $\vv_i = \nabla \ERMfunction(\weights(i))$ (we consider variants below). The update is then computed as
\begin{equation}\label{eq:FrankWolfeUpdate}
\weights(i+1) = \weights(i) + \alpha_i (\vd_i - \weights(i)),
\end{equation}
for $\alpha_i\in [0,1]$. Note that, by convexity, we are assured that $\weights(i+1) \in \parameterset$ as long as $\weights(0)\in \parameterset$. In many applications, $\parameterset$ is polyhedral, which makes \eqref{eq:FrankWolfeLP} a linear program. Moreover, for simple sets $\parameterset$, problem \eqref{eq:FrankWolfeLP} admits closed-form solutions.\\

In the context of deep neural network training, two notable applications of Frank-Wolfe have appeared. Firstly, the Deep Frank Wolfe algorithm, by  \cite{berrada2018deep}, modifies iteration \eqref{eq:gradientupdate} with an optimization problem that can be solved using Frank-Wolfe in its dual. Secondly, the use of a stochastic version of Frank-Wolfe in the training problem \eqref{eq:ERMproblem} by \cite{pokutta2020deep} and \cite{xie2020efficient}, which enforces structure in the neural network weights directly. We review these next, starting with the latter.

\subsubsection{Stochastic Frank-Wolfe}
Note that problem \eqref{eq:ERMproblem} is of the form \eqref{eq:unconstrained} with
\[\ERMfunction(\weights) = \frac{1}{\samplesize} \sum_{i=1}^\samplesize \loss(f(\sampleinput_i, \weights), \sampleoutput_i).  \]
We remind the reader that we are omitting the biases in this section to simplify notation, as they can be incorporated as part of $\weights$.

Usually, some structure of the weights is commonly desired (e.g., sparsity or boundedness), which traditionally have been incorporated as regularizing terms in the objective, as mentioned above.
The recent work by \cite{xie2020efficient} and \cite{pokutta2020deep}, on the other hand, enforce structure on $\parameterset$ directly using Frank-Wolfe ---more precisely, stochastic versions of it.\\

\cite{xie2020efficient} use a stochastic Frank-Wolfe approach to impose an $\ell_1$-norm constraint on the weights and biases $\weights$ when training a neural network with 1 hidden layer. Note that $\ell_1$ constraints are polyhedral.
 Their algorithm is designed for a general Online Convex Optimization setting, where ``losses'' are revealed in each iteration. However, in their computational experiments, they included tests in an offline setting given by a DNN training problem. 
 
 The approach closely follows the Frank-Wolfe method described above. The key difference lies in the estimation of the stochastic gradient they use, which is not standard, and it is one of the most important aspects of the algorithm. Instead of using $\vv_i = \nabla \ERMfunction(\weights(i))$ in \eqref{eq:FrankWolfeLP}, the following \emph{stochastic recursive estimator} of the gradient is used:
 \begin{align*}
 \vv_0 =& \tilde{\nabla}\ERMfunction(\weights(0)) \\
 \vv_i =& \tilde{\nabla}\ERMfunction(\weights(i)) + (1-\rho_i)(v_{i-1} - \tilde{\nabla}\ERMfunction(\weights(i-1)))
 \end{align*}
 where $\tilde{\nabla}\ERMfunction$ is a stochastic gradient, and $\rho_i$ is a parameter. The authors show that the gradient approximation error of this estimator converges to 0 at a sublinear rate, with high probability. This is important for them to analyze the ``regret bounds'' they provide for the online setting.

The experimental results in \cite{xie2020efficient} in DNN training are very positive. They test their approach in the MNIST and CIFAR10 datasets and outperform existing state-of-the-art approaches in terms of suboptimality, training accuracy, and test accuracy.\\

\cite{pokutta2020deep} implement and test several variants of stochastic versions of Frank-Wolfe in the training of neural networks, including the approach by \cite{xie2020efficient}. \cite{pokutta2020deep} focus their experiments on their main proposed variant, which they refer to simply as Stochastic Frank-Wolfe (SFW). This variant uses
\[\vv_i = (1-\rho_i)\vv_{i-1} + \rho_i \tilde{\nabla}\ERMfunction(\weights(i)),\]
where $\rho_i$ is a momentum parameter. The authors propose many different options for $\parameterset$ including $\ell_1, \ell_2$ and $\ell_\infty$ balls, and $K$-sparse polytopes. Of these, only the $\ell_2$ ball is non-polyhedral.

Overall, the computational experiments are promising for SFW. The authors advocate for this algorithm, arguing that it provides excellent computational performance while being simple to implement and competitive with other state-of-the-art algorithms.

\subsubsection{Deep Frank-Wolfe}

Another application of Frank-Wolfe within DNN training was proposed by \cite{berrada2018deep}. While this approach does not make heavy use of linear programming techniques, the application of Frank-Wolfe is quite novel, and they do rely on one linear program needed when performing an update as \eqref{eq:FrankWolfeUpdate}. \\

The authors note that \eqref{eq:gradientupdate} can also be written as the solution to the following \emph{proximal} problem \citep{bubeck2015convex}:
\begin{align}\label{eq:proximalupdate}
\begin{split}
    \weights(i+1) = \arg\min_{\weights}\,\biggl\{ & \frac{1}{2\alpha_i}\|\weights - \weights(i)\|^2  + \mathcal{T}_{\weights(i)}(\ERMfunction(\weights)) \biggr\}
    \end{split}
\end{align}
where $\mathcal{T}_{\weights(i)}$ represents the first-order Taylor expansion at $\weights(i)$. We are omitting regularizing terms since they do not play a fundamental role in the approach; all this discussion can be directly extended to include regularizers. \cite{berrada2018deep} note that \eqref{eq:proximalupdate} linearizes the loss function, and propose the following \emph{loss-preserving proximal} problem to replace \eqref{eq:proximalupdate}: 
\begin{align}\label{eq:lossproximalupdate}
\begin{split}
    \weights(i+1) = & \arg\min_{\weights}\,\biggl\{  \frac{1}{2\alpha_i}\|\weights - \weights(i)\|^2 + \frac{1}{\samplesize} \sum_{i=1}^\samplesize \loss(\mathcal{T}_{\weights(i)}(f(\sampleinput_i, \weights)), \sampleoutput_i) \biggr\}
    \end{split}
\end{align}

Using the results by \cite{lacoste2013block}, the authors argue that \eqref{eq:lossproximalupdate} is amenable to Frank-Wolfe in the dual when $\loss$ is piecewise linear and convex (e.g., the hinge loss). To be more specific, the authors show that in this case, and assuming $\alpha_i = \alpha$, there exists $\mA,\vb$ such that the dual of \eqref{eq:lossproximalupdate} is simply
\begin{subequations} \label{eq:dualDFW}
\begin{align}
 \max_{\mathbf{\beta}} \quad  &  \frac{-1}{2\alpha} \|\mA \mathbf{\beta} \|^2 + \vb^\top \mathbf{\beta} \\
 \text{s.t.} \quad & \mathbf{1}^\top \mathbf{\beta} = 1 \\
 & \mathbf{\beta} \geq 0
\end{align}
\end{subequations}
The authors consider applying Frank-Wolfe to this last problem, and to recover the primal solution using the primal-dual relation $\weights = -\mA\mathbf{\beta}$, which is a consequence of KKT. The Frank-Wolfe iteration \eqref{eq:FrankWolfeUpdate} in the notation of \eqref{eq:dualDFW} would look like
\begin{equation}\label{eq:FrankWolfeLP-dual}
\mathbf{\beta}_{i+1} = \mathbf{\beta}_{i} + \gamma_i (\vd_i - \mathbf{\beta}_i).
\end{equation}
Here, $\vd_i$ is feasible for \eqref{eq:dualDFW} and obtained using a linear programming oracle, and $\gamma_i$ is the Frank-Wolfe step-length. Note that the feasible region of \eqref{eq:dualDFW} is a simplex: exploiting this, the authors show that an optimal $\gamma_i$ can be computed in closed-form: here, ``optimal'' refers to a minimizer of \eqref{eq:dualDFW} when restricted to points of the form $\mathbf{\beta}_{i} + \gamma_i (\vd_i - \mathbf{\beta}_i)$.

With all these considerations, the bottleneck in this application of Frank-Wolfe is obtaining $\vd_i$; recall that this Frank-Wolfe routine is embedded within a single iteration of the overall training algorithm; therefore, in each iteration of the training algorithm, possibly multiple computations of $\vd_i$ would be required in order to solve \eqref{eq:dualDFW} to optimality. To alleviate this, the authors propose to perform only one iteration of Frank-Wolfe: they set $\vd_0$ to be the stochastic gradient and compute a closed-form expression for $\mathbf{\beta}_{1}$. This is the basic ingredient of the Deep Frank Wolfe (DFW). It is worth noting that this algorithm is not guaranteed to \change{converge; however,} its empirical performance is competitive.

Other two important considerations are taken into account the implementation of this algorithm: smoothing of the loss function (as the Hinge loss is piecewise linear) and the adaptation of Nesterov's Momentum to this new setting. We refer the reader to the corresponding article for these details. One of the key features of DFW is that it only requires one hyperparameter ($\alpha$) to be tuned.

The authors test DFW in image classification and natural language inference. Overall, the results obtained by DFW are very positive: in most cases, it can outperform adaptive gradient methods, and it is competitive with SGD while converging faster. 

\subsection{Polyhedral encoding of multiple training problems}

One of the questions raised by \cite{arora2018understanding} (see Section \ref{sec:arora}) was whether the dependency on $\samplesize$ of their algorithm could be improved since it is typically the largest coefficient in a training problem. This question was studied by \cite{TrainingLP}, who show that, in an approximation setting, a more ambitious goal is achievable: there is a polyhedral encoding of multiple training problems whose size has a mild dependency on $D$.

As in the previous section, we omit the biases $\biases$ to simplify notation, as all parameters can be included in $\weights$. Let us assume the class of neural networks $F$ in \eqref{eq:ERMproblem} are restricted to have bounded parameters (we assume they lie in the interval $[-1,1]$), and let us assume the sample has been normalized in such a way that $(\sampleinput_i, \sampleoutput_i) \in [-1,1]^{\inputdimension+\outputdimension}$. Furthermore, let $\parameternumber$ be the dimension of $\parameterset$ (the number of parameters in the neural network). With this notation, we define the following.

\begin{definition}
Consider the ERM problem \eqref{eq:ERMproblem} with parameters $\samplesize, \parameterset, \loss, f$ --- sample size, parameter space, loss function, network architecture, respectively. For a function $g$, let $\mathcal{K}_\infty(g)$ be the Lipschitz constant of $g$ using the infinity norm. We define the \emph{Architecture Lipschitz Constant} $\mathcal{K}(\samplesize,\parameterset,\loss, f)$ as
\begin{equation}
    \label{eq:Lipschitz-Arch}
    \mathcal{K}(\samplesize,\parameterset,\loss, f) \doteq \mathcal{K}_\infty(\loss(f(\cdot , \cdot), \cdot ))
\end{equation}
over the domain $[-1,1]^\inputdimension  \times \parameterset \times [-1,1]^\outputdimension$.
\end{definition}

Using this definition and the boundedness of parameters, a straightforward approximate training algorithm can be devised whose running time is linear in $D$. Simply do a grid search in the parameters' space, and evaluate all data points in each possible parameter. It is not hard to see that, to achieve $\epsilon$-optimality, such an algorithm would run in time which is linear in $D$ and exponential in $ \mathcal{K}(\samplesize,\parameterset,\loss, f) / \epsilon$.
What was proved by \cite{TrainingLP} is that one can take a step further and represent multiple training problems at the same time.

\begin{theorem}[\cite{TrainingLP}]\label{theorem:bienstocketal}
  Consider the ERM problem \eqref{eq:ERMproblem} with parameters $\samplesize, \parameterset, \loss, f$, and let $\mathcal{K}:= \mathcal{K}(\samplesize,\parameterset,\loss, f)$ be the corresponding network architecture. Consider $\epsilon > 0$ arbitrary. There exists a polytope $P_{\epsilon}$ of size\footnote{Here, the size of the polytope is the number of variables and constraints describing it.}
  \(O(\samplesize \left(2\mathcal{K}/\epsilon\right)^{\inputdimension+\outputdimension+\parameternumber}) \)
  with the following properties:
  \begin{enumerate}
      \item $P_{\epsilon}$ can be constructed in time $O(\left(2\mathcal{K}/\epsilon\right)^{\inputdimension+\outputdimension+\parameternumber} \samplesize)$ plus the time required for $O(\left(2\mathcal{K}/\epsilon\right)^{\inputdimension+\outputdimension+\parameternumber})$ evaluations of the loss function $\loss$ and $f$.
      \item For \emph{any} sample $(\tilde{X}, \tilde{Y}) = (\sampleinput_i, \sampleoutput_i)_{i=1}^\samplesize$, $(\sampleinput_i,\sampleoutput_i) \in [-1,1]^{\inputdimension+\outputdimension}$, there is a face $\mathcal{F}_{\tilde{X},\tilde{Y}}$ of $P_{\epsilon}$ such that optimizing a linear function over $\mathcal{F}_{\tilde{X},\tilde{Y}}$ yields an $\epsilon$-approximation to the ERM problem \eqref{eq:ERMproblem}.
      \item  The face $\mathcal{F}_{\tilde{X},\tilde{Y}}$ arises by simply substituting-in actual data for the data-variables $x,y$, which is used to fixed variables in the description of $P_{\epsilon}$.
  \end{enumerate}
\end{theorem}

This result is very abstract in nature but possesses some interesting features. Firstly, it encodes (approximately) \emph{every} possible training problem arising from data in $[-1,1]^{\inputdimension+\outputdimension}$ using a benign dependency on $\samplesize$: the polytope size depends only linearly on $\samplesize$, while a discretized enumeration of all the possible samples of size $\samplesize$ would be exponential in $\samplesize$. Secondly, every possible ERM problem appears in a \emph{face} of the polytope; this suggests a strong geometric structure across different ERM problems. Lastly, this result is applicable to a wide variety of network architectures; in order to obtain an architecture-specific result, it suffices to compute the corresponding value of $\mathcal{K}$ and plug it in. Regarding this last point, the authors computed the constant $\mathcal{K}$ for various well-known architectures and obtained the results of Table \ref{tab:results}.

\begin{table}
  \caption{Summary of polyhedral encoding sizes for various architectures. DNN refers to a fully-connected Deep Neural Network, CNN to a Convolutional Neural Network, and ResNet to a Residual Network. $G$ is the graph defining the Network, $\Delta$ is the maximum in-degree in $G$, $\hiddenlayers$ is the number of hidden layers, and $\maxlayerwidth$ is the maximum width of a layer.}
  \label{tab:results}
\vskip 0.15in
\begin{adjustbox}{max width=\textwidth}
  \begin{tabular}[h]{llll}
\hline
    Type  & Loss & Size of polytope & Notes  \\
\hline \hline
DNN & Absolute/Quadratic/Hinge & $O\big(\big( \outputdimension \maxlayerwidth^{O(\hiddenlayers^2)} /\epsilon \big)^{\inputdimension+\outputdimension+\parameternumber} \samplesize \big)$ & $\parameternumber\in O(|E({G})|)$ \\
DNN & Cross Entropy w/ Soft-Max & $O\big(\big( \outputdimension \log (\outputdimension) \maxlayerwidth^{O(k^2)} /\epsilon\big)^{\inputdimension+\outputdimension+\parameternumber} \samplesize\big)$ & $\parameternumber\in O(|E({G})|)$ \\
CNN & Absolute/Quadratic/Hinge & $O\big(\big( \outputdimension \maxlayerwidth^{O(\hiddenlayers^2)} /\epsilon\big)^{\inputdimension+\outputdimension+\parameternumber} \samplesize\big)$ & $\parameternumber \ll |E({G})|$ \\
ResNet & Absolute/Quadratic/Hinge & $O\big(\big( \outputdimension \Delta^{O(\hiddenlayers^2)} /\epsilon\big)^{\inputdimension+\outputdimension+\parameternumber} \samplesize\big)$  \\
ResNet & Cross Entropy w/ Soft-Max & $O\big(\big( \outputdimension \log(\outputdimension) \Delta^{O(\hiddenlayers^2)} /\epsilon\big)^{\inputdimension+\outputdimension+\parameternumber} \samplesize\big)$  \\
\hline
  \end{tabular}
  \end{adjustbox}
\vskip -0.1in
\end{table}

The proof of this result relies on a graph-theoretical concept called \emph{treewidth}. This parameter is used for measuring structured sparsity, and in \cite{bienstock2018lp} it was proved that any optimization problem admits an approximate polyhedral reformulation whose size is exponential only in the treewidth parameter. On a high level, the neural network result is obtained by noting that \eqref{eq:ERMproblem} connects different sample points only through a sum; therefore, the following reformulation of the optimization problem can be considered, which decouples the different data points:
\begin{align}
\label{ERMepigraph}
\begin{split}
    \min_{\weights \in \parameterset,\mL} \quad & \frac{1}{\samplesize} \sum_{d = 1}^{\samplesize} \mL_d   \\
    \mbox{s.t.} \quad &\mL_d \, = \, \loss(f(\sampleinput_d, \weights), \sampleoutput_d) \quad \forall\, d\in [\samplesize] 
    \end{split}
  \end{align}
This reformulation does not seem useful at first; however, it has a \emph{treewidth} that does not depend on $\samplesize$, even if the data points are considered variables. From this point, the authors are able to obtain the polytope whose size does not depend exponentially on $\samplesize$, and which is capable of encoding all possible ERM problems. The face structure the polytope has is more involved, and we refer the reader to \cite{TrainingLP} for these details.

It is worth mentioning that the polytope size provided by \cite{TrainingLP} in the setting of \cite{arora2018understanding} is 
\begin{equation}\label{eq:oursize}
    O( ( 2\mathcal{K}_\infty(\loss) \layerwidth_1^{O(1)} / \epsilon )^{(\inputdimension+1)(\layerwidth_1)} \samplesize )
\end{equation}
where $\mathcal{K_\infty}(\loss)$ is the Lipschitz constant of the loss function with respect to the infinity norm over a specific domain. These two results are not completely comparable, but they give a good idea of how good the size of the polytope constructed in \cite{TrainingLP} is. The dependency on $\samplesize$ is better in the polytope size, the polytope encodes multiple training problems, and the result is more general (it applies to almost any architecture); however, the polytope only gives an approximation, and its construction requires boundedness.

\subsection{Backpropagation through MILP}

In the work by \cite{goebbelstraining2021}, a novel use of Mixed-Integer Linear Programming is proposed in training ReLU networks: to serve as an alternative to SGD. This new algorithm works as backpropagation, as it updates the weights of the neural network iteratively, starting from the last layer. The key difference is that each update in a layer amounts to solving a MILP.

Let us focus only on one hidden layer at a time (of width $\layerwidth$), so we can assume we have an architecture as in Figure \ref{fig:HyperplaneToNN}. Furthermore, we assume we have some target output vectors $\{\mT_d\}_{d=1}^\samplesize$ (when processing the last hidden layer in the backpropagation, this corresponds to $\{\sampleoutput_d\}_{d=1}^\samplesize$) and some layer input $\{\mI_d\}_{d=1}^\samplesize$ (when processing the last hidden layer, this corresponds to evaluating the neural network on $\{\sampleinput_d\}_{d=1}^\samplesize$ up to the second-to-last hidden layer). The algorithm proposed by \cite{goebbelstraining2021} solves the following optimization problem to update the weights $\weights$ and biases $\biases$ of the given layer:
\begin{subequations}
\label{eq:backpropmilp}
\begin{align}
    \min_{\weights,\hat{\vh},\biases,\vh,\vz} \quad & \sum_{d=1}^\samplesize \sum_{j=1}^\layerwidth |\mT_{d,j} - \vh_{d,j}| \\
    \text{s.t.}\quad & \forall d\in [D],\, j\in[n],\, \nonumber\\ &\hat{\vh}_{d,j} = (\weights \mI_d)_j + \biases_j   \\
    &\hat{\vh}_{d,j} \leq M\vz_{d,j} \\ 
    &\hat{\vh}_{d,j} \geq -M(1-\vz_{d,j})  \\
    &|\hat{\vh}_{d,j}-\vh_{d,j}| \leq M(1-\vz_{d,j}) \\ 
    &\vh_{d,j} \leq M\vz_{d,j} \\
    & \vh_{d,j} \geq 0 \\
    & \vz_{d,j} \in \{0,1\}.
\end{align}
\end{subequations}
Here, $M$ is a large constant that is assumed to bound the input to any neuron. Note that problem \eqref{eq:backpropmilp} can easily be linearized. 
This optimization problem finds the weights ($\weights$) and biases ($\biases$) that minimize the difference between the ``real'' output of the network for each sample ($\vh_d$) and the target output ($\mT_d$). The auxiliary variables $\hat{\vh}_{d,j}$ represent the input to the each neuron ---so $\vh_{d,j} = \sigma(\hat{\vh}_{d,j})$--- and $\vz_{d,j}$ indicates if the $j$-th neuron is activated on input $\mI_{d}$.

When processing intermediate layers, the definition $\mI_d$ can easily be adapted from what we mentioned above. However, the story is different for the case of $\mT_d$. When processing the last layer, as previously mentioned, $\mT_d$ simply corresponds to $\sampleoutput_d$. For intermediate layers, to define $\mT_d$, the author proposes to use a similar optimization problem to \eqref{eq:backpropmilp}, but leaving $\weights$ and $\biases$ fixed and having $\mI_d$ as variables; this defines ``optimal inputs'' of a layer. These optimal inputs are then used as target outputs $\mT_d$ when processing the preceding layer, and thus the algorithm is iterated. For details, see \cite{goebbelstraining2021}.

The computational results in that paper show that a similar level of accuracy to that of gradient descent can be achieved. However, the use of potentially expensive MILPs impairs the applicability of this approach to large networks. Nonetheless, it shows an interesting new avenue for training whose running times may be improved in future implementations.

\subsection{Training binarized neural networks using MILP}

As mentioned before, the training problem of a DNN is an unrestricted non-convex optimization
problem, which is typically continuous as the weights and biases are frequently allowed to have any real value.
Nonetheless, if the weights and biases are required to be integer-valued, the training
problem becomes a discrete optimization problem, for which gradient-descent-based methods may find some difficulties in their applicability.

In this context, \cite{icarte2019training} proposed a MILP formulation for the training problem of binarized neural networks (BNNs): these are neural networks where the weights and biases are restricted to be in $\{-1,0,1\}$ and where the activations are LTU (i.e., sign functions). Later on, \cite{thorbjarnarson2020training,thorbjarnarson2023optimal} used a similar technique to allow more general integer-valued weights. We review the core feature in these formulations that yield a \emph{linear} formulation of the training problem. 

Let us focus on an intermediate layer $i$ with width $\layerwidth$, and let us omit biases to simplify the discussion. Using a DNN's layer-wise architecture, one usually aims \change{to describe}:
\begin{subequations} \label{eq:basicsystemtraining}
\begin{align}
    \hat{\vh}^i_{d,j} & = (\weights^{i} \vh^{i-1}_d)_j  && d=1,\ldots,D,\, j=1,\ldots,n \label{eq:basicsysinput}\\
    \vh^{i}_{d,j} & = \sigma(\hat{\vh}^i_{d,j}) && d=1,\ldots,D,\, j=1,\ldots,n. \label{eq:basicsysactivation}
\end{align}
\end{subequations}
We remind the reader that $D$ is the cardinality of the training set.  Additionally, for each data point indexed by $d$ and each layer $i$, each variable $\vh_{d}^{i}$ is the output vector of all the neurons of the layer, and each variable $\hat{\vh}_{d,j}^i$ is the input of neuron $j$.
Besides the difficulty posed by the activation function, one important issue with system \eqref{eq:basicsystemtraining} is the non-linearity of the products between the $\weights$ and $\vh$ variables. Nonetheless, this issue disappears when each entry of $\weights$ and $\vh$ is bounded and integer, as in the case of BNNs. 

Let us begin with reformulating \eqref{eq:basicsysactivation}. We can introduce auxiliary variables $\vu_{d,j}^i\in \{0,1\}$ that will indicate if the neuron is active. We also introduce a tolerance $\varepsilon > 0$ to determine the activity of a neuron. Using this, we can (approximately) reformulate \eqref{eq:basicsysactivation} \emph{linearly} using big-M constraints:
\begin{subequations}\label{eq:activationreform}
\begin{align}
    \vh^{i}_{d,j} & = 2\vu_{d,j}^i -1     && d=1,\ldots,D,\, j=1,\ldots,n \label{eq:oandu}\\
    \hat{\vh}^i_{d,j} & \geq -M (1-\vu_{d,j}^i) && d=1,\ldots,D,\, j=1,\ldots,n \\
     \hat{\vh}^i_{d,j} & \leq -\varepsilon + M \vu_{d,j}^i && d=1,\ldots,D,\, j=1,\ldots,n 
\end{align}
\end{subequations}
where $M$ is a large constant. 
As for \eqref{eq:basicsysinput}, note that 
\[ \hat{\vh}^i_{d,j} = \sum_{k=1}\weights^{i}_{j,k} \vh^{i-1}_{d,k}. \]
Therefore, using \eqref{eq:oandu}, we see that it suffices to describe each product $\weights_{j,k}^{i}\vu_{d,k}^{i-1}$ linearly. We can introduce new variables $\vz_{j,k,d}^{i}$ and note that
\[\vz_{j,k,d}^{i-1}= \weights_{j,k}^{i}\vu_{d,k}^{i-1}\]
if and only if the three variables satisfy
\begin{align*}
    | \vz_{j,k,d}^{i-1} | & \leq \vu_{d,k}^{i-1} \\
    |\vz_{j,k,d}^{i-1} - \weights_{j,k}^{i}| &  \leq 1- \vu_{d,k}^{i-1} \\
    \vu_{d,k}^{i-1} & \in \{0,1\}.
\end{align*}
This last system can be easily converted to a linear system, and thus the training problem in this setting can be cast as a mixed-integer linear optimization problem.

Other works have also relied on similar formulations to train neural networks. \cite{icarte2019training} introduce different objective functions that can be used along
with the linear system to produce a MILP that can train BNNs. They also introduce a Constraint-Programming-based model and a hybrid model and then compare all of them computationally. \cite{thorbjarnarson2020training} introduce more MILP-based training models that leverage piecewise linear approximations of well-known non-linear loss functions and that can handle integer weights beyond $\{-1,0,1\}$. A similar setting is studied by \cite{sildir2022mixed}, where piecewise linear approximations of non-linear activations are used, and integer weights are exploited to formulate the training problem as a MILP. Finally,  \cite{bernardelli2022bemi} rely on a multi-objective MIP model for training BNNs; from here, they create a BNN ensemble to produce robust classifiers. 

From these articles, we can conclude that the MILP-based approach to training their neural networks can result in high-quality neural networks, especially in terms of generalization. However, many of these MILP-based methods currently do not scale well, as opposed to gradient-descent-based methods. We believe that, even though there are some theoretical limitations to the efficiency of MILP-based methods, there is considerable practical improvement potential in using them in neural network training.

\section{Conclusion}

The rapid advancement of neural networks and their ubiquity has given rise to numerous new challenges and opportunities in deep learning: we need to design them in more reliable ways, to better understand their limits, and to test their robustness, among other challenges.
While, traditionally, continuous optimization has been the predominant technology used in the optimization tasks in deep learning, some of these new challenges have made discrete optimization tools gain a remarkable importance.

In this survey, we have reviewed multiple areas where polyhedral theory and linear optimization have played a critical role. For example, in understanding the expressiveness of neural networks, in optimizing trained neural networks (e.g. for verification purposes), and even in designing new training algorithms. 
We hope this survey can provide perspective in a rapidly-changing field, and motivate further developments in both deep learning and discrete optimization. There is still much to be explored in the intersection of these fields.

\paragraph{Acknowledgments} We thank Christian Tjandraatmadja and Toon Tran for early feedback on the manuscript and asking questions that helped shaping it. 
We also thank the anonymous reviewers for their corrections and suggestions.

Thiago Serra was  supported by
the National Science Foundation (NSF) award IIS 2104583. 
Calvin Tsay was supported by the Engineering \& Physical Sciences Research Council (EPSRC) grant EP/T001577/1. 

\bibliography{references}

\end{document}